\documentclass{jfm}
\usepackage{graphicx}
\usepackage{epstopdf, epsfig}

\usepackage{amsmath}
\usepackage{xcolor}

\definecolor{blue}{rgb}{0,0,0}
\definecolor{green}{rgb}{0,0,0}
\definecolor{cyan}{rgb}{0,0,0}
\definecolor{red}{rgb}{0,0,0}

\definecolor{brown}{rgb}{0,0,0}
\definecolor{pink}{rgb}{0, 0, 0}

\usepackage{cleveref}
\usepackage{bm}

\usepackage{tabularx}

\DeclareMathOperator*{\argmin}{arg\,min}

\title{Sparsity-promoting algorithms for the discovery of informative Koopman-invariant subspaces}

\shortauthor{S. Pan, N. Arnold-Medabalimi, and K. Duraisamy}

\shorttitle{Extraction of informative Koopman-invariant subspaces}

\author{Shaowu Pan\aff{1}\corresp{\email{shawnpan@umich.edu}},
        Nicholas Arnold-Medabalimi\aff{1} and
        Karthik Duraisamy\aff{1}
}

\affiliation{\aff{1}Department of Aerospace Engineering, University of Michigan, Francois-Xavier Bagnoud Building, 1320 Beal Ave, Ann Arbor, MI, USA}

\begin{document}

\maketitle

\begin{abstract}

Koopman decomposition is a non-linear generalization of eigen-decomposition, and is being increasingly utilized in the analysis of spatio-temporal dynamics. Well-known techniques such as the dynamic mode decomposition (DMD) and its linear variants provide approximations to the Koopman operator, and have been applied extensively in many fluid dynamic problems. Despite being endowed with a richer dictionary of nonlinear observables, nonlinear variants of the DMD, such as extended/kernel dynamic mode decomposition (EDMD/KDMD) are seldom applied to large-scale problems primarily due to the difficulty of discerning the Koopman invariant subspace from thousands of resulting Koopman eigenmodes. To address this issue, we propose a framework based on multi-task feature learning to extract the most informative Koopman invariant subspace by removing redundant and spurious Koopman triplets. {\color{brown} In particular, we develop a pruning procedure that penalizes departure from linear evolution.} These algorithms can be viewed as sparsity promoting extensions of EDMD/KDMD. Further, we extend KDMD to a continuous-time setting and show a relationship between the present algorithm, sparsity-promoting DMD, and an empirical criterion from the viewpoint of non-convex optimization. The effectiveness of our algorithm is demonstrated on examples ranging from simple dynamical systems to two-dimensional cylinder wake flows at different Reynolds numbers and a three-dimensional turbulent ship air-wake flow. The latter two problems are designed such that very strong nonlinear transients are present, thus requiring an accurate approximation of the Koopman operator. {\color{brown}Underlying physical mechanisms are analyzed, with an emphasis on characterizing transient dynamics. The results are compared to existing theoretical expositions and numerical approximations.}

\end{abstract}

\begin{keywords}
Koopman decomposition; Sparsity promoting techniques; Kernel methods; feature learning; Data-driven analysis; Operator theory.
\end{keywords}

\section{Introduction}
\label{intro}

\textcolor{red}{Complex unsteady flow phenomena such as turbulence~\citep{pope2001turbulent}, flow instability~\citep{drazin2004hydrodynamic,lietz2017plasma}, and fluid-structure interactions~\citep{dowell2001modeling} are prevalent in many physical systems.
To analyze and understand such phenomena, it is useful to extract coherent modes associated with important dynamical mechanisms {\em and} track their evolution.} Koopman operator theory~\citep{budivsic2012applied} offers an elegant framework to  reduce spatio-temporal fields associated with non-linear dynamics as a linear combination of time evolving modes \textcolor{red}{ordered by frequency and growth rates}~\citep{rowley2009spectral}. 

Consider a general continuous nonlinear dynamical system $\dot{\mathbf{x}} = \mathbf{F}(\mathbf{x})$, where the system state $\mathbf{x}$ {\color{red} evolves on a manifold $\mathcal{M} \subset \mathbb{R}^N$.} $\mathbf{F}: \mathcal{M} \mapsto T\mathcal{M}$ is a vector-valued smooth function and $T\mathcal{M}$ is the tangent bundle, i.e., $\forall p \in \mathcal{M}, \mathbf{F}(p) \in T_p \mathcal{M} $. The aforementioned task of decomposition is equivalent to finding a set of $L$ observables associated with the system, $\{g_i\}_{i=1}^{L}$ $(g_i: \mathcal{M} \mapsto \mathbb{C})$, that evolve linearly in time while spanning the system state $\mathbf{x}$. The significance of the above statement is that it represents  a {\em global} linearization of the nonlinear dynamical system~\citep{budivsic2012applied,brunton2016koopman}. 

Formally, this idea can be traced back to Koopman operator theory for Hamiltonian dynamical systems introduced by \citet{koopman1931hamiltonian} to study the evolution of observables $g$ on $L^2(\mathcal{M})$, \textcolor{pink}{which is a vector space of square integrable functions defined on the manifold $\mathcal{M}$}. \textcolor{pink}{More generally, } for \textcolor{pink}{a certain vector space of observable functions $\mathcal{F}$ defined on the manifold $\mathcal{M}$} and $t>0$, the \textcolor{pink}{Koopman semigroup (parameterized by $t$)\footnote{\textcolor{pink}{Note that for a discrete dynamical system resulting from the discretization (with time step $\Delta t$) of a continuous dynamical system, the corresponding discrete-time Koopman operator is an element of this semigroup: $\mathcal{K}_{\Delta t}$.}} is the set $\{\mathcal{K}_t\}_{t\in\mathbb{R}^+}$:} $\mathcal{K}_t: \mathcal{F} \mapsto \mathcal{F}$ that governs the dynamics of observables in the form $\mathcal{K}_t g(\mathbf{x}(0)) \triangleq \textcolor{red}{g}(S(t, \mathbf{x}_0)) = g (\mathbf{x}(t))$ where $S(t,\cdot)$ is the flow map of the dynamical system~\citep{mezic2013analysis}. \textcolor{pink}{$\mathcal{K} = \lim_{t\rightarrow 0} (\mathcal{K}_t f - f) / t$ is referred to as the continuous-time Koopman operator, i.e., Koopman generator.} While the Koopman operator is linear over the space of observables, $\mathcal{F}$ is most often infinite-dimensional, e.g., $L^2(\mathcal{M})$, which makes the approximation of the Koopman operator a difficult problem. Throughout this work, we assume $\mathcal{F} = L^2(\mathcal{M})$. Readers interested in a  detailed discussion of the choice of $\mathcal{F}$ are referred to~\citet{bruce2019koopman}. 

\begin{figure}
\centering
\includegraphics[width=\textwidth]{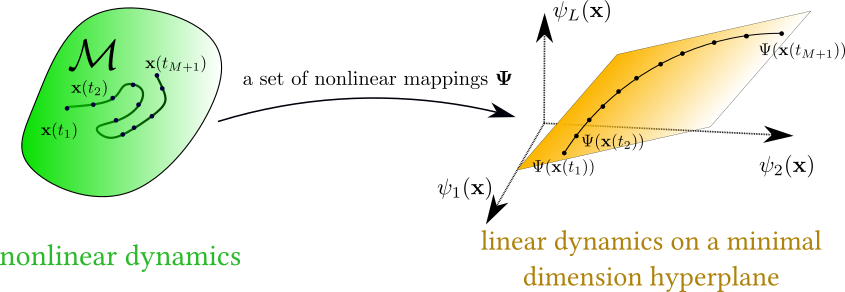}
\caption{\textcolor{red}{Schematic of transformation of a non-linear system to a linear $L$-dimensional space and the extraction of a minimal Koopman-invariant subspace.}}
\label{fig:mkisp}
\end{figure}

A special subspace of $\mathcal{F}$, referred to as a \emph{minimal Koopman-invariant subspace} $\mathcal G$, has the following property: for any $\phi \in \mathcal{G}$, for any $t \in \mathbb{R}^{+}$, $\mathcal{K}_t \phi \in \mathcal{G}$ \textcolor{pink}{and $\mathbf{x} \in \mathcal{G}$}. \textcolor{pink}{ Existing techniques such as the extended dynamic mode decomposition~\citep{williams2015data} are capable of approximating Koopman operators, but typically yield high dimensional spaces ($L$-dimensional space in \cref{fig:mkisp}). In this work, we are interested in accurately approximating such a subspace - but with the minimal possible dimension and Koopman invariance  -  using learning techniques, as illustrated in \cref{fig:mkisp}.  This can yield useful coordinates for a multitude of applications including modal analysis and control~\citep{arbabi1804data}.  }

In general, physical systems governed by PDEs, e.g., fluid dynamics, are infinite-dimensional. From a numerical viewpoint, the number of degrees of freedom can be related to the spatial discretization  (for example, the number of grid points). Although a finite-dimensional manifold can be extracted  \citep{holmes2012turbulence}, e.g., $O(10)$-$O(10^3)$ dimensions via Proper Orthogonal Decomposition (POD), finding a Koopman-invariant subspace on such manifolds is still challenging.
Currently, the most popular method to approximate the Koopman operator  is the dynamic mode decomposition (DMD) \citep{schmid2010dynamic,rowley2017model} mainly for two reasons. First, it is  straightforward and computationally efficient compared to nonlinear counterparts such as Extended DMD (EDMD) \citep{williams2015data} and Kernel DMD (KDMD)~\citep{williams2014kernel}. Second, the essence of DMD is to decompose a spatio-temporal field into several temporally growing/decaying travelling/stationary harmonic waves, which are prevalent in fluid mechanics.   
However, the accuracy of DMD is limited by the assumption that the Koopman-invariant subspace lies in the space spanned by snapshots of the state $\mathbf{x}$. Thus, DMD  is  used to mainly identify and visualize  coherent structures. Indeed, DMD can be interpreted as a $L^2$ projection of the action of the Koopman operator on the linear space spanned by snapshots of the system state~\citep{korda2018convergence}. 

To overcome the above limitations, it is natural to augment the observable space with either the history of the state~\citep{arbabi2017ergodic,brunton2017chaos,kamb2018time,le2017higher} or nonlinear observables of the state~\citep{williams2014kernel,williams2015data}. \textcolor{cyan}{Time-delay embedding can be very useful in reduced order modeling of systems   for which sparse measurements can be easily obtained, assuming  the inputs and outputs are not high dimensional~\citep{korda2018linear}.} Although time-delay embedding is simple to implement and has strong connections to Takens' embedding~\citep{kamb2018time,pan2019structure}, the main practical issue arises in reduced order modeling of high fidelity simulations in a predictive setting due to the requirement of a large number of snapshots of the full order model.
Further, if one is only interested in the post-transient dynamics of the system state  \emph{on} an attractor, linear observables with time delays are sufficient to extract an informative Koopman-invariant subspace \citep{arbabi2017ergodic,pan2019structure,brunton2017chaos,arbabi2017study,mezic2005spectral,rojsel2017koopman}. However, if one is interested in the \textcolor{cyan}{strongly nonlinear} transient dynamics leading to an attractor or reduced order modeling \textcolor{cyan}{for high fidelity numerical simulations~\citep{xu2019multi,xu2020reduced,huang2018challenges,parish2020adjoint,carlberg2013gnat}}, time-delay embedding may become less appropriate as several delay snapshots \textcolor{cyan}{of the full order model} are required to initialize the model. In that case, nonlinear observables may be more appropriate. 

Driven by the interest in modal analysis and control of  transient flow phenomena,
we consider augmentation of the observable space with nonlinear functions of the state, e.g.,  EDMD~\citep{williams2015data}/KDMD~\citep{williams2014kernel}. 
Although it has been reported that KDMD allows for a set of more interpretable Koopman eigenvalues~\citep{williams2014kernel} and better accuracy~\citep{rojsel2017koopman}, issues such as mode selection, spurious modes~\citep{kaiser2017data,zhang2017evaluating}, and choice of dictionary/kernel in EDMD/KDMD  remain open. In fact, the choice of kernel type and hyperparameter can significantly affect the resulting eigenmodes,  distribution of eigenvalues~\citep{kutz2016dynamic}, and the accuracy of predictions~\citep{zhang2017evaluating}.  

Searching for an accurate and informative Koopman-invariant subspace has long been a  pursuit in the DMD community. \citet{rowley2009spectral} and \citet{schmid2012decomposition} considered selecting dominant DMD modes in  order of their amplitudes. However, following such a criterion~\citep{tu2013dynamic,kou2017improved}, may result in the selection of \textcolor{brown}{irrelevant} modes that may have large amplitudes, but  decay rapidly.  As a result, \citet{tu2013dynamic} considered weighting the loss term by the magnitude of eigenvalues to penalize the retention of fast decaying modes.
Sparsity-promoting DMD (\textcolor{cyan}{referred as ``spDMD'' throughout the paper}) developed by  \citet{jovanovic2014sparsity} recasts mode selection in DMD as an optimization problem with a $\ell_1$ penalty. With a preference of stable modes over fast decaying ones, \citet{tissot2014model} proposed a simpler criterion based on time-averaged-eigenvalue-weighted amplitude. This was followed by \citet{kou2017improved} who used  a similar criterion but computed the ``energy" of each mode, yielding  similar performance to spDMD at a  lower computational cost. Based on the orthonormal property of pseudo-inverse, \citet{hua2017high} proposed an ordering of Koopman modes by defining a new ``energy". Compared with previous empirical criteria, the ``energy" for each mode involves a pseudo-inverse which combines the influence from all eigenmodes, and therefore  the contribution from each mode cannot be isolated. 
Instead of selecting modes from a ``reconstruction'' perspective,  \citet{zhang2017evaluating} studied the issue of spurious modes by evaluating the deviation of the identified eigenfunctions from linear evolution in an a priori sense. 
Further, Optimized DMD~\citep{chen2012variants,askham2018variable} combines DMD with mode selection simultaneously, which
is the forerunner of recently proposed neural network-based models for Koopman eigenfunctions in spirit~\citep{takeishi2017learning,lusch2018deep,otto2019linearly,pan2019physics,li2017extended,yeung2019learning}. Regardless of the above issues related to non-convex optimization~\citep{dawson2016characterizing}, extension of optimized DMD to EDMD/KDMD is not straightforward. Further, neural network-based models require large amounts of data, are prone to local minima,  and lack interpetability.
%
There have been a few attempts towards  mode selection in EDMD/KDMD. \citet{brunton2016koopman} present an iterative method that augments the dictionary of EDMD until a convergence criterion is reached for the subspace. This is effectively  a recursive implementation of EDMD.  Recently, \citet{haseli2019efficient} showed that given a sufficient amount of data, if there is any accurate Koopman eigenfunction spanned by the dictionary, it must correspond to one of the obtained eigenvectors. Moreover, they proposed the idea of mode selection by checking if the reciprocal of identified eigenvalue also appears when the temporal sequence of data is reversed, which is similar to the idea of comparing eigenvalues on the complex plane from different trajectories, as proposed by \citet{hua2017high}. In contrast to the  ``bottom-up" method of \citet{brunton2016koopman} with subspace augmentation, \citet{hua2017high} proposed a ``top-down" subspace subtraction method relying on iteratively projecting the features onto the null space. 
A similar idea can be traced back to \citet{kaiser2017data} who propose a search for the sparsest vector in the null space. 

\textcolor{red}{As the main contribution of this work,} we propose a novel EDMD/KDMD framework equipped with the following strategy to extract an accurate yet minimal Koopman-invariant subspace: 
\begin{enumerate}
\item We first evaluate the normalized maximal deviation of the evolution of each eigenfunction from linear evolution in a posteriori fashion.
\item Using the above criteria,  we select a user-defined number of accurate EDMD/KDMD modes;
\item Among the accurate EDMD/KDMD modes obtained above, informative modes are selected using multi-task feature learning~\citep{argyriou2008convex,bach2012optimization}. 
\end{enumerate}

To the best of our knowledge, this is the first \textcolor{red}{model-agnostic} attempt to address sparse identification of \textcolor{red}{an accurate minimal} Koopman-invariant subspace. \textcolor{red}{Our contribution also extends the classical sparsity-promoting DMD~\citep{jovanovic2014sparsity} in a more general setting that includes EDMD/KDMD.} \textcolor{brown}{ The applications are focused on strongly transient flows, and new classes of stable Koopman modes are identified.}
    

The organization of the paper is as follows: In \Cref{sec:review}, we  provide a review of EDMD/KDMD in  discrete-time and present corresponding extensions to continuous-time. Following this, we discuss  current challenges in  standard EDMD/KDMD. In \Cref{sec:select}, we propose the  framework of sparse identification of  informative Koopman-invariant subspaces for EDMD/KDMD  with hyperparameter selection and provide implementation details. \textcolor{brown}{A novel optimization procedure that penalizes departure from linear evolution is detailed in \Cref{sec:pruning,sec:multi-task}}. Differences and connections between spDMD, Kou's emprical criterion, and our proposed framework is shown from the viewpoint of optimization. In \Cref{sec:result}, numerical verification and modal analysis on examples from a fixed point attractor to a transient cylinder wake flow, and a transient ship-airwake are provided. In \Cref{sec:conclusion}, conclusions are drawn.

\section{Review of EDMD/KDMD}
\label{sec:review}

In this section, we provide a summary of our framework to extract the Koopman operator using EDMD/KDMD in both discrete-time and continuous-time.
Although the original EDMD is seldom used for data-driven Koopman analysis in fluid flows due to its poor scaling, it has recently been reported that a scalable version of EDMD can be applied to high-dimensional systems~\citep{degennaro2019scalable}. Since our algorithm to extract the Koopman-invariant subspace is agnostic to the type of technique to compute the original set of Koopman modes, we include both EDMD and KDMD for completeness.

\subsection{Extended Dynamic Mode Decomposition}

\subsubsection{Discrete-time Formulation}

For simplicity, consider $M$ sequential snapshots  sampled uniformly in time with $\Delta t$ on a trajectory, $\{\mathbf{x}_i\}_{i=1}^{M}$. The EDMD algorithm~\citep{williams2015data} assumes a dictionary of $L$ linearly independent functions $i=1,\ldots,L$, $\psi_i(\mathbf{x}): \mathcal{M} \mapsto \mathbb{C}$ that approximately spans a Koopman-invariant subspace $\mathcal{F}_L$ 
\begin{equation}
\label{eq:F_L}
\mathcal{F}_L = \textrm{span} \{ \psi_1(\mathbf{x}), \ldots, \psi_L(\mathbf{x}) \}.
\end{equation}

Thus we can write for any 
$g \in \mathcal{F}_{L}$ as $g(\mathbf{x}) = \mathbf{\Psi}(\mathbf{x}) \mathbf{a}$ with $\mathbf{a} \in \mathbb{C}^{L}$, where the feature vector $\mathbf{\Psi}(\mathbf{x})$ is 
\begin{equation}
    \label{eq:feature_vector}
    \mathbf{\Psi}(\mathbf{x}) = 
    \begin{bmatrix}
    \psi_1(\mathbf{x}) & \ldots & \psi_L(\mathbf{x})    
    \end{bmatrix}.
\end{equation}
Consider a set of $L^{'}$ observables as $\{g_l\}_{l=1}^{L^{'}} = \{\mathbf{\Psi} \mathbf{a}_l\}_{l=1}^{L^{'}}$,  where \textcolor{pink}{$\mathbf{a}_l \in \mathbb{C}^L$} is arbitrary. After the discrete-time Koopman operator $\mathcal{K}_{\Delta t}$ is applied on each $g_l$, given data $\{\mathbf{x}_i\}_{i=1}^{M}$, we have the following \textcolor{cyan}{for $l=1,\ldots,L^{\prime}$}, \textcolor{red}{and $i=1,\ldots,M-1$}, 
\begin{equation}
     \label{eq:discrete_koopman_linear_relation}
     \mathcal{K}_{\Delta t} g_l(\mathbf{x}_i) = g_l(\mathbf{x}_{i+1})  = \mathbf{\Psi}(\mathbf{x}_{i+1}) \mathbf{a}_l = \mathbf{\Psi}(\mathbf{x}_{i}) \mathbf{K} \mathbf{a}_l + r_{i,l},
\end{equation} 
\textcolor{pink}{where $r_{i,l}$ is simply the residual for the $l$-th function on the $i$-th data point.}
The standard EDMD algorithm seeks \textcolor{red}{a constant matrix} $\mathbf{K} \in \mathbb{C}^{L \times L}$ \textcolor{red}{that governs the linear dynamics in the observable space} such that \textcolor{pink}{the sum of the square of the residuals $r_{i,l}$ from \cref{eq:discrete_koopman_linear_relation}} over all samples and functions, 
\begin{align}
    {\label{eq:cost_func_edmd}}
    J(\mathbf{K}, \{\mathbf{a}_l\}_{l=1}^{L^{'}}) = \sum_{l=1}^{L^{'}}\sum_{m=1}^{M-1} |( \mathbf{\Psi}(\mathbf{x}_{m+1}) - \mathbf{\Psi}(\mathbf{x}_m) \mathbf{K}  )\mathbf{a}_l|^2 = \lVert (\mathbf{\Psi}_{\mathbf{x}^{'}} - \mathbf{\Psi}_{\mathbf{x}}\mathbf{K})
     \mathbf{A}^{'} \rVert^2_F,
\end{align}
is minimized over $\{\mathbf{x}_i\}_{i=1}^{M+1}$.
In the above equation, 
\begin{equation}
    \mathbf{\Psi}_{\mathbf{x}} = 
    \begin{bmatrix}
        \psi_1(\mathbf{x}_1) & \ldots & \psi_L(\mathbf{x}_1) \\
        \vdots & \vdots & \vdots \\
        \psi_1(\mathbf{x}_{M-1}) & \ldots & \psi_L(\mathbf{x}_{M-1})  
    \end{bmatrix}, 
    \mathbf{\Psi}_{\mathbf{x}^{'}} = 
    \begin{bmatrix}
        \psi_1(\mathbf{x}_2) & \ldots & \psi_L(\mathbf{x}_2) \\
        \vdots & \vdots & \vdots \\
        \psi_1(\mathbf{x}_{M}) & \ldots & \psi_L(\mathbf{x}_{M})  
    \end{bmatrix}
\end{equation}
\begin{equation}
    \mathbf{A}^{'} = \begin{bmatrix}\mathbf{a}_1 & \ldots & \mathbf{a}_{L^{'}}\end{bmatrix}.
\end{equation}

Considering $\partial J/\partial \mathbf{K} = \mathbf{0}$, one obtains
$\mathbf{\Psi}_{\mathbf{x}}^{\textrm{H}} \mathbf{\Psi}_{\mathbf{x}^{'}}  \mathbf{A}^{'} \mathbf{A}^{'\textrm{H}} = \mathbf{\Psi}_{\mathbf{x}}^{\textrm{H}} \mathbf{\Psi}_{\mathbf{x}} \mathbf{K} \mathbf{A}^{'}  \mathbf{A}^{'\textrm{H}}$. Thus, the corresponding  minimizer $\mathbf{K}_{\textrm{opt}}$ is,
\begin{equation}
    \label{eq:discrete_EDMD_general}
    \mathbf{K}_{\textrm{opt}} = \mathbf{G}^{+} \mathbf{A} ( \mathbf{A}^{'} \mathbf{A}^{'\textrm{H}}) ( \mathbf{A}^{'} \mathbf{A}^{'\textrm{H}})^{+},
\end{equation}
where $+$ denotes the pseudoinverse and
\begin{align}
\label{eq:G_EDMD}
    \mathbf{G} &=  \sum_{m=1}^{M-1} \mathbf{\Psi}(\mathbf{x}_m)^{\textrm{H}} \mathbf{\Psi}(\mathbf{x}_m) = \mathbf{\Psi}_{\mathbf{x}}^{\textrm{H}} \mathbf{\Psi}_{\mathbf{x}} , \\
    \mathbf{A} &=  \sum_{m=1}^{M-1} \mathbf{\Psi}(\mathbf{x}_m)^{\textrm{H}} \mathbf{\Psi}(\mathbf{x}_{m+1}) =  \mathbf{\Psi}_{\mathbf{x}}^{\textrm{H}} \mathbf{\Psi}_{\mathbf{x}^{'}},
\end{align}
where $\textrm{H}$ denotes conjugate transpose. Note that when the set of observable fully span $\mathcal{F}_L$, i.e., $\mathbf{A}^{'}$ is full rank, $( \mathbf{A}^{'} \mathbf{A}^{'\textrm{H}}) ( \mathbf{A}^{'} \mathbf{A}^{'\textrm{H}})^{+}$ reduces to identity. Then we have the more familiar $\mathbf{K}_{\textrm{EDMD}}$ as,
\begin{equation}
    \label{eq:discrete_EDMD}
    \mathbf{K}_{\textrm{EDMD}} = \mathbf{G}^{+} \mathbf{A}.
\end{equation}
which is independent of \emph{the choice of $\{\mathbf{a}_l\}_{l=1}^{L^{'}}$}. 

Assuming that all eigenvalues of $\mathbf{K}_{\textrm{EDMD}}$ are simple\footnote{This is an immediate consequence for ergodic system~\citep{parry2004topics}.}, for $i=1,\ldots,L$, the corresponding Koopman eigenfunctions $\mathbf{\varphi}_{i}$ associated with Koopman eigenvalues $\lambda_i $ are 
\begin{equation}
   \label{eq:edmd_eigen_functions}
    \varphi_i(\mathbf{x}) = \mathbf{\Psi} (\mathbf{x}) \mathbf{v}_i, 
\end{equation}
where $\mathbf{K}_{\textrm{EDMD}} \mathbf{v}_i = \lambda_{i} \mathbf{v}_i$. Finally, the Koopman modes of a vector-valued $Q$ dimensional observable function $\mathbf{g}: \mathcal{M} \mapsto \mathbb{C}^Q$, 
are the vectors of weights $\{\mathbf{b}_i\}_{i=1}^{L}$ assoicated with the expansion in terms of eigenfunctions $\{\varphi_i\}_{i=1}^{L}$ as, 
\begin{equation}
    \mathbf{g}(\mathbf{x}) = \sum_{i=1}^{L} \varphi_i(\mathbf{x}) \mathbf{b}_i,
\end{equation}
where $\mathbf{b}_i$ is often numerically determined in practice by ordinary least squares,  
\begin{equation}
    \mathbf{B} = 
    \begin{bmatrix}
        \mathbf{b}_1 \\
        \vdots \\
        \mathbf{b}_L
    \end{bmatrix}
    = 
    \left(\mathbf{\Psi}_{\mathbf{x}} \mathbf{V}\right)^{+}
    \begin{bmatrix} 
        \mathbf{g}(\mathbf{x}_1) \\ \vdots \\ \mathbf{g}(\mathbf{x}_{M})
    \end{bmatrix},
\end{equation}
with $\mathbf{V} = \begin{bmatrix} \mathbf{v}_1 & \ldots & \mathbf{v}_L \end{bmatrix}$.

\subsubsection{Continuous-time Formulation}

Consider $M$ data snapshots of the dynamical system with state $\mathbf{x}$ sampled over $\mathcal{M}$ as $\{\mathbf{x}_i, \dot{\mathbf{x}}_{i}\}_{i=1}^{M}$ where $\dot{\mathbf{x}}_{i} = \mathbf{F}(\mathbf{x}_i)$.
Recall the generator of the semigroup of Koopman operators $\mathcal{K}: \mathcal{D}(\mathcal{K}) \mapsto \mathcal{F}, \mathcal{K} = \lim_{t\rightarrow 0} (\mathcal{K}_t f - f) / t$ where $\mathcal{D}(\mathcal{K})$  is the domain in which the aforementioned limit is well-defined and the closure of ${\mathcal{D}(\mathcal{K})}$ is $\mathcal{F}$. One can have the evolution of \textcolor{pink}{any observable $\mathbf{g} = \mathbf{\Psi} \mathbf{a} \in \mathcal{F}_L$} as,
\begin{equation}
 \mathcal{K} \mathbf{g} = \dot{\mathbf{g}}  = \mathbf{F} \cdot \nabla_{\mathbf{x}} \mathbf{\Psi} \mathbf{a} = \mathbf{\Psi} \mathbf{K} \mathbf{a} + r,
\end{equation}
where $r$ is  the residual.
Similarly, one can find a $\mathbf{K}$ that minimizes the sum of the square of residual $r$ minimized solution, 
\begin{equation}
\label{eq:continuous_koopman_linear_relation}
\mathbf{K}_{\textrm{EDMD}} =  \mathbf{G}^{+} \mathbf{A},
\end{equation}
where 
\begin{align}{\label{eq:edmd_G}}
\mathbf{G} &=  \sum_{m=1}^{M}\mathbf{\Psi}(\mathbf{x}_m)^{\textrm{H}} \mathbf{\Psi}(\mathbf{x}_m),\\
{\label{eq:edmd_A}}
\mathbf{A} &= \sum_{m=1}^{M}\mathbf{\Psi}(\mathbf{x}_m)^{\textrm{H}} (\mathbf{F} \cdot \nabla_{\mathbf{x}}) \mathbf{\Psi}(\mathbf{x}_m) = \sum_{m=1}^{M}\mathbf{\Psi}(\mathbf{x}_m)^{\textrm{H}} (\dot{\mathbf{x}}_m \cdot \nabla_{\mathbf{x}}) \mathbf{\Psi}(\mathbf{x}_m).
\end{align}

Consider eigenvalues $\{\mu_i\}_{i=1}^L$ and eigenvectors $\{\mathbf{v}_i \}_{i=1}^L$ of $\mathbf{K}_{\textrm{EDMD}}$. Koopman eigenfunctions are in the same form as that in discrete-time formulations while continuous-time Koopman eigenvalues $\mu_i$ can be converted to the aforementioned discrete-time sense as $\lambda_i = e^{\mu_i\Delta t}$.
\subsection{Kernel Dynamic Mode Decomposition}

\subsubsection{Discrete-time Formulation}

Instead of explicitly expressing a dictionary of candidate functions, one can instead implicitly define a dictionary of candidate functions via the kernel trick, which is essentially the dual form of  EDMD~\citep{williams2014kernel}.  Note that, from the EDMD formulation in \cref{eq:cost_func_edmd}, any vector in the range of $\mathbf{K}$ orthogonal to the range of $\mathbf{\Psi}^{\textrm{H}}_{\mathbf{x}}$ is annihilated, and therefore cannot be inferred~\citep{williams2014kernel}. \textcolor{red}{Assuming $\mathbf{\Psi}_{\mathbf{x}}$ is of rank $r$, we can obtain a full SVD  $\mathbf{\Psi}_{\mathbf{x}}=\mathbf{Q}\mathbf{\Sigma}\mathbf{Z}^{\textrm{H}}$ and the corresponding economical-SVD as $\mathbf{Q}_r \mathbf{\Sigma}_r \mathbf{Z}_r^{\textrm{H}}$. With the aforementioned definitions, we  have }$\mathbf{K} = \mathbf{Z}_r \mathbf{\hat{K}} \mathbf{Z}_r^{\textrm{H}}$ without loss~\citep{otto2019linearly}. 

Since the multiplication by a unitary matrix preserves the Frobenius \textcolor{red}{norm}, we  have \begin{align}
J(\mathbf{K}, \{\mathbf{a}_l\}_{l=1}^{L^{'}})  &= \lVert (\mathbf{\Psi}_{\mathbf{x}^{'}} - \mathbf{\Psi}_{\mathbf{x}} \mathbf{K}  )\mathbf{A}^{'} \rVert^2_F \\ &= \lVert (\mathbf{Q}^{\textrm{H}} \mathbf{\Psi}_{\mathbf{x}^{'}} - \mathbf{Q}^{\textrm{H}} \mathbf{Q}_r \mathbf{\Sigma }_r  \mathbf{\hat{K}} \mathbf{Z}^{\textrm{H}}_r )\mathbf{A}^{'} \rVert^2_F \\ &= \lVert (\mathbf{Q}_r^{\textrm{H}} \mathbf{\Psi}_{\mathbf{x}^{'}} - \mathbf{\Sigma}_r \mathbf{\hat{K}} \mathbf{Z}_{r}^{\textrm{H}}) \mathbf{A}^{'}\rVert_F^2 + \lVert \mathbf{Q}^{\textrm{H}}_{r,\perp} \mathbf{\Psi}_{\mathbf{x}^{'}} \mathbf{A}^{'} \rVert_F^2,\end{align} where $\mathbf{Q}^{\textrm{H}}_{r,\perp}$ 
are the last $m-r$ rows of $\mathbf{Q}^{\textrm{H}}$. By taking derivatives with respect to $\mathbf{\hat{K}}$, one can obtain the minimal-norm minimizer as,
\begin{equation}
    \label{eq:kdmd_discrete_time_general}
    \mathbf{\hat{K}}_{\textrm{opt}} = \mathbf{\Sigma}_r^{+} \mathbf{Q}_r^{\textrm{H}} \mathbf{\Psi}_{\mathbf{x}^{'}} \mathbf{A}^{'} \mathbf{A}^{'\textrm{H}} \mathbf{Z}_r (\mathbf{Z}^{\textrm{H}}_r  \mathbf{A}^{'} \mathbf{A}^{'\textrm{H}} \mathbf{Z}_r)^{+}.
\end{equation}

Notice that, since any column in the span of $\mathbf{A}^{'}$ that is orthogonal to the span of $\mathbf{Z}_r$ will be annihilated by $\mathbf{Z}_r^{\textrm{H}}$ and thus cannot be utilized to determine $\mathbf{\hat{K}}$, it is reasonable to restrict $\mathbf{a}_l$ within the column space of $\mathbf{Z}_r$. Assuming  $L^{'}$ is sufficiently large such that the column space of $\mathbf{A}^{'}$ fully spans $\mathbf{Z}_r$, \cref{apdx:equation}  can be proved (details  in \cref{apdx:idtt}),  
\begin{equation}
   \label{apdx:equation}
    \mathbf{A}^{'} \mathbf{A}^{'\textrm{H}} \mathbf{Z}_r (\mathbf{Z}^{\textrm{H}}_r  \mathbf{A}^{'} \mathbf{A}^{'\textrm{H}} \mathbf{Z}_r)^{+} = \mathbf{Z}_r.
\end{equation}
With \cref{apdx:equation}, we can rewrite \cref{eq:kdmd_discrete_time_general} as the familiar KDMD formulation,  
\begin{equation}
    \mathbf{\hat{K}}_{\textrm{KDMD}} = \mathbf{\Sigma}_r^{+} \mathbf{Q}_r^{\textrm{H}} \mathbf{\Psi}_{\mathbf{x}^{'}} \mathbf{Z}_{r} = \mathbf{\Sigma}_r^{+} \mathbf{Q}_r^{\textrm{H}} \mathbf{\Psi}_{\mathbf{x}^{'}} \mathbf{\Psi}_{\mathbf{x}}^{\textrm{H}} \mathbf{Q}_r \mathbf{\Sigma}_r^{+},
\end{equation}
where $\mathbf{\Psi}_{\mathbf{x}} \mathbf{\Psi}_{\mathbf{x}}^{\textrm{H}} = \mathbf{Q}_r \mathbf{\Sigma}^2_r \mathbf{Q}_r^{\textrm{H}}$ with $(\mathbf{\Psi}_{\mathbf{x}} \mathbf{\Psi}_{\mathbf{x}}^{\textrm{H}})_{ij} = k(\mathbf{x}_i, \mathbf{x}_j)$,  $(\mathbf{\Psi}_{\mathbf{x}^{'}} \mathbf{\Psi}_{\mathbf{x}}^{\textrm{H}})_{ij} = k(\mathbf{x}_{i+1}, \mathbf{x}_{j})$ for $1 \le i,j\le M-1$. Again, such a minimizer is independent of \emph{the choice of $\mathbf{A}^{'}$}. 

Notice that, to compute Koopman eigenvalues and eigenfunctions, one would only need access to $\mathbf{\Psi}_{\mathbf{x}^{'}} \mathbf{\Psi}_{\mathbf{x}}^{\textrm{H}}$ and $\mathbf{\Psi}_{\mathbf{x}} \mathbf{\Psi}_{\mathbf{x}}^{\textrm{H}}$, i.e., the inner product among features on all data points. Fortunately, on a compact domain $\mathcal{M}$, it is well-known from Mercer's theorem~\citep{mercer1909functions} that once a suitable non-negative kernel function $k(\cdot, \cdot): \mathcal{M} \times \mathcal{M} \mapsto \mathbb{R}$ is defined, $k(\mathbf{x}, \mathbf{y})$ is the inner product among a potentially infinite dimensional feature vector $\mathbf{\Psi}$ evaluated at $\mathbf{x}, \mathbf{y} \in \mathcal{M}$. Note that the choice of such a feature vector might not be unique but the corresponding reproducing kernel Hilbert space (RKHS) is~\citep{aronszajn1950theory}. In the case of a Gaussian kernel, one can choose the canonical feature vector $k(\cdot,\mathbf{x})$ which are ``bumps" of a certain band-width distributed on $\mathcal{M}$. From the view point of kernel PCA~\citep{williams2014kernel}, $\mathbf{Q}_r$ resulting from the finite dimensional rank truncation on the Gram matrix $\mathbf{\Psi}_{\mathbf{x}} \mathbf{\Psi}_{\mathbf{x}}^{\textrm{H}}$ is a numerical approximation to the $r$ most dominant variance-explained mode shapes in the RKHS evaluated on the data points~\citep{rasmussen2003gaussian}, and $\mathbf{Z}_r$ represents   the  $r$ dominant variance-explaining directions in terms of the feature vector in the RKHS.

Similar to EDMD, given $\mathbf{\hat{K}}_{\textrm{KDMD}} = \mathbf{\hat{V}} \mathbf{\hat{\Lambda}} \mathbf{\hat{V}}^{-1}$, $\mathbf{\hat{V}} = \begin{bmatrix}\mathbf{\hat{v}}_1 \ldots \mathbf{\hat{v}}_r \end{bmatrix}$,  for $i=1,\ldots,r$, the corresponding $i$-th Koopman eigenfunctions $\varphi_i$ and Koopman modes for a vector observable $\mathbf{g}$ are, 
\begin{align}
    \varphi_i(\mathbf{x}) &= \mathbf{\Psi}(\mathbf{x}) \mathbf{\Psi}_{\mathbf{x}}^{\textrm{H}} \mathbf{Q}_r \mathbf{\Sigma}_r^{+} \mathbf{\hat{v}}_i, \\
    \mathbf{B} &=  ( \mathbf{\Psi}_{\mathbf{x}} \mathbf{\Psi}_{\mathbf{x}}^{\textrm{H}} \mathbf{Q}_r \mathbf{\Sigma}_r^{+} \mathbf{\hat{V}} )^{+} 
    \begin{bmatrix}
        \mathbf{g}(\mathbf{x}_1) \\
        \vdots \\
        \mathbf{g}(\mathbf{x}_M) 
    \end{bmatrix}.
\end{align}

\subsubsection{Continuous-time Formulation}

\textcolor{red}{To the best of our knowledge, continuous-time KDMD has not been previously reported in the literature. This can be helpful when non-uniform multi-scale samples are collected.} For the kernel trick to be applied in the continuous formulation,  we write $\mathbf{\Psi}_{\mathbf{x}^{'}}$ as, 
\begin{equation}
    \mathbf{\Psi}_{\mathbf{x}^{'}} = 
    \begin{bmatrix}
        \mathbf{F}(\mathbf{x}_1) \cdot \nabla_{\mathbf{x}} \mathbf{\Psi}(\mathbf{x}_1) \\
        \vdots
        \\
        \mathbf{F}(\mathbf{x}_M) \cdot \nabla_{\mathbf{x}} \mathbf{\Psi}(\mathbf{x}_M)
    \end{bmatrix} =
    \begin{bmatrix}
        \dot{\mathbf{x}}_1 \cdot \nabla_{\mathbf{x}} \mathbf{\Psi}(\mathbf{x}_1) \\
        \vdots
        \\
        \dot{\mathbf{x}}_M \cdot \nabla_{\mathbf{x}} \mathbf{\Psi}(\mathbf{x}_M)
    \end{bmatrix}.
\end{equation}
To compute $\mathbf{\Psi}_{\mathbf{x}^{'}} \mathbf{\Psi}_{\mathbf{x}}^{\textrm{H}}$,  denoting the $q$-th component of $\mathbf{F}$ as $f_q$, 
\begin{align}
\nonumber
(\mathbf{\Psi}_{\mathbf{x}^{'}} \mathbf{\Psi}^{\textrm{H}}_{\mathbf{x}} )_{ij}&= \mathbf{F}(\mathbf{x}_i) \cdot \nabla_{\mathbf{x}} \mathbf{\Psi}(\mathbf{x}_i) \mathbf{\Psi}^{\textrm{H}}(\mathbf{x}_j) \\
\nonumber
&= \sum_{l=1}^L \sum_{q=1}^{N}  \left .  \left( f_q(\mathbf{x})  \dfrac{\partial \psi_l (\mathbf{x})}{\partial x_q}\right)\right\vert_{\mathbf{x} = \mathbf{x}_i}  \left.\overline{\psi_l(\mathbf{x})}\right\vert_{\mathbf{x}=\mathbf{x}_j} \\
\nonumber
&= \sum_{q=1}^{N} f_q(\mathbf{x}_i) \dfrac{\partial}{\partial x_q} \sum_{l=1}^L \left. \left( \psi_l(\mathbf{x})  \overline{\psi_l(\mathbf{x}')} \right) \right\vert_{\mathbf{x}=\mathbf{x}_i, \mathbf{x}'=\mathbf{x}_j} \\
\nonumber
&= \mathbf{F}(\mathbf{x}_i) \cdot \nabla_{\mathbf{x}} k(\mathbf{x}, \mathbf{x}')\vert_{\mathbf{x}=\mathbf{x}_i, \mathbf{x}'=\mathbf{x}_j} \\
& = \dot{\mathbf{x}_i} \cdot \nabla_{\mathbf{x}} k(\mathbf{x}, \mathbf{x}')\vert_{\mathbf{x}=\mathbf{x}_i, \mathbf{x}'=\mathbf{x}_j},
\end{align}
where \textcolor{brown}{the overline symbol is the complex-conjugate operator}, and the appearance of Jacobian indicates that a differentiable kernel function is required for the extension to continuous-time. For common kernels used in Koopman analysis, the kernel function, Jacobian, and hyperparameters are listed in \cref{tab:kernel_table}.

\begin{table}
\begin{center}
\def~{\hphantom{0}}
 \begin{tabular}{lccc}
 kernel type & $k(\mathbf{x}, \mathbf{x}')$ & $\nabla_{\mathbf{x}} k(\mathbf{x}, \mathbf{x}')$ &  hyper para. \\ [3pt]
  linear & $\mathbf{x}\mathbf{x}'^{\textrm{H}}$ & $\mathbf{x}'^{\textrm{H}}$ \\
  polynomial & $ (1+\mathbf{x}\mathbf{x}'^{\textrm{H}})^{\alpha}$ & $\alpha(1+\mathbf{x}\mathbf{x}'^{\textrm{H}})^{\alpha-1} \mathbf{x}'^{\textrm{H}}$ & $\alpha$ \\
  isotropic Gaussian & $\exp(-\lVert \mathbf{x}-\mathbf{x}' \rVert^2/\sigma^2)$ & $\frac{-2(\mathbf{x} - \mathbf{x}')^{\textrm{H}}}{\sigma^2} \exp(-\lVert \mathbf{x}-\mathbf{x}' \rVert^2/\sigma^2) $ & $\sigma$ \\ 
 \end{tabular}
 \caption{Common choice of differentiable kernel functions.}
\label{tab:kernel_table}
 \end{center}
\end{table}

\subsection{Challenges in EDMD/KDMD}
\label{sec:challenges}
In this section, we briefly discuss two broad challenges in the use of EDMD and KDMD for Koopman analysis.
\subsubsection{Mode Selection}

The number of approximated Koopman tuples (eigenfunction, eigenvalue, modes) from EDMD grows with the dictionary size, whereas the KDMD grows with the number of snapshots. However, in most cases, a significant number of the eigenfunctions fail to evolve linearly, or are redundant in contribution to the reconstruction of the state $\mathbf{x}$. For example \textcolor{cyan}{as shown by \cite{budivsic2012applied}}, the Koopman eigenfunctions that vanish nowhere form an Abelian group under pointwise products of functions, while polynomial observables evolve linearly for a general linear system. These eigenfunctions, associated with the polynomial observables, are redundant in terms of providing an intrinsic coordinate for the linear dynamics.

When the number of features is larger than the number of data snapshots, EDMD eigenvalues can be misleading~\citep{otto2019linearly} and often plagued with spurious eigenfunctions that do not evolve linearly even when the number of data snapshots is sufficient. Analytically, it is clear that a Koopman eigenfunction in the span of the dictionary will be associated with one of the eigenvectors obtained from EDMD, given $\mathbf{\Psi}_{\mathbf{x}}$ is full rank, and contains sufficient snapshots $M$~\citep{haseli2019efficient}. Indeed, the EDMD is a $L_2$ projection of the Koopman operator under the empirical measure~\citep{korda2018convergence}. As a result, we \textcolor{brown}{seek a Koopman-invariant subspace following the standard EDMD/KDMD}. Since KDMD \textcolor{brown}{can be viewed as an} efficient way of populating a dictionary of nonlinear features in  high dimensional spaces, \textcolor{cyan}{the above} arguments apply to KDMD as well. It should be noted that \textcolor{brown}{ numerical conditioning} can play a critical role since full rank matrices can be ill-conditioned.

\subsubsection{Choice of Dictionary (for EDMD) or Kernel (for KDMD)}

Although the use of a kernel defines an infinite-dimensional feature space, the resulting finite number of effective features can still be affected by both the type of the kernel and the hyperparameters in the kernel as clearly shown by \citet{kutz2016dynamic}.  Compared to EDMD/KDMD, which are based on a fixed dictionary of features, neural network approaches~\citep{otto2019linearly,pan2019physics,lusch2018deep} have the potential to be more expressive in searching for a larger Koopman-invariant subspace. From a kernel viewpoint~\citep{cho2009kernel}, feedforward neural networks enable adaptation of the kernel function to the data. Such a characteristic could become significant when the underlying Koopman eigenfunction is discontinuous. From an efficiency standpoint, a kernel-guided scalable EDMD~\citep{degennaro2019scalable} may be  pursued. This can be achieved by generating kernel-consistent random Fourier features or approximating a few components of the feature vector constructed from Mercer's theorem, i.e., the eigenfunctions of the Hilbert--Schmidt integral operator on the RKHS.

\section{Sparse Identification of Informative Koopman-invariant Subspace}
\label{sec:select}

\textcolor{red}{To address the challenges described in \cref{sec:challenges}}, we \textcolor{red}{develop} a \textcolor{red}{novel} framework that uses EDMD/KDMD modes to identify a sparse, accurate, and informative Koopman-invariant subspace. \textcolor{green}{Our framework first prunes spurious, inaccurate eigenmodes and second determines a sparse representation of the system state $\mathbf{x}$ from the accurate eigenmodes.} In addition to the training data, as required in standard EDMD/KDMD, a validation trajectory data-set is required to avoid overfitting on training data. The terms spEDMD/spKDMD will refer to filtered mode selections of EDMD and KDMD, respectively. 

\subsection{Pruning spurious modes by a posteriori error analysis}
\label{sec:pruning}

Given a validation trajectory $\mathbf{x}(t)$ where $t \in [0,T]$ associated with the nonlinear dynamical system, for $i =1,\ldots,L$, we define the goodness of $i$-th eigenfunctions in a posteriori way as the maximal normalized deviation from linear evolution conditioned on trajectory $\mathbf{x}(t)$ as $Q_i$ in the form
\begin{align}
\label{eq:linear_evolve_1}
e_{i, \mathbf{x}(0)}(t) &= \dfrac{\lvert \varphi_i(\mathbf{x}(t)) - e^{\lambda_i t} \varphi_i(\mathbf{x}(0)) \rvert}{  \lVert \varphi_i(\mathbf{x}) \rVert_2}, \\
\label{eq:linear_evolve_2}
Q_i \triangleq e^{max}_{i, \mathbf{x}(0)} &= \max_t e_{i, \mathbf{x}(0)}(t),
\end{align}
where $\lVert \varphi_i(\mathbf{x}) \rVert_2 \triangleq \sqrt{\frac{1}{T} \int_{0}^{T} \varphi_i^*(\mathbf{x}(t))\varphi_i(\mathbf{x}(t)) dt}$.
In practice, we evaluate the above integral terms discretely in time.
A similar a priori and less restrictive method has been previously proposed~\citep{zhang2017evaluating}. In contrast, in the proposed method, the maximal error is evaluated in an a posteriori way to better differentiate spurious modes from accurate ones. For any $1 \le \hat{L} \le L$, we can always select top $\hat{L}$ accurate eigenmodes out of $L$ eigenmodes denoting their index in eigen-decomposition as $\{i_1,i_2,\ldots,i_{\hat{L}}\}$, i.e., $Q_{i_1} \le \ldots \le Q_{i_{\hat{L}}} \le \ldots \le Q_{i_L}$. Then, for the next sparse reconstruction step, we simply use $\bm{\varphi}$ defined as follows to reconstruct the state $\mathbf{x}$,
\begin{equation}
    \bm{\varphi}_{\hat{L}}(\mathbf{x}(t)) \triangleq  \begin{bmatrix}\varphi_{ i_1  }(\mathbf{x}(t)) & \ldots & \varphi_{ i_{\hat{L}}  }(\mathbf{x}(t)) \end{bmatrix} \in \mathbb{C}^{\hat{L}}.
\end{equation}
 To choose an appropriate $\hat{L}$ to linearly reconstruct the system state $\mathbf{x}$, we monitor the normalized reconstruction error for the aforementioned set of top $\hat{L}$ accurate eigenmodes in the following form
\begin{equation}
\label{eq:rec_loss}
R_{\hat{L}} \triangleq \dfrac{   \lVert (\mathbf{I}  - \mathbf{\Psi}_{\hat{L}}\mathbf{\Psi}_{\hat{L}}^{+}) \mathbf{X} \rVert_F}{ \lVert \mathbf{X} \rVert_F },
\end{equation}
where $\mathbf{I}$ is the identity matrix, and
\begin{equation}
\mathbf{X} = 
\begin{bmatrix}
\mathbf{x}_1\\
\vdots\\
\mathbf{x}_M
\end{bmatrix}, \quad 
\mathbf{\Psi}_{\hat{L}} = 
\begin{bmatrix}
\bm{\varphi}_{\hat{L}} (\mathbf{x}_1)\\
\vdots\\
\bm{\varphi}_{\hat{L}} (\mathbf{x}_M)
\end{bmatrix}.
\end{equation}
As a result, the evaluation of \cref{eq:rec_loss} for each $\hat{L}$ is of similar expense to least-square regression. For an increasing number of selected eigenfunctions $\hat{L}$, the reconstruction error $R_{\hat{L}}$ decreases, while the largest linear evolution error $Q_{i_{\hat{L}}}$ increases. Then, a truncation $\hat{L}$ can be defined by the user to strike a balance between linear evolution error $Q_{i_{\hat{L}}}$ and reconstruction error $R_{\hat{L}}$. \textcolor{cyan}{In the next subsection, we will further} select a subset of eigenmodes for spanning the minimal Koopman-invariant subspace.

\subsection{Sparse reconstruction via multi-task feature learning}
\label{sec:multi-task}
Numerical experiments revealed that, in the  selected set of $\hat{L}$ most accurate eigenfunctions,  two types of redundant eigenfunctions were found:
\begin{enumerate}
	\item Nearly constant eigenfunctions with eigenvalues close to zero,
	\item Pointwise products of Koopman eigenfunctions introduced by nonlinear observables, not useful in linear reconstruction.
\end{enumerate}

To filter the above modes, \textcolor{green}{we} consider sparse regression with $\hat{L}$ most accurate eigenfunctions as features and the system state $\mathbf{x}$ as target. 
Note that, since we have guaranteed the accuracy of selected eigenmodes, one can either choose features a priori: $\varphi_i(\mathbf{x}(t))$ or a posteriori \textcolor{cyan}{(multi-step prediction)} $e^{\lambda_i t } \varphi_i(\mathbf{x}(0))$. Here we choose the latter since it is directly related to prediction, and can actually be reused from the previous step without additional computational cost. We denote the corresponding \textcolor{cyan}{multi-step prediction} feature matrix as $\mathbf{\hat{\Psi}}_{\hat{L}}$ ,  
\begin{equation}
\label{eq:aposter_feature}
\mathbf{\hat{\Psi}}_{\hat{L}} = 
\begin{bmatrix}
\bm{\varphi}_{\hat{L}} (\mathbf{x}_1) \\
\bm{\varphi}_{\hat{L}} (\mathbf{x}_1) e^{\Delta t\Lambda_{\hat{L}}}\\
\vdots\\
\bm{\varphi}_{\hat{L}} (\mathbf{x}_1) e^{(M-1)\Delta t\Lambda_{\hat{L}}}
\end{bmatrix} \in \mathbb{C}^{M \times \hat{L}},
\end{equation}
where $\Lambda_{\hat{L}} = \textrm{diag}(\lambda_{i_1},\ldots,\lambda_{i_{\hat{L}}})$. Note that similar features $\mathbf{\hat{\Psi}}_{\hat{L}}$ were also considered in sparsity-promoting DMD~\citep{jovanovic2014sparsity} and optimized DMD~\citep{chen2012variants}. Finally, the fact that there is no control over the magnitudes of the implicitly defined features in the standard KDMD may cause unequal weighting between different features. Thus, we consider scaling the initial value of all eigenfunctions to be unity in \cref{eq:scaling},
\begin{equation}
\label{eq:scaling}
\mathbf{\hat{\Psi}}_{\hat{L},\textrm{scaled}} = \mathbf{\hat{\Psi}}_{\hat{L}} \mathbf{\Lambda}^{-1}_{ini} = 
\begin{bmatrix}
1 & \ldots & 1 \\
e^{\Delta t\lambda_{i_1}} & \ldots & e^{\Delta t\lambda_{i_{\hat{L}}} }\\
\vdots & \vdots & \vdots \\
e^{({M-1})\Delta t\lambda_{i_1}} & \ldots & e^{({M-1})\Delta t\lambda_{i_{\hat{L}}} }
\end{bmatrix},
\end{equation}
where 
\begin{equation}
\mathbf{\Lambda}_{ini} = \textrm{diag}(\begin{bmatrix} \varphi_{i_1}(\mathbf{x}_1) & \ldots &  \varphi_{i_{\hat{L}}}(\mathbf{x}_1) \end{bmatrix}).
\end{equation}

Since $\mathbf{x}$ is finite-dimensional, searching for a sparse combination of $\mathbf{\hat{\Psi}}_{\hat{L}}$ to reconstruct $\mathbf{x}$ is \emph{equivalent} to the solution of a multi-task feature learning problem with preference over a relatively small size of features. Note that this type of problem has been studied extensively in the machine learning community~\citep{argyriou2008convex,zhao2015multi,argyriou2008spectral}.
In this work, given $\mathbf{X}$ and $\mathbf{\hat{\Psi}}_{\hat{L},\textrm{scaled}}$, we \textcolor{red}{leverage} the multi-task ElasticNet~\citep{scikit-learn} to search for a row-wise sparse $\mathbf{B}^{'}$, which solves the following convex optimization problem:
\begin{equation}
\label{eq:multi_task_enet}
\mathbf{B}^{'*} = \argmin_{\mathbf{B}^{'} \in \mathbb{C}^{\hat{L} \times N}} \frac{1}{2M}\lVert \mathbf{X} - \mathbf{\hat{\Psi}}_{\hat{L},\textrm{scaled}}\mathbf{B}^{'} \rVert_{F}^2 + \alpha \rho \lVert \mathbf{B}^{'} \rVert_{2,1} + \frac{\alpha (1-\rho)}{2}  \lVert \mathbf{B}^{'} \rVert_F^2,
\end{equation}
and 
\begin{equation}
\mathbf{B} = \mathbf{\Lambda}^{-1}_{ini}\mathbf{B}^{'*},
\end{equation}
where $\lVert \cdot \rVert_{2,1}$ defined in \cref{eq:def_l21_norm} is the so-called $\ell_1/\ell_2$ norm for a matrix $\mathbf{W}$,
\begin{equation}
\label{eq:def_l21_norm}
\lVert \mathbf{W} \rVert_{2,1} \triangleq \sum_{i} \sqrt{\sum_{j}\mathbf{W}_{ij}^2} = \sum_{i} \lVert \mathbf{w}_i \rVert_2,
\end{equation}
and $\mathbf{W}_i$ is $i$-th row of $\mathbf{W}$. This norm is special in that it controls the number of shared features learned across all tasks, i.e., $i$-th Koopman mode $\mathbf{b}_i$ is either driven to a zero vector or not while the standard $\ell_1$ only controls the number of features for each task independently. 

\begin{figure}
\centering
\includegraphics[width=0.8\textwidth]{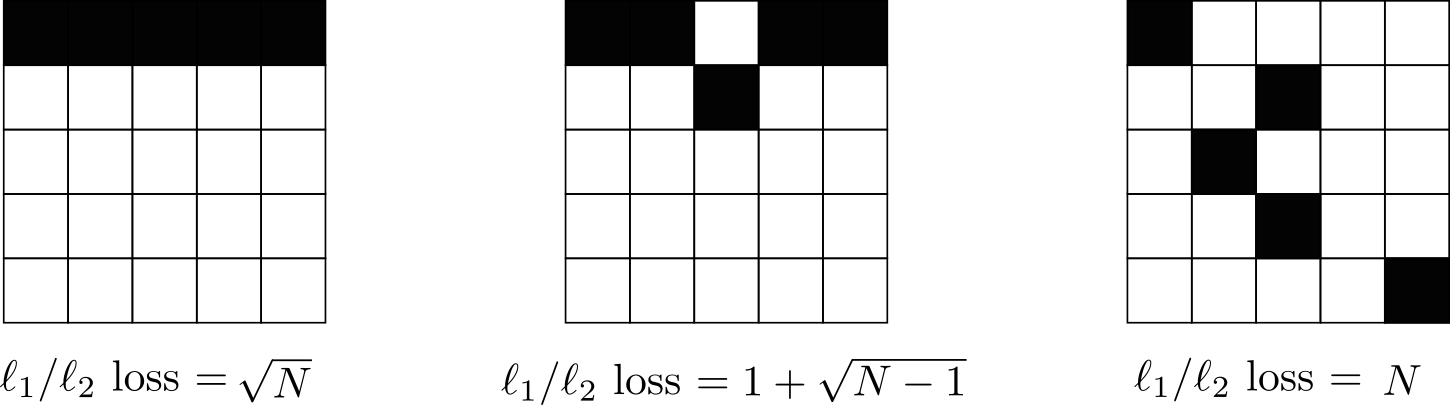}
\caption{Illustration of $\ell_1/\ell_2$ norm \textcolor{brown}{(defined in \cref{eq:def_l21_norm})} for different $N\times N$ \textcolor{brown}{0-1} binary matrices.}
\label{fig:l21}
\end{figure}

As a simple illustration, the $\ell_1/\ell_2$ norm for three different $N \times N$ square matrices (here $N=5$) with 0-1 binary entries is displayed in \cref{fig:l21}. \textcolor{brown}{Since $\sqrt{N} \le 1 + \sqrt{N-1} \le N$, minimizing the $\ell_1/\ell_2$ norm leads to a penalty on the number of rows. As shown in the second term on the right hand side of \cref{eq:multi_task_enet}, minimizing the $\ell_1/\ell_2$ norm penalizes the number of Koopman eigenmodes.}

The above procedure not only serves the purpose of selecting modes that explain the behavior of all components in the state, but is also particularly natural for EDMD/KDMD since Koopman modes are obtained via regression. 
$\alpha$ is a penalty coefficient that controls the amount of total regularization in the $\ell_1/\ell_2$ and $\ell_2$ norms, while $\rho$ is the ElasticNet mixing parameter~\citep{zou2005regularization} that ensures uniqueness of the solution when highly correlated features exist. In our case, we choose $\rho = 0.99$ and sweep $\alpha$ over a certain range with $L_r$ non-zero features denoted as $\mathbf{\hat{\Psi}}_{L_r}$ for each $\alpha$, while monitoring the normalized residual $\min_{\mathbf{B}}\lVert \mathbf{X} - \mathbf{\hat{\Psi}}_{L_r}\mathbf{B} \rVert_F/\lVert \mathbf{X} \rVert_F$ to choose an appropriate $\alpha$. 
It has to be mentioned that, sometimes the sparsest solution from a multi-task ElasticNet was found to shrink to a small number instead of zero. This is a consequence of the insufficiency of the current optimization algorithm which employs coordinate descent~\citep{scikit-learn}. Hence for each target component, we consider an additional hard-thresholding step  by setting the corresponding magnitude of the coefficient, i.e., contribution of any mode, to zero if it is smaller than a certain threshold $\epsilon \in [10^{-2},10^{-3}]$. 

\begin{figure}
\centering
\includegraphics[width=0.6\textwidth]{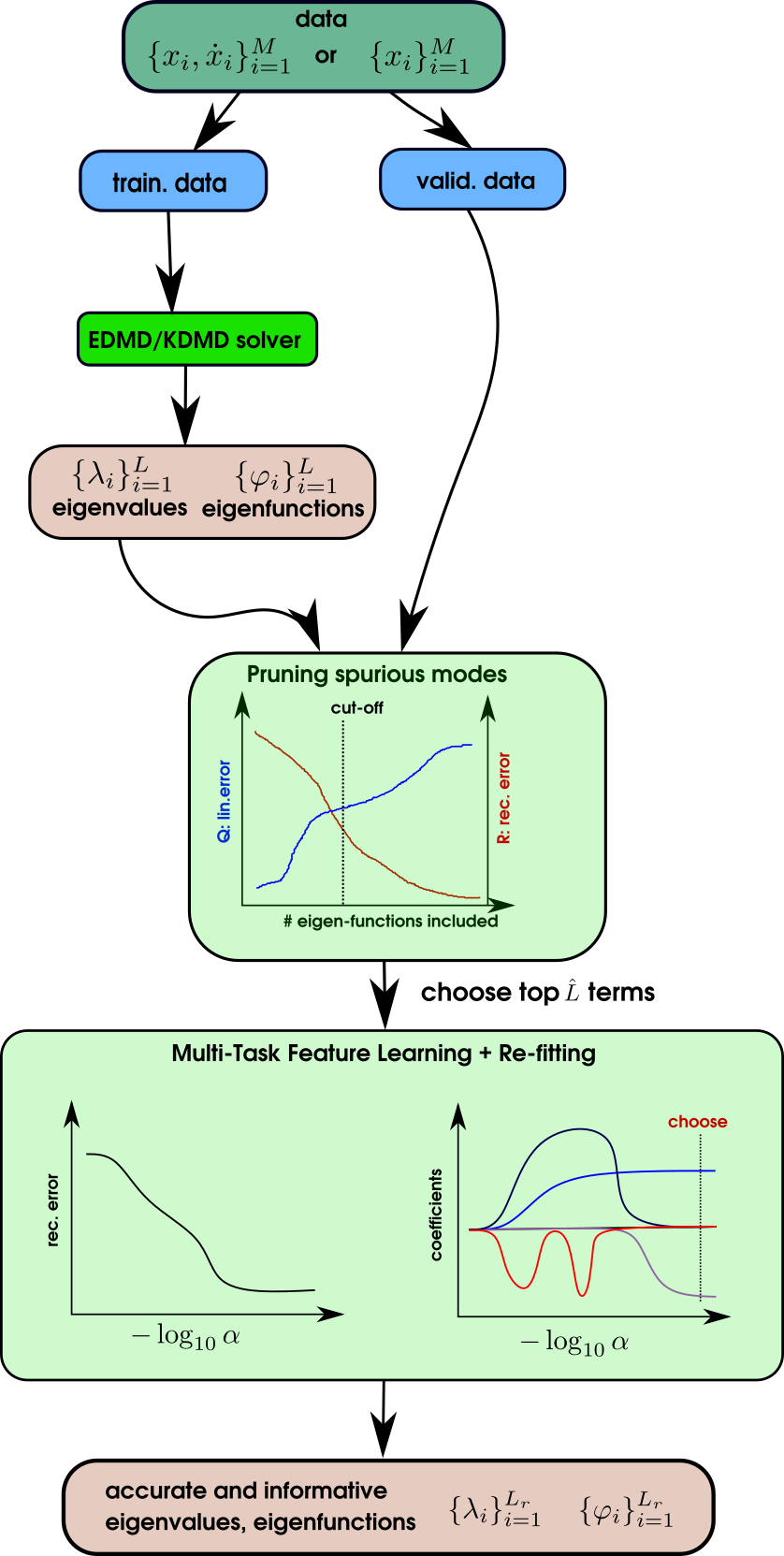}
\caption{Schematic illustrating the idea of sparse identification of  Koopman-invariant subspaces for EDMD and KDMD.}
\label{fig:framework}
\end{figure}

Finally, we refit the Koopman modes as $\mathbf{B}_{L_r} = \mathbf{\hat{\Psi}}_{L_r}^{+}\mathbf{X}$ which avoids the bias introduced by the penalty term~\footnote{spDMD does not refit $\mathbf{B}$} in \cref{eq:multi_task_enet}.
To summarize, the general idea of the framework is illustrated in \cref{fig:framework}. As a side note for interested readers, if one only performs multi-task feature learning without  hard-thresholding and refitting, one would obtain a smooth ElasticNet path instead of a discontinuous one with hard-thresholding and refitting. However, the smooth ElasticNet can lead to difficulties in choosing the proper $\alpha$ visually, especially when the given dictionary of EDMD/KDMD is not rich enough to cover an informative Koopman-invariant subspace. \textcolor{red}{Further discussion on the computational complexity of our framework is presented in \cref{apdx:cc}.}

\textcolor{red}{Thus far, we have presented our main contribution: a novel optimization-based  framework to search for an accurate and minimal Koopman-invariant subspace from data. An appealing aspect of our framework is the \emph{model agnostic property}, which makes the extension easy from the standard EDMD/KDMD to more advanced approximation methods~\citep{jungers2019non,mamakoukas2019local,azencot2019consistent}. In the following subsection, we present two mathematical insights: 1) multi-task feature learning generalizes spDMD under a specific constraint; 2) a popular empirical criterion can be viewed as a single step of proximal gradient descent.}

\subsubsection{Relationship between sparsity-promoting DMD, Kou's criterion and multi-task feature learning}

For simplicity, neglecting the ElasticNet part (i.e., using $\rho=1$),  \cref{eq:multi_task_enet} with $L$ modes leads to a multi-task Lasso problem, 
\begin{equation}
\label{eq:multi_task_enet_simpler}
\min_{\mathbf{B}^{'}\in\mathbb{C}^{L \times N}} \frac{1}{2M}\lVert \mathbf{X} - \mathbf{\hat{\Psi}}_{L,\textrm{scaled}}\mathbf{B}^{'} \rVert_{F}^2 + \alpha  \lVert \mathbf{B}^{'} \rVert_{2,1}.
\end{equation}
Recall that in spDMD~\citep{jovanovic2014sparsity}, DMD modes $\phi_1,\ldots,\phi_L$ with $\lVert \phi_i \rVert_2 = 1$ remain the same as standard DMD. Similariy, if we posit a structural constraint on $\mathbf{B}^{'}$ in \cref{eq:multi_task_enet_simpler} by enforcing the modes as those from DMD, then there exist $\alpha_1,\ldots,\alpha_{L}$ such that, 
\begin{equation}
\mathbf{B}^{'} = 
\begin{bmatrix}
\alpha_1 & & \\
& \ddots & \\
& & \alpha_{L}
\end{bmatrix}
\begin{bmatrix}
\phi_1^{\textrm{T}} \\
\vdots\\
\phi_{L}^{\textrm{T}}
\end{bmatrix}.
\end{equation}
Note the fact that $\lVert \mathbf{B}^{'} \rVert_{2,1} = \sum_{i}^{L} \lvert \alpha_i \rvert$. Hence, we recover the $\ell_1$ optimization step in the spDMD~\citep{jovanovic2014sparsity} from \cref{eq:multi_task_enet_simpler} as, 
\begin{align}
\label{eq:spdmd}
&\min_{\alpha_1,\ldots,\alpha_{L} \in \mathbb{C}} \frac{1}{2M}\Bigg\Vert \mathbf{X} - 
\begin{bmatrix}
1 & \ldots & 1 \\
\vdots & \vdots & \vdots \\
e^{({M-1})\Delta t\lambda_{i_1}} & \ldots & e^{({M-1})\Delta t\lambda_{i_{L}} }
\end{bmatrix} 
\begin{bmatrix}
\alpha_1 & & \\
& \ddots & \\
& & \alpha_{L}
\end{bmatrix}
\begin{bmatrix}
\phi_1^{\textrm{T}} \\
\vdots\\
\phi_{L}^{\textrm{T}}
\end{bmatrix}
\Bigg\Vert^2_F + \alpha \sum_{i=1}^{L}\lvert \alpha_i \rvert,
\end{align}
where $\phi_1,\ldots,\phi_{L}$ are the DMD modes corresponding to eigenvalues as $\lambda_1,\ldots,\lambda_{L}$. Hence, multi-task feature learning solves a less constrained optimization than spDMD in the context of DMD. 

\citet{kou2017improved} proposed an empirical criterion for mode selection by ordering modes with ``energy" $I_i$ defined as 
\begin{equation}
\label{eq:kou_criterion}
I_i = \sum_{j=1}^{M} |\alpha_i e^{(j-1)\Delta t \lambda_i}| = 
\begin{cases}
      \frac{|\alpha_i|(1-\lvert e^{\Delta t \lambda_i}\rvert^M)}{1-\lvert e^{\Delta t \lambda_i}\rvert},& \text{if } \lvert e^{\Delta t \lambda_i} \rvert \neq 1,\\
    M|\alpha_i| ,              & \text{otherwise},
\end{cases}.
\end{equation} 
From an optimization viewpoint, consider a posteriori prediction matrix $\mathbf{X}_{\textrm{DMD}}$  from DMD
\begin{equation}
\label{eq:DMD_matrix}
\mathbf{X} \approx \mathbf{X}_{\textrm{DMD}} = 
\begin{bmatrix}
1 & \ldots & 1 \\
\vdots & \vdots & \vdots \\
e^{({M-1})\Delta t\lambda_{1}} & \ldots & e^{({M-1})\Delta t\lambda_{{L}} }
\end{bmatrix} 
\begin{bmatrix}
\alpha_1 & & \\
& \ddots & \\
& & \alpha_{L}
\end{bmatrix}
\begin{bmatrix}
\phi_1^{\textrm{T}} \\
\vdots\\
\phi_{L}^{\textrm{T}}
\end{bmatrix},
\end{equation}
where $\mathbf{X}_{\textrm{DMD}}$ is determined from DMD using the snapshot  pair  $(\mathbf{X},\mathbf{X}^{'})$.  $\mathbf{X}_{\textrm{DMD}}$ is a rank-1 summation of contributions from different modes~\citep{schmid2010dynamic}.
Hence, a general  mode selection technique with a user-defined preference weighting $\mathbf{w}$ is the following weighted $\ell_0$ nonconvex optimization problem:
\begin{align}
\label{eq:dmd_kou}
& \min_{\mathbf{a} \in \mathbb{C}^{L}}\Bigg\lVert \mathbf{X}_{\textrm{DMD}} -  \begin{bmatrix}
1 & \ldots & 1 \\
\vdots & \vdots & \vdots \\
e^{({M-1})\Delta t\lambda_{1}} & \ldots & e^{({M-1})\Delta t\lambda_{{L}} }
\end{bmatrix} 
\textrm{diag}(\mathbf{a})
\begin{bmatrix}
\phi_1^{\textrm{T}} \\
\vdots\\
\phi_{L}^{\textrm{T}}
\end{bmatrix}\Bigg\rVert^2_F + \lambda \lVert \mathbf{a} \rVert_{\mathbf{w},0}
\end{align}
where $\lVert \mathbf{a} \rVert_{\mathbf{w},0} \triangleq \sum_{i} w_i|a_i|^{0}$, $|a_i|^0$ is one if $a_i \neq 0$ and zero otherwise. 
Note that this \textcolor{cyan}{pseudo-norm} can be viewed as a limiting case of a weighted composite sine function smoothed $\ell_0$ regularization~\citep{wang2019re}.

To solve this non-convex optimization problem, compared to the popular $\ell_1$ relaxation method such as the one in sparsity-promoting DMD, a less-known but rather efficient way is iterative least-squares hard thresholding. This  has been used in sparse identification of dynamical systems (SINDy)~\citep{brunton2016discovering}, and convergence to a local minimum has been proved~\citep{zhang2019convergence}. Indeed, a more rigorous framework that generalizes such an algorithm is the proximal gradient method~\citep{parikh2014proximal}.  Much like Newton's method is a standard tool for unconstrained smooth optimization, the proximal gradient method is the standard tool for constrained non-smooth optimization. Here, it is straightforward to derive the iterative algorithm that extends to the weighted $\ell_0$ norm from step $k$ to step $k+1$ as 
\begin{equation}
\label{eq:opt_iht}
\mathbf{a}^{k+1} = \textrm{prox}_{\frac{\lambda}{2}\eta_k \lVert \cdot \rVert_{\mathbf{w},0}} (\mathbf{a}^k - \eta_k \nabla_{\mathbf{a}}\mathcal{Q}(\mathbf{a}^{k}) ),
\end{equation}
where 
\begin{equation}
\mathcal{Q}(\mathbf{a}) = \frac{1}{2}\Bigg\lVert \mathbf{X}_{\textrm{DMD}} -  \begin{bmatrix}
1 & \ldots & 1 \\
\vdots & \vdots & \vdots \\
e^{({M-1})\Delta t\lambda_{1}} & \ldots & e^{({M-1})\Delta t\lambda_{{L}} }
\end{bmatrix} 
\textrm{diag}(\mathbf{a})
\begin{bmatrix}
\phi_1^{\textrm{T}} \\
\vdots\\
\phi_{L}^{\textrm{T}}
\end{bmatrix}\Bigg\rVert^2_F,
\end{equation}
and $\eta_k$ is the step-size at step $k$.
Notice that the weighted $\ell_0$ norm is a separable sum of $a_i$. After some algebra, we have the proximal operator as
\begin{equation}
\textrm{prox}_{\frac{\lambda}{2}\eta_k \lVert \cdot \rVert_{\mathbf{w},0}} (\mathbf{a}) = 
\begin{bmatrix}
\mathcal{H}_{\sqrt{\lambda \eta_k}}(a_1/\sqrt{w_1}) & 
\ldots &
\mathcal{H}_{\sqrt{\lambda \eta_k}}(a_L/\sqrt{w_L})
\end{bmatrix}^{\textrm{T}},
\end{equation}
where $\mathcal{H}_{\sqrt{\lambda \eta_k} }(a)$ is an element-wise hard thresholding operator defined as $a$ if $|a| < \sqrt{\lambda \eta_k}$ and zero otherwise. 

Particularly, if one considers the initial step-size to be extremely small $\eta_1 \ll 1$, then the second term in \cref{eq:opt_iht} can be neglected. Thus, for $i=1,\ldots,L$, with the following weighting scheme that penalizes fast decaying modes:
\begin{equation}
\label{eq:analytic_final_relation}
    w_i = 1/\beta_i^{2}, \quad \beta_i = 
\begin{cases}
     \frac{1-\lvert e^{\Delta t \lambda_i}\rvert^M}{1-\lvert e^{\Delta t \lambda_i}\rvert},& \text{if } \lvert e^{\Delta t \lambda_i} \rvert \neq 1,\\
    M,              & \text{otherwise},
\end{cases}
\end{equation}
one immediately realizes the thresholding criterion for $i$-th entry of $\mathbf{a}$ becomes 
\begin{equation}
\sqrt{\lambda \eta_k} > |\alpha_i/\sqrt{w_i}| = |\alpha_i\beta_i|. 
\end{equation}
Then \textcolor{cyan}{plugging \cref{eq:analytic_final_relation} in \cref{eq:opt_iht}},  the first iteration in \cref{eq:opt_iht} reduces to mode selection with Kou's criterion in \cref{eq:kou_criterion}. Normally, $\beta_i$ is very large for unstable modes and small for decaying modes. It is important to note that a) such a choice of $\mathbf{w}$ preferring unstable/long-lasting modes over decaying modes is still user-defined; 2) Optimization is in an a priori sense to obtain DMD. Thus, the insufficiency of the a priori formulation to account for  temporal evolution  is indeed \emph{compensated} by this criterion, while DMD in an a posteriori formulation (e.g., sparsity-promoting DMD) includes such a effect implicitly in the optimization. Hence, it is possible that in some circumstances spDMD and Kou's criterion could achieve similar performance~\citep{kou2017improved}. 

\begin{figure}
\centering
\includegraphics[width=0.9\textwidth]{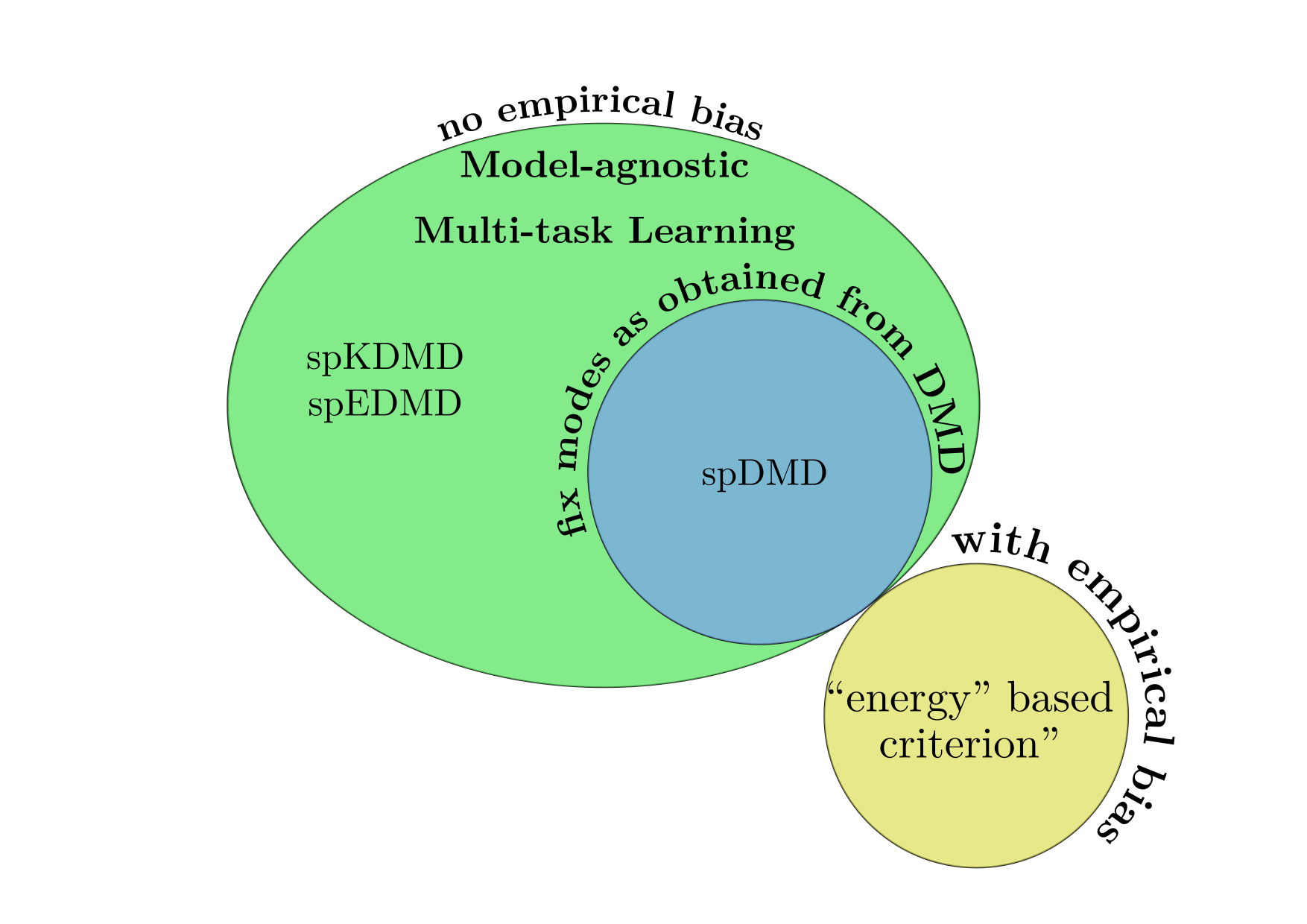}
\caption{Differences and similarities among existing mode selection methods.}
\label{fig:analytic_rel}
\end{figure}

Lastly, \textcolor{cyan}{as summarized in \cref{fig:analytic_rel},} it is important to mention the similarities and differences between spKDMD and spDMD: 1) 
spKDMD will refit Koopman modes while spDMD does not; and 2) The amplitude for all the modes in spKDMD is fixed as unity while it has to be determined in spDMD. 

\subsection{Hyper-parameter selection}
\label{sec:hyp_cv}

 For simplicity, hyper-parameter selection for KDMD is only discussed in this section. To fully determine a kernel in KDMD, one would have to choose the following:
\begin{enumerate}
    \item kernel type,
    \item kernel parameters, e.g., scale parameters $\sigma$,
    \item rank truncation $r$.
\end{enumerate}
In this work, for simplicity, we fix the kernel type to be an isotropic Gaussian.
Motivated by previous work on error evaluation in Koopman modes by \citet{zhang2017evaluating}, we consider evaluation with cross validation on a priori mean normalized accuracy defined in \cref{eq:aprior_loss_dis,eq:aprior_loss_cont} for $i$-th eigenfunction,
\begin{equation}
\label{eq:aprior_loss_dis}
\textrm{discrete form: } {Q}_i^{a} = \frac{1}{M-1}\sum_{j=1}^{M-1} \frac{|\varphi_i(\mathbf{x}_{j+1}) - \lambda_i \varphi_i(\mathbf{x}_j) |}{\sqrt{\frac{1}{M}\sum_{k=1}^M \varphi_i^*(\mathbf{x}_k) \varphi_i(\mathbf{x}_k) } } 
\end{equation}
\begin{equation}
\label{eq:aprior_loss_cont}
\textrm{continuous form: } {Q}_i^{a} =  \frac{1}{M}\sum_{j=1}^{M} \frac{|\dot{\mathbf{x}}_{j} \cdot \nabla_{\mathbf{x}} \varphi_i(\mathbf{x}_{j}) - \lambda_i \varphi_i(\mathbf{x}_j) |}{\sqrt{\frac{1}{M}\sum_{k=1}^M \varphi_i^*(\mathbf{x}_k) \varphi_i(\mathbf{x}_k) } }   
\end{equation}
on validation data for different number of rank truncation and kernel parameters.

Note that evaluation on maximal instead of mean normalized accuracy would lead to the error metric to be strongly dependent on the local sparsity of training data in the feature space. This is particularly true for a single trajectory for a high-dimensional dynamical system is used, and random shuffled cross validation is performed~\citep{pan2018long}. 

For each set of hyperparameters, we first compute the number of eigenfunctions of which the error defined in \cref{eq:aprior_loss_dis,eq:aprior_loss_cont} is below a certain threshold on both training and validation data for each fold of cross validation. Then we compute the average number of such eigenfunctions over all folds which indicates the quality of the corresponding subspace. Finally, we plot the average number versus rank truncation $r$ and kernel scale parameters $\sigma$ to select hyperparameters.


%

\subsection{Implementation}
\label{sec:imple}

We implement the described framework in Python with moderate parallelism in each module. We use \texttt{scipy.special.hermitenorm}~\citep{scipy} to generate normalized Hermite polynomials and \texttt{MultiTaskElasticNet} in the \texttt{scikit-learn}~\citep{scikit-learn} for multi-task feature learning where we set the maximal iteration as $10^5$ and tolerance as $10^{-12}$.
MPI parallelism using \texttt{mpi4py}~\citep{mpi4py} is used for the grid search in hyperparameter selection. To prepare data with hundreds of gigabytes collected from high fidelity simulations, a  distributed SVD written in C++ named Parallel Data Processing (PDP) Tool is developed for dimension reduction. A brief description of this tool is given in \cref{apdx:pdp}.

\section{Results and Discussion}
\label{sec:result}

\subsection{2D fixed point attractor with the known finite Koopman-invariant subspace}

We first consider a simple fixed point nonlinear dynamical system which has an  analytically determined, finite-dimensional non-trivial Koopman-invariant subspace~\citep{brunton2016koopman,kaiser2017data} to show the effectiveness of proposed method. We consider a continuous-time formulation. The governing equation for the dynamical system is given as follows, 
\begin{align}
    \dot{x}_1 &= \mu x_1, \\
    \dot{x}_2 &= \lambda (x_2 - x_1^2),
\end{align}
where $\mu = -0.05, \lambda = -1$. One natural choice of the minimal Koopman eigenfunctions that \emph{linearly} reconstructs the state  is~\citep{brunton2016koopman}
\begin{equation}
\label{eq:2f_analytic}
\varphi_1(\mathbf{x}) = x_2 - \lambda x_1^2 /(\lambda-2\mu), \quad \varphi_2(\mathbf{x})= x_1, \quad \varphi_3(\mathbf{x}) = x_1^2
\end{equation}
with eigenvalues $\lambda=-1$, $\mu=-0.05$, $2\mu=-0.1$ respectively. 

\begin{figure}
\centering
\includegraphics[width=0.6\textwidth]{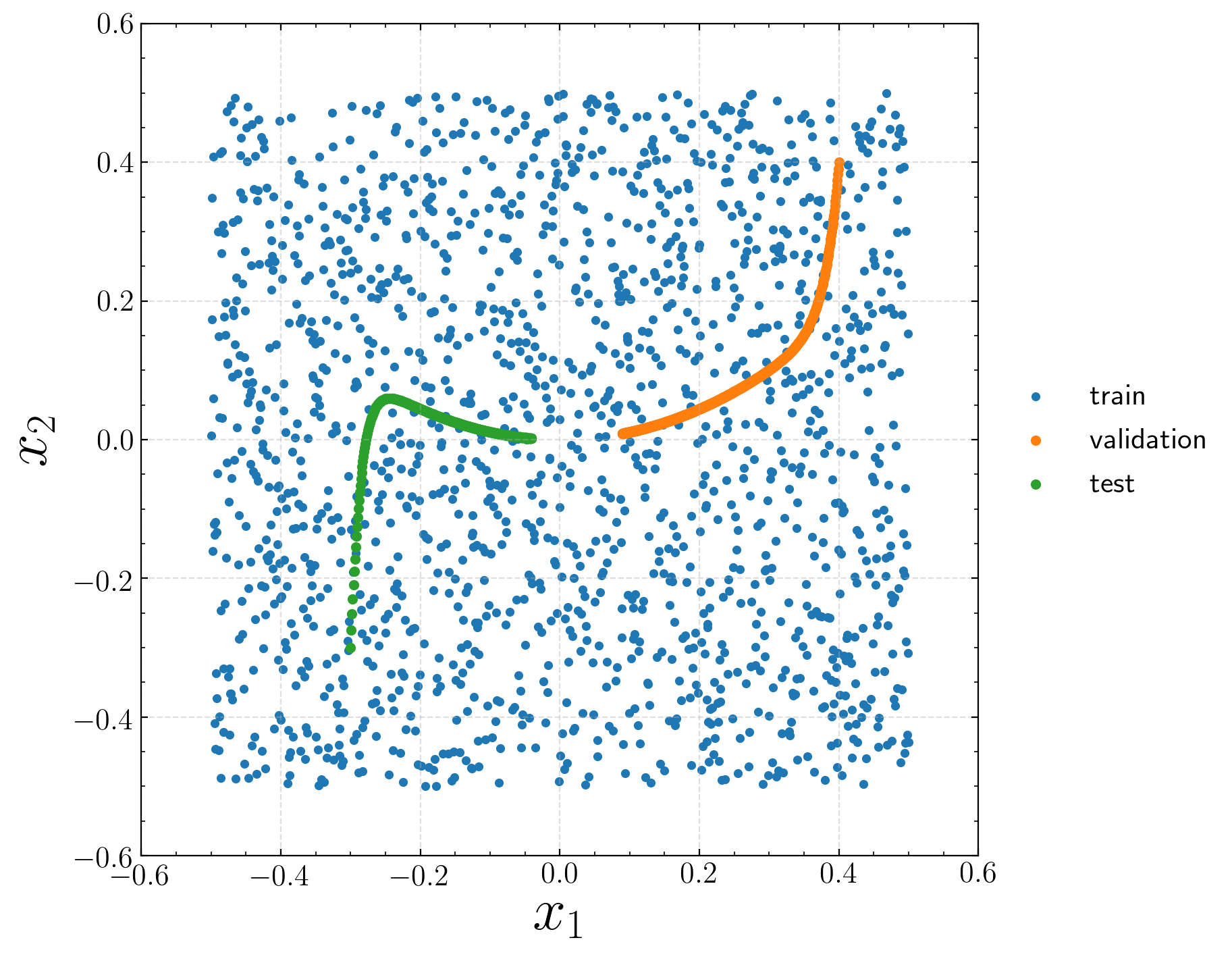}
\caption{Data distribution for 2D fixed point attractor.}
\label{fig:2d_fp_data}
\end{figure}

The way we generate training, validation, and testing data is described below with distribution of the data shown in \cref{fig:2d_fp_data}, 
\begin{enumerate}
\item \emph{training data}: a point cloud with $M = 1600$ pairs of $\{\mathbf{x}^{(i)}, \mathbf{\dot{x}}^{(i)}\}_{i=1}^{M}$, is generated by Latin hypercube sampling~\citep{baudin2015pydoe} within the domain $x_1,x_2 \in [-0.5,0.5]$. 
\item \emph{validation data}: a single trajectory with initial condition as $x_1(0)=x_2(0)=0.4$, sampling time interval $\Delta t = 0.03754$ from $t=0$ to $t=30$.

\item \emph{testing data}: a single trajectory with initial condition as $x_1(0)=x_2(0) = -0.3$, sampling time interval $\Delta t =0.06677$ from $t=0$ to $t=40$.

\end{enumerate}

As an illustration, we consider two models to approximate the Koopman operator from training data: 
\begin{enumerate}
\item a continuous-time EDMD with normalized Hermite polynomials up to fifth order with $L=36$ features, 
\item a continuous-time KDMD with isotropic Gaussian kernel $\sigma = 2$ with reduced rank $r = L = 36$.
\end{enumerate}

Details of the above choices based on the steps of hyperparameter selection in \cref{sec:hyp_cv} are given in \cref{apdx:hyper_2d}. 
\subsubsection{Results for continuous-time EDMD with mode selection}

\begin{figure}
\centering
\includegraphics[width=\textwidth]{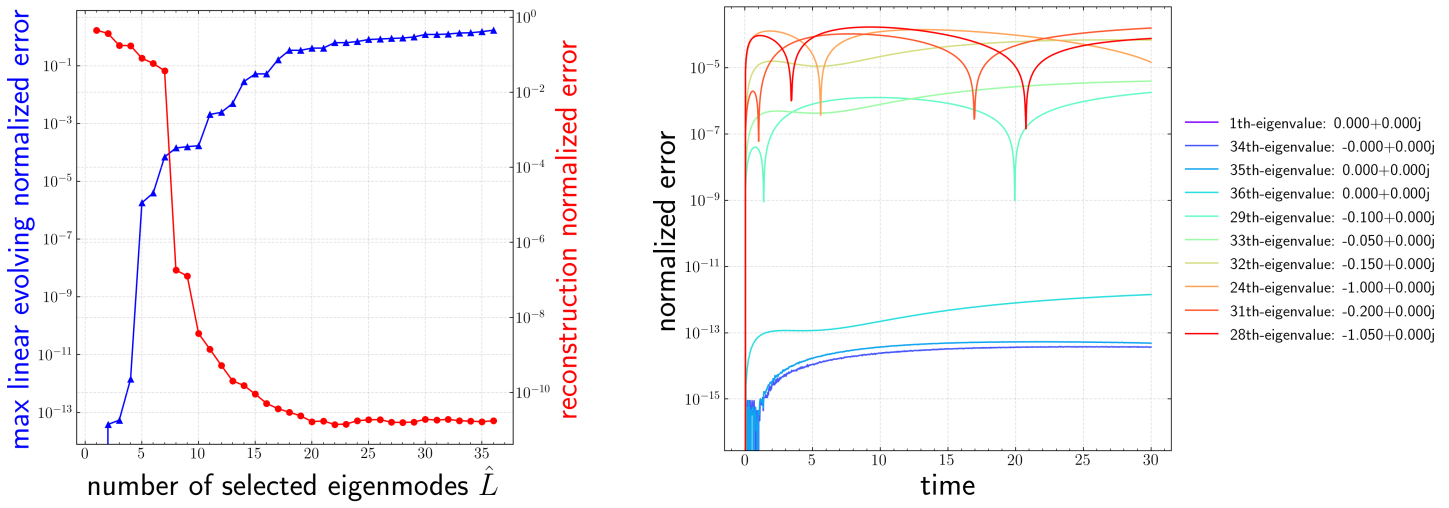}
\caption{Error analysis of 36 eigenmodes from continuous-time EDMD for the 2D fixed point attractor. Left: trends of linearly evolving error $Q$ and reconstruction error $R$. Right: temporal evolution of relative error for top $\hat{L}=10$ accurate eigenmodes.}
\label{fig:2d_fp_edmd_error_analysis}
\end{figure}

As displayed in \cref{fig:2d_fp_edmd_error_analysis}, we begin with an error analysis of all of the eigenmodes on validation data in \cref{fig:2d_fp_data} according to linearly evolving error $Q$ defined in \cref{eq:linear_evolve_2} and $R$ defined in \cref{eq:rec_loss}. From \textcolor{cyan}{the left subfigure in} \cref{fig:2d_fp_edmd_error_analysis}, considering both the linearly evolving error and the quality of the reconstruction, we choose the cut-off threshold at $\hat{L}=10$. We observe a sharp cut-off in the left subfigure in \cref{fig:2d_fp_edmd_error_analysis} around the number of selected eigenmodes $\hat{L} = 8$. This is a reasonable choice, since from the eigenvalues in the right subfigure in  \cref{fig:2d_fp_edmd_error_analysis}, we notice the analytic Koopman eigenmodes are not covered until first 8 accurate eigenmodes are selected. \textcolor{cyan}{Note that the legend in the right subfigure is ordered by the maximal deviation from linear evolution, e.g., the second most accurate mode is 34-th mode with zero eigenvalue.}  Indeed, the first four eigenfunctions (index=$1,34,35,36$) are redundant in terms of reconstruction in this problem~\footnote{This could be interesting if the system is instead Hamiltonian.}. The fifth (index=$29$) and sixth (index=$33$) eigenmodes  correspond to two of the analytic eigenfunctions that span the system, and the seventh (index=$32$) eigenmode is indeed the product of the fifth and sixth eigenfunctions. Similarly, the ninth and tenth eigenfunctions (index=$31,28$)  also appear to be the polynomial combination of the true eigenfunctions. 

\begin{figure}
\centering
\includegraphics[width=\textwidth]{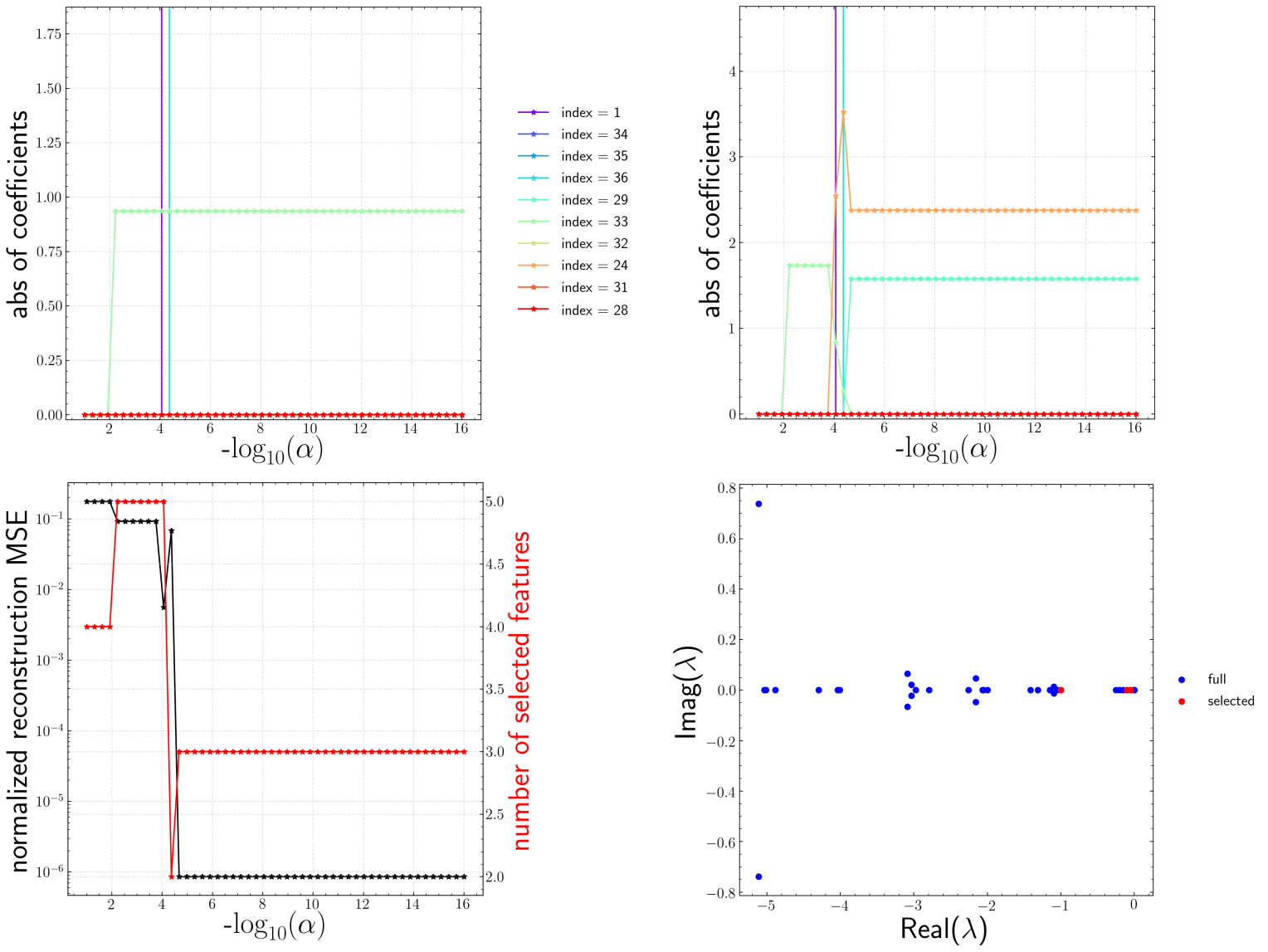}
\caption{Result of multi-task feature learning on top $\hat{L}=10$ accurate eigenmodes from continuous-time EDMD for the 2D fixed point attractor. Top left: ElasticNet path for $x_1$. Top right: ElasticNet path for $x_2$. Bottom left: trends of normalized reconstruction error and number of non-zero terms versus $\alpha$. Bottom right: selected continuous-time eigenvalues.}
\label{fig:2d_fp_edmd_sparse_rec}
\end{figure}

According to \cref{eq:multi_task_enet},  to further remove redundant modes, we perform multi-task feature learning on the $\hat{L}=10$ eigenmodes . The corresponding ElasticNet path is shown in \cref{fig:2d_fp_edmd_sparse_rec}. Note that each $\alpha$ corresponds to a minimizer of  \cref{eq:multi_task_enet}. To choose a proper $\alpha$ so as to find a proper Koopman-invariant subspace, it is advisable to check the trend of the normalized reconstruction error and number of non-zero features. Given the dictionary, for simple problems for which there exists an \emph{exact} Koopman-invariant subspace that also spans system state,  a proper model can be obtained by selecting $\alpha \approx 10^{-6}$ which ends up with only 3 eigenfunctions  as shown in \cref{fig:2d_fp_edmd_sparse_rec}. Moreover, as  is common for EDMD with polynomial basis \citep{williams2014kernel,williams2015data}, a pyramid of eigenvalues appears in \cref{fig:2d_fp_edmd_sparse_rec}. 

\begin{figure}
\centering
\includegraphics[width=\textwidth]{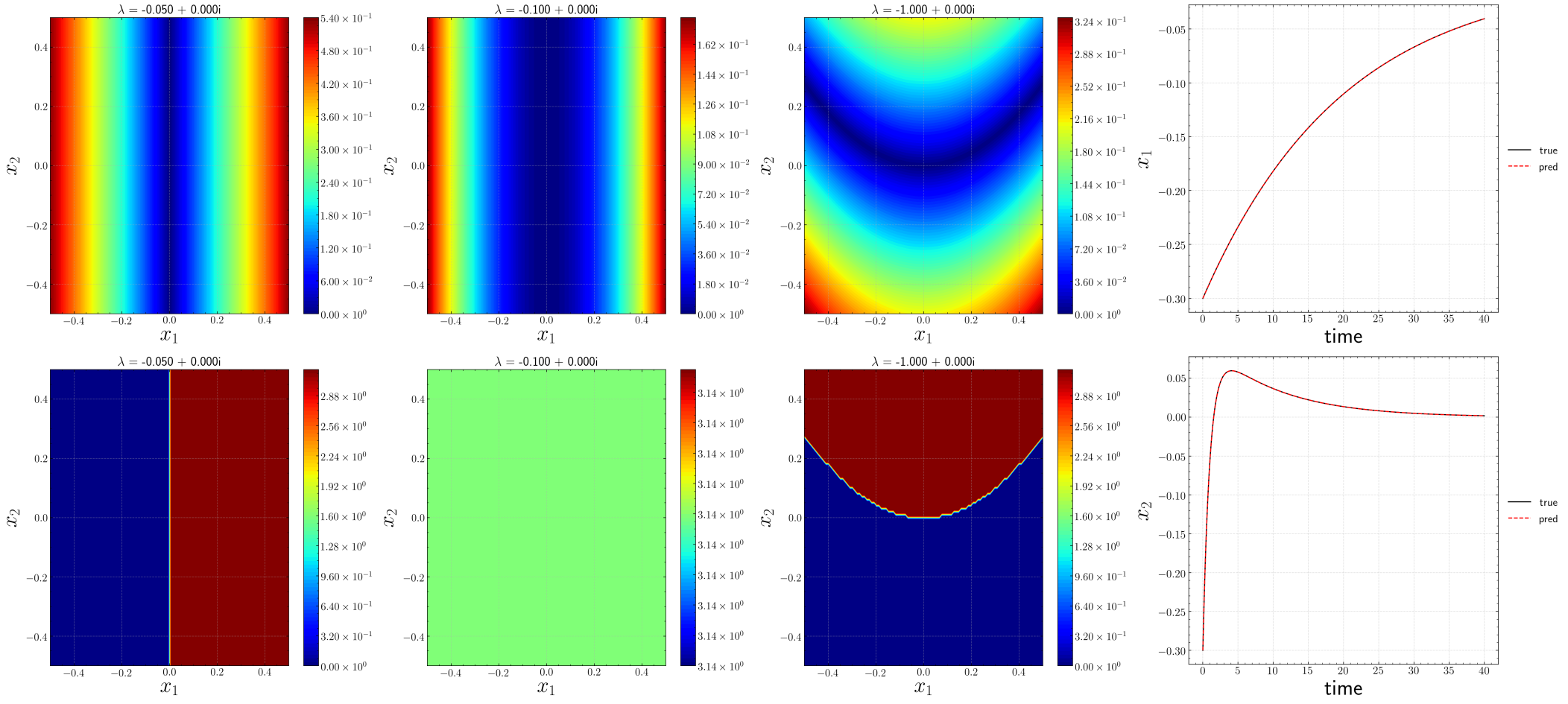}
\caption{Sparsely selected eigenfunctions and eigenvalues from continuous-time EDMD for 2D fixed point attractor  with corresponding prediction on testing data with an unseen initial condition $x_1(0)=x_2(0)=-0.3$. From left to right, the top three figures show contours of magnitude of eigenfunctions, while the bottom three figures are those of phase angle of eigenfunctions. Last column: comparison between prediction and ground truth for an unseen testing trajectory.}
\label{fig:2d_fp_edmd_eigens}
\end{figure}

As shown in \cref{fig:2d_fp_edmd_eigens}, both the identified eigenvalues, and contour of the phase angle and magnitude of selected eigenfunctions from spEDMD match the analytic eigenfunctions given in \cref{eq:2f_analytic} very well. As expected, the prediction on unseen testing data is also excellent. 
Note that the indices of true eigenfunctions $\varphi_1$, $\varphi_2$ and $\varphi_3$ ordered by Kou's criterion in \cref{eq:kou_criterion} are 8, 5 and 6. In this case, all of the true eigenfunctions are missing in the top 3 modes chosen by Kou's criterion. Indeed, the top 3 modes chosen by Kou's criterion have nearly zero eigenvalues.

%

\subsubsection{Results of continuous-time KDMD with mode selection}


\begin{figure}
\centering
\includegraphics[width=\textwidth]{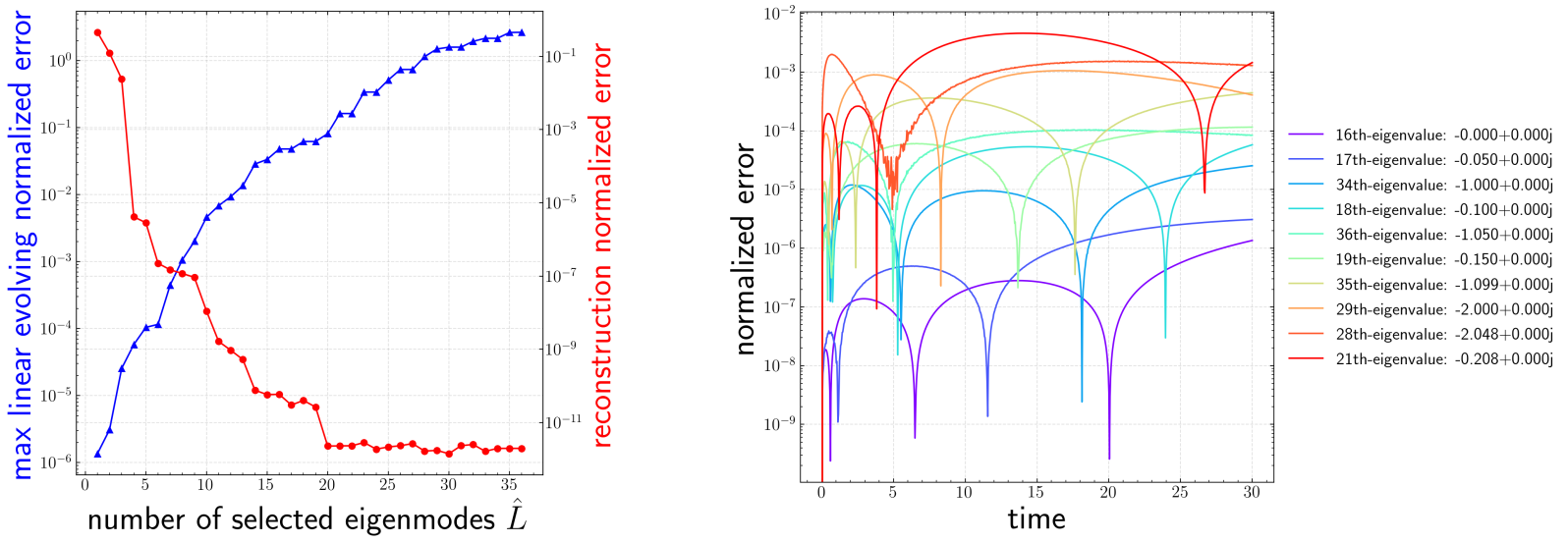}
\caption{Error analysis of 36 eigenmodes from continuous-time KDMD for the 2D fixed point attractor. Left: trends of linearly evolving error $Q$ and reconstruction error $R$. Right: temporal evolution of relative error for top $\hat{L}=10$ accurate eigenmodes.}
\label{fig:2d_fp_kdmd_error_analysis}
\end{figure}

\begin{figure}
\centering
\includegraphics[width=\textwidth]{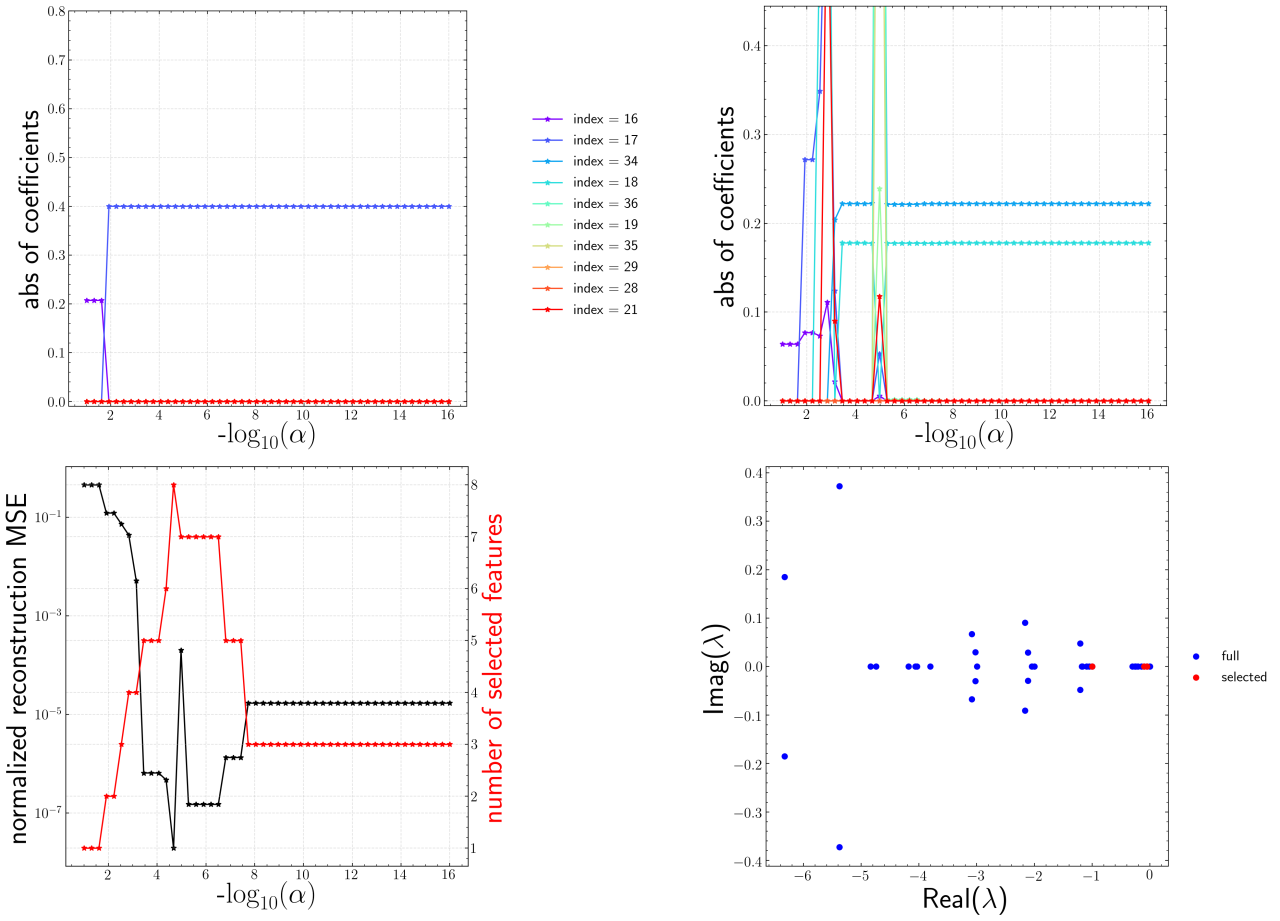}
\caption{Result of multi-task feature learning on top $\hat{L}=10$ accurate eigenmodes from continuous-time KDMD for the 2D fixed point attractor. Top left: ElasticNet path for $x_1$. Top right: ElasticNet path for $x_2$. Bottom left: trends of normalized reconstruction error and number of non-zero terms versus $\alpha$. Bottom right: selected continuous-time eigenvalues.}
\label{fig:2d_fp_kdmd_sparse_rec}
\end{figure}

\begin{figure}
\centering
\includegraphics[width=\textwidth]{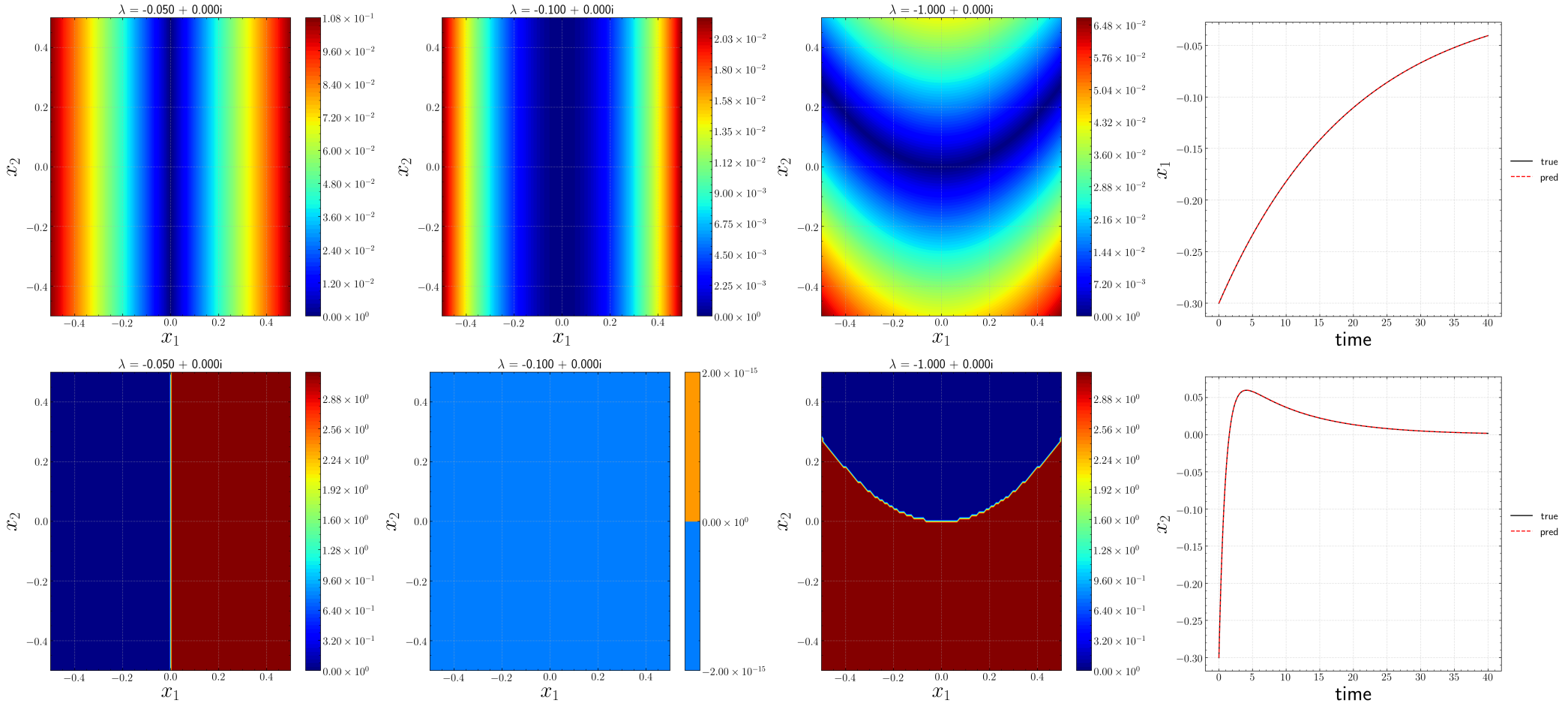}
\caption{Sparsely selected eignefunctions and eigenvalues from continuous-time KDMD for 2D fixed point attractor with corresponding prediction on testing data with an unseen initial condition $x_1(0)=x_2(0)=-0.3$. From left to right, the top three figures show contours of the magnitude of eigenfunctions, while the bottom three figures are those of phase angle of eigenfunctions. Last column: comparison between predictions and ground truth for an unseen testing trajectory.}
\label{fig:2d_fp_kdmd_eigens}
\end{figure}

The  mode selection algorithm presented above can be applied in precisely the same form to KDMD, given a set of eigenfunctions and eigenvalues. 
Error analysis of eigenfunctions is shown in \cref{fig:2d_fp_kdmd_error_analysis}, from which we choose $\hat{L}=10$ as well. \textcolor{cyan}{As before, eigenvalues ordered by maximal deviation from linear evolution are shown in the legend in the right subfigure in \cref{fig:2d_fp_kdmd_error_analysis}.} Again, in the left subfigure in \cref{fig:2d_fp_kdmd_error_analysis}, we observe a sharp decrease in the reconstruction error after the 4 most accurate modes are included.  This is expected, as  the second to fourth most accurate modes are analytically exact from the right subfigure. 
As shown in \cref{fig:2d_fp_kdmd_sparse_rec,fig:2d_fp_kdmd_eigens}, it is confirmed that that both spEDMD and spKDMD arrive at the same analytic eigenfunctions with difference up to a constant factor.
It should be noted that, although polynomials are not analytically in the RKHS~\citep{minh2010some}, good approximations can still be achieved \emph{conditioned on the data} we have, i.e., $x_1,x_2\in[-0.5,0.5]$.
Again, the indices of true eigenfunctions $\varphi_1$ to $\varphi_3$ ordered by Kou's criterion are 8, 2 and 3. Hence, $\varphi_1$ is missing in the top 3 modes chosen by Kou's criterion. Similarly, the first mode chosen by Kou's criterion has a zero eigenvalue.

\subsubsection{Effect of SVD regularization}

SVD truncation is a standard regularization technique in the solution of a potentially ill-conditioned linear system. In the standard EDMD in \cref{eq:discrete_EDMD_general} - for example - $\mathbf{G}$ could be potentially ill-conditioned, leading to   spurious eigenvalues in $\mathbf{K}$. Hence, \citet{williams2014kernel} recommend  SVD truncation  in \cref{eq:G_EDMD} to obtain a robust solution of $\mathbf{K}$. Effectively, it shrinks the number of EDMD/KDMD modes. It has to be recognized, however, that the mode reduction from SVD truncation is not the same as mode selection. Most importantly, one should not confuse \emph{numerical} spuriousness from poor numerical conditioning with \emph{functional} spuriousness from the orthogonal projection error of the Koopman operator~\citep{korda2018convergence}. Indeed, SVD truncation does not always lead to better approximation of a Koopman-invariant subspace. It is rather a linear dimension reduction that optimally preserves the variance in the feature space conditioned on the training data without knowing the linear evolution property of each feature. 

\begin{figure}
\centering
\includegraphics[width=\textwidth]{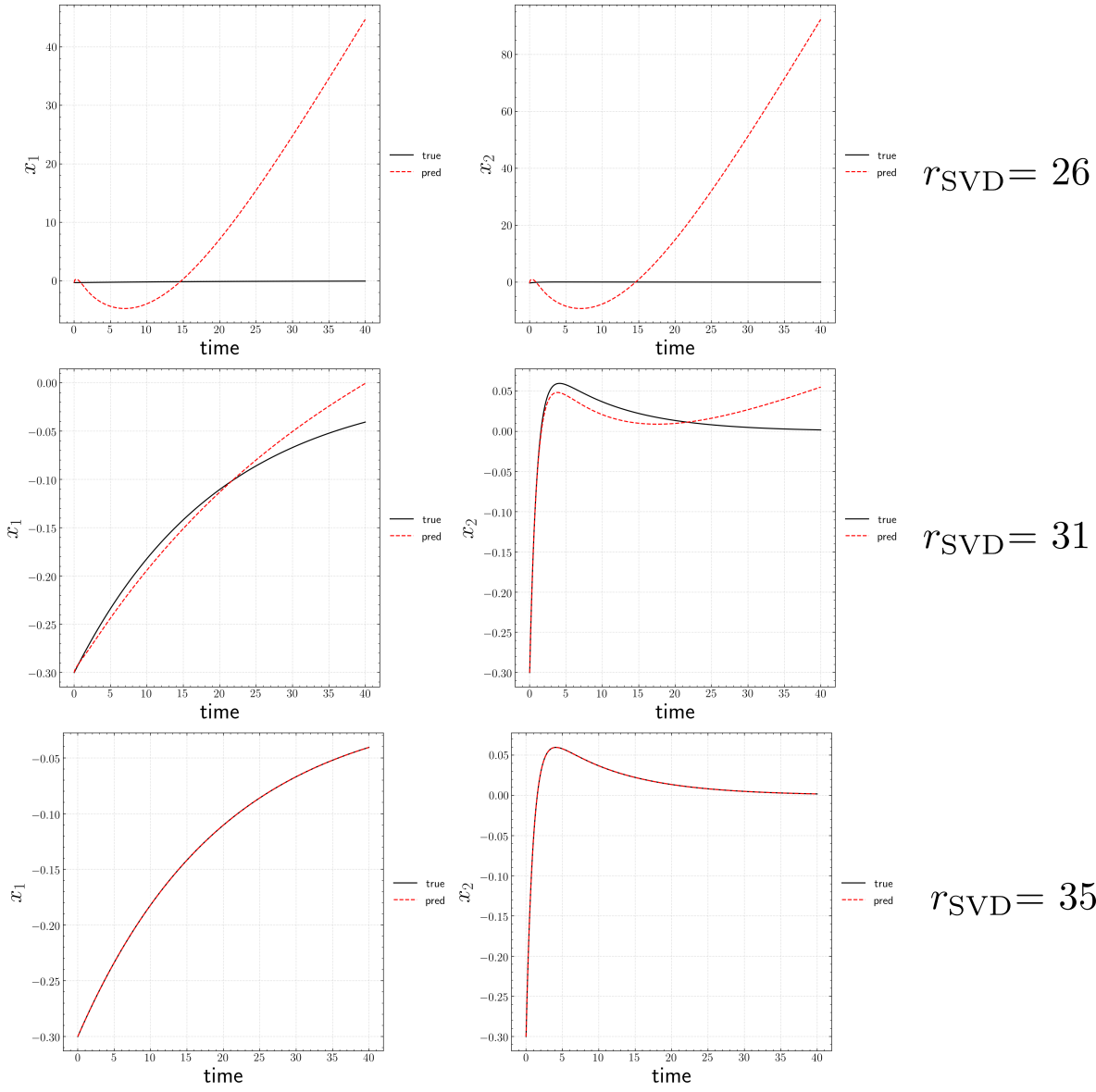}
\caption{Standard EDMD prediction on  unseen trajectory with different SVD truncations for fixed point attractor.}
\label{fig:2d_fp_edmd_svd}
\end{figure}

For demonstration, we take the above fixed point attractor system where we use the same data and standard EDMD algorithm with the same order of Hermite polynomials. The results of prediction on the unseen testing data shown in \cref{fig:2d_fp_edmd_svd} indicate that even though only 3 eigenfunctions (indeed 3 features in the Hermite polynomial) are required, standard EDMD fails to identify the correct eigenfunctions with 26 SVD modes while the results improve with  31 modes retained. The sensitivity of standard EDMD with respect to SVD truncation is likely a result of the use of normalized Hermite polynomials where SVD truncation would lead to a strong preference over the subspace spanned by the higher order polynomials. We did not observe such a sensitivity for KDMD, unless the subspace is truncated below 10 modes.

\subsection{Two-dimensional, transient  flow past a cylinder}

\begin{figure}
\centering
\includegraphics[width=\textwidth]{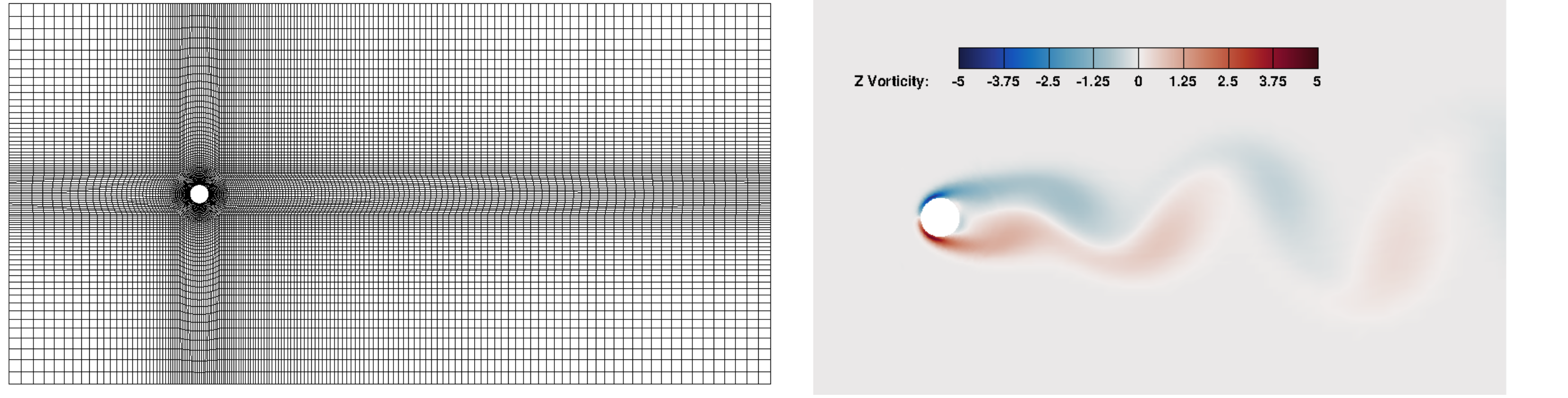}
\caption{Left: illustration of computational mesh for 2D cylinder wake problem (coarsest). Right: contour of vorticity $\omega_z$ for $\Rey=70$ when vortex shedding is fully developed (t = 175).}
\label{fig:cyd_geom}
\end{figure}

\textcolor{cyan}{As a classical example for Koopman analysis in fluid dynamics~\citep{bagheri2013koopman,williams2014kernel,otto2019linearly}}, transient two-dimensional flow past cylinder (\cref{fig:cyd_geom}) is considered at  different Reynolds numbers ($\Rey = {U_{\infty}D}/{\nu}$), where $U_{\infty}=1 $ is the freestream velocity, $D=2$ is the diameter of the cylinder, and $\nu$ is the kinematic viscosity. The two-dimensional incompressible Navier--Stokes equations govern the dynamics with far-field boundary conditions for pressure and velocity and no-slip velocity on the cylinder surface. 
Numerical simulations are performed using the \texttt{icoFoam} solver in \texttt{OpenFOAM}~\citep{jasak2007openfoam} solving the 2D incompressible Navier-Stokes equations. We explore $Re = 70, 100, 130$ by changing the viscosity. 
The pressure field is initialized with i.i.d Gaussian noise $\mathcal{N}(0, 0.3^2)$. The velocity is initialized with a uniform freestream velocity superimposed with i.i.d Gaussian noise $\mathcal{N}(0, 0.3^2)$. It should be noted that the noise is generated on the coarsest mesh shown in \cref{fig:cyd_geom}, and interpolated to  the finer meshes.
Grid convergence with increasing mesh resolution is assessed by comparing the temporal evolution of the drag coefficient $C_D$ and lift coefficient $C_L$. 

Note that the dynamics of a cylinder wake involves \textcolor{cyan}{four} regimes: near-equilibrium linear dynamics, \textcolor{cyan}{nonlinear algebraic interaction between equilibrium and the limit cycle, exponential relaxation rate to the limit cycle}, and periodic limit cycle dynamics~\citep{chen2012variants,bagheri2013koopman}. Instead of considering  data only from each of these regimes separately\textcolor{cyan}{~\citep{chen2012variants,taira2019modal} or with only the last two regimes where exponential linear dynamics is expected~\citep{bagheri2013koopman}}, we  start collecting data immediately after the flow field becomes  unstable,  and stop after the flow field  experiences several limit cycles. \textcolor{cyan}{Note that the regime with algebraic interaction is non-modal~\citep{schmid2007nonmodal}, and therefore cannot be expressed as individual exponential terms~\citep{bagheri2013koopman}. This becomes a challenging problem for DMD~\citep{chen2012variants}}. The sampling time interval is $\Delta t = 0.1 t_{ref}$ where $t_{ref} = D/U_{\infty}$. 

\begin{figure}
\centering
\includegraphics[width=0.6\textwidth]{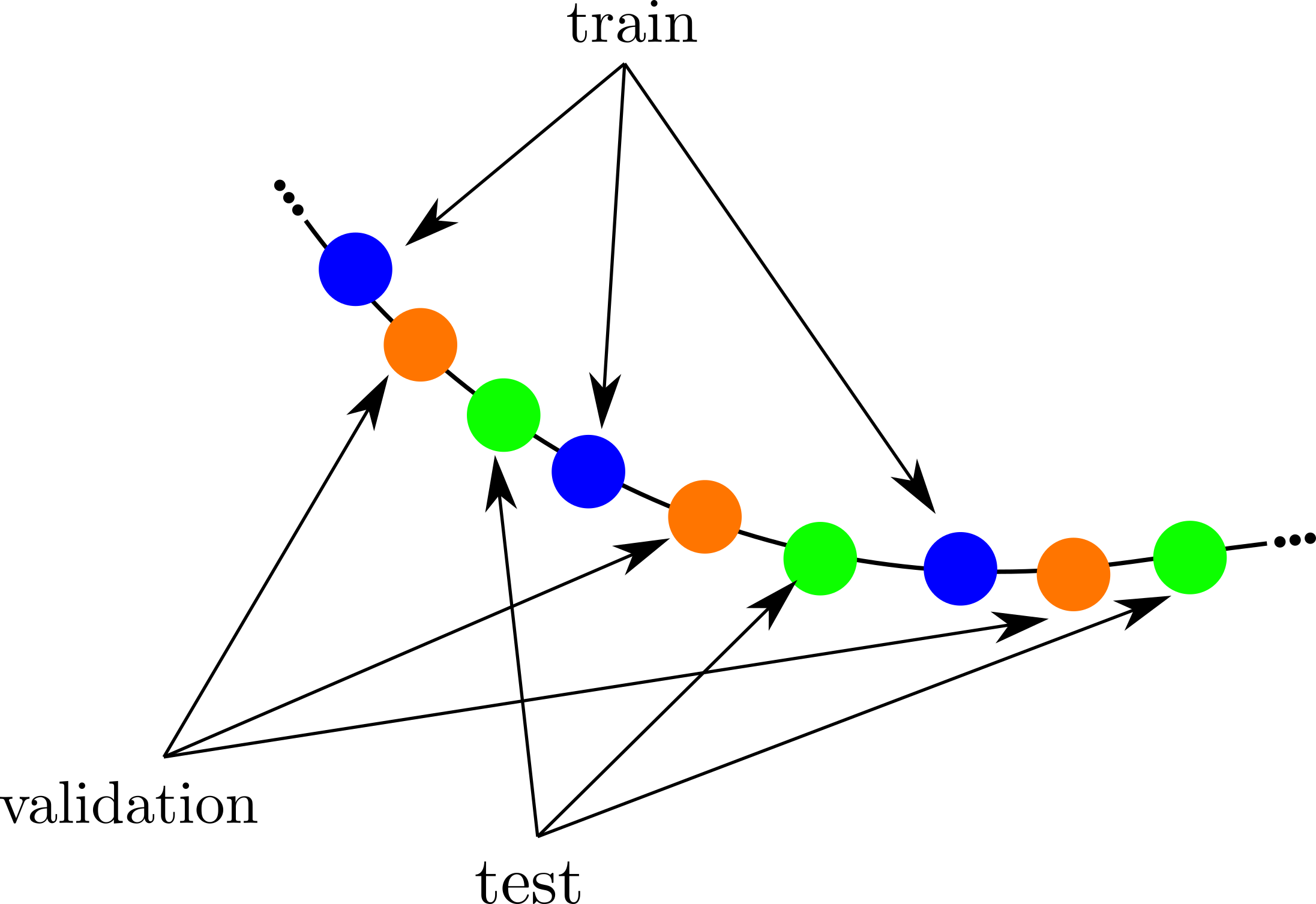}
\caption{Illustration of splitting a uniformly sampled single trajectory in high dimensional phase space into training, validation and testing sets.}
\label{fig:data_triplets}
\end{figure}

For each $\Rey$, 891 snapshots of full velocity field $U$ and $V$ with sampling time interval $\Delta t$ are collected as two matrices of size $N_{\textrm{grid}}\times N_{\textrm{snapshots}}$. Following this, each velocity component is  shifted and scaled (normalized) between $[-1,1]$. 
Since space-filling sampling in any high-dimensional space would be extremely difficult, we split the trajectory into training, validation, and testing data by sampling with strides similar to the ``even-odd" sampling scheme previously proposed by \citet{otto2019linearly}. As illustrated in \cref{fig:data_triplets}, given index $i$, if $ i \mod 3 = 0$, the $i$-th point belongs to training set while $ i \mod 3 = 1$ corresponds to validation, and $ i \mod 3 = 2$ for testing data. Consequently, the time interval in the training, testing, validation trajectory is tripled as $3\Delta t = 0.3 t_{ref}$.
Thus, training, validation, and testing data are split into 297 snapshots each. Finally, we stack data matrices along the first axis corresponding to the number of grid points, and perform a distributed SVD described in \cref{sec:imple}. For all three cases, the top 20 POD modes are retained, corresponding to 99\% of kinetic energy. Next, we apply our algorithm to discrete-time KDMD \textcolor{blue}{with isotropic Gaussian kernel} on this reduced-order nonlinear system. We choose the hyperparameters $\sigma = 3$ and $r=180$. Further details are given in \cref{apdx:hyper_cyd}. 

\subsubsection{Results of discrete-time KDMD with mode selection}
\label{sec:kdmd_cyd}

\begin{figure}
\centering
\includegraphics[width=\textwidth]{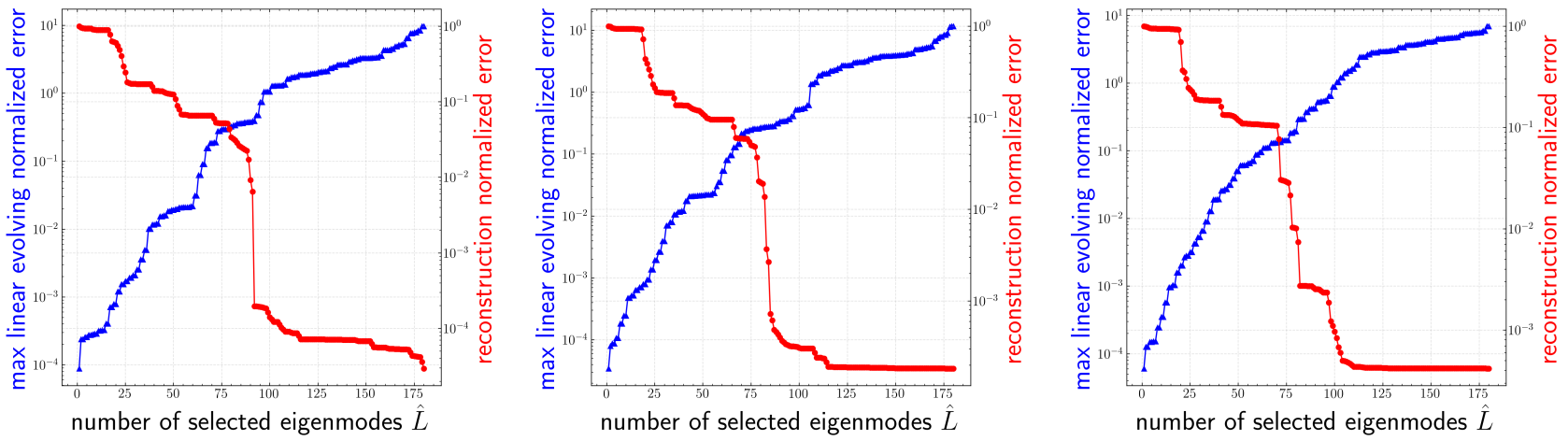}
\caption{Trend of linear evolution error $Q$ and reconstruction error $R$ from discrete-time KDMD for the 20D cylinder wake flow case. Left: $\Rey=70$. Center: $\Rey =100$. Right: $\Rey=130$.}
\label{fig:20d_cyd_qr}
\end{figure}

For all three Reynolds numbers, a posteriori error analysis is shown in \cref{fig:20d_cyd_qr}. A good choice of the number of accurate modes $\hat{L}$ retained for reconstruction is around 60 since the corresponding maximal deviation from linear evolution is still around  5\% while the reconstruction error reaches a plateau after $\hat{L}>60$.

\begin{figure}
\centering
\includegraphics[width=\textwidth]{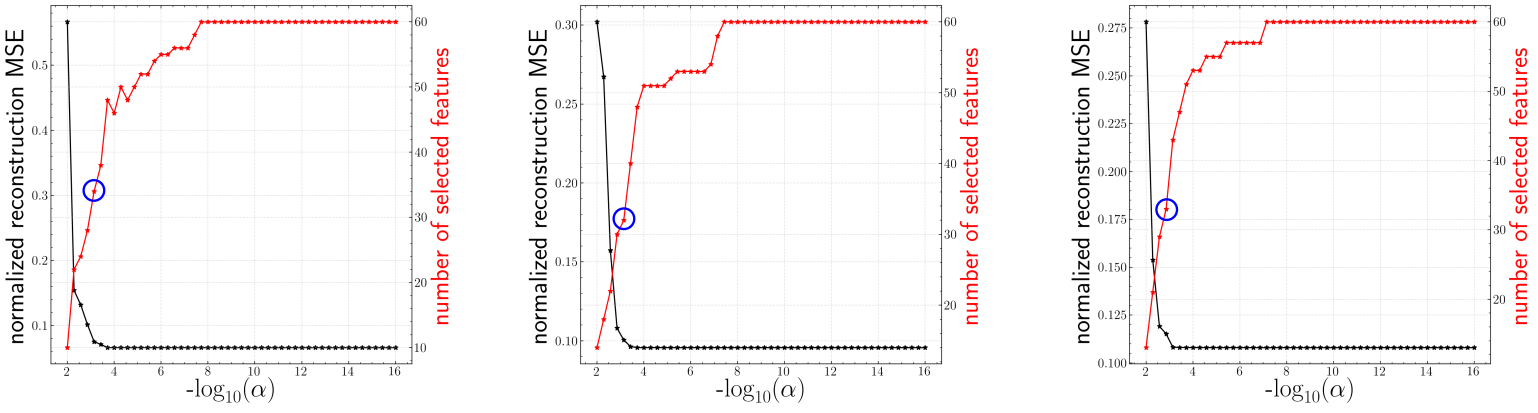}
\caption{Variation of reconstruction error $R$ and number of non-zero terms for the 20D cylinder wake flow. Left: $\Rey=70$. Center: $\Rey =100$. Right: $\Rey=130$. Blue circle corresponds to selected $\alpha$.} 
\label{fig:20d_cyd_mse}
\end{figure}

After the mode selection on validation data, a $\alpha$-family of solutions is obtained with corresponding reconstruction error and the number of non-zeros terms as shown in \cref{fig:20d_cyd_mse}. Note that the  chosen solution is highlighted as blue circles. As shown in  \cref{tab:modes_selection_cyd}, nearly half of the accurate KDMD eigenmodes identified are removed with the proposed sparse feature selection. Note that for all three cases, the number of selected modes (around 32 to 34) is still larger than that required in neural network models (around 10) \citep{pan2019physics,otto2019linearly}. This is because the subspace spanned by KDMD/EDMD relies on a pre-determined dictionary rather than being data-adaptive like neural network models. Nevertheless, due to the additional expressiveness from non-linearity, we will see  in \cref{sec:compare} that spKDMD performs significantly better than DMD~\citep{schmid2010dynamic} and spDMD~\citep{jovanovic2014sparsity}, while enjoying the property of convex optimization at a much lower computational cost than the inherently non-convex and computationally intensive neural network counterparts.

\begin{table}
 \begin{center}
\def~{\hphantom{0}}
 \begin{tabular}{lccc}
  & $Re=70$ & $Re=100$ &  $Re=130$ \\ [3pt] 
  $\alpha_{\textrm{select}}$ & $7.19\times 10^{-4}$ & $7.19\times 10^{-4}$ & $1.39\times 10^{-3}$\\
  number of selected modes & 34 & 32 & 33 \\
  number of total modes & 297 & 297 & 297 \\
  normalized reconstruction error $R$ & 0.075 & 0.105& 0.113 \\ 
 \end{tabular}
\caption{Summary of mode selection for discrete-time KDMD on 20D cylinder wake flow.}
\label{tab:modes_selection_cyd}
 \end{center}
\end{table}

\begin{figure}
\centering
\includegraphics[width=\textwidth]{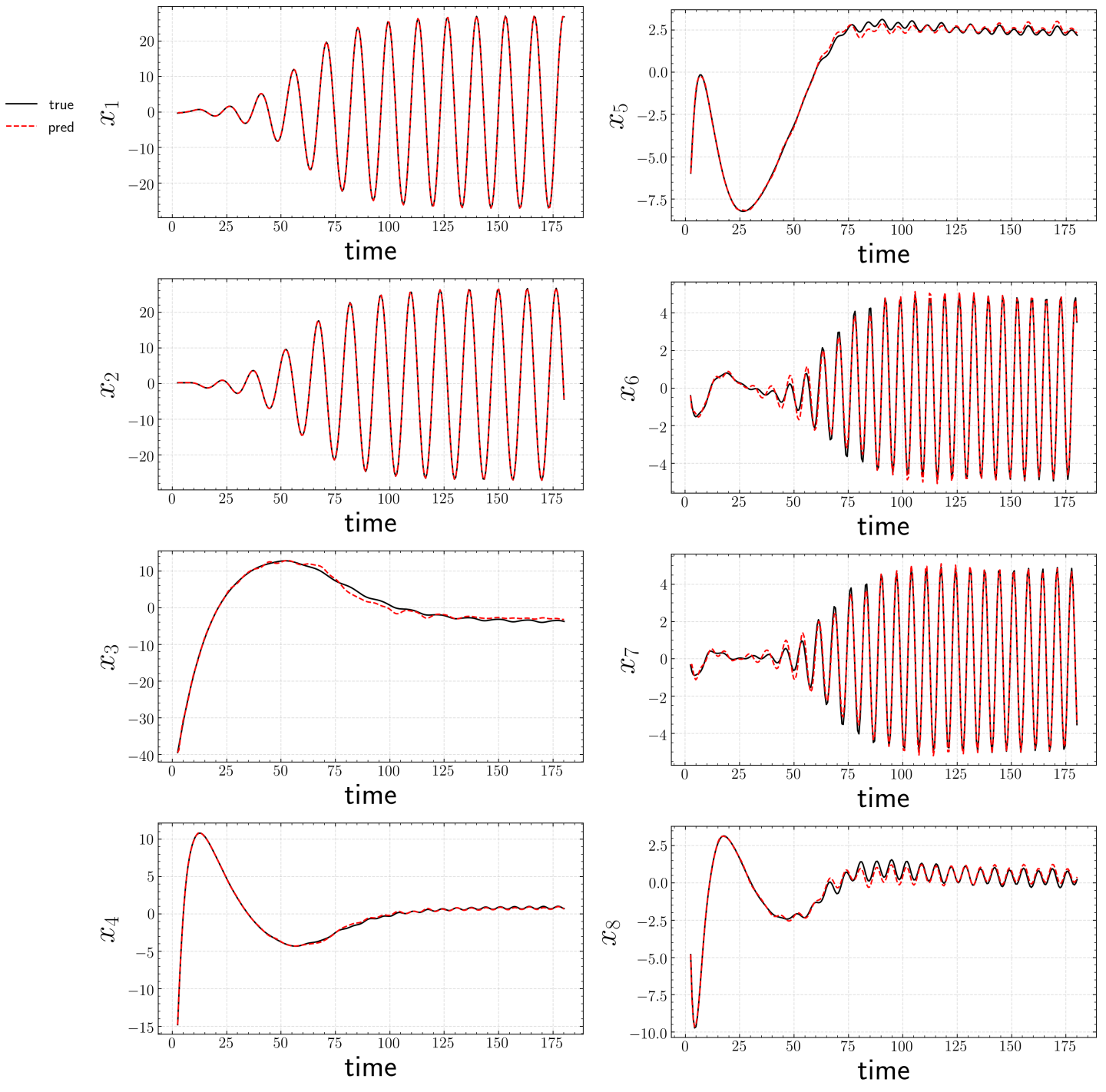}
\caption{A posteriori prediction of testing trajectory for $\Rey = 70$ in terms of top 8 POD coefficients with sparsity-promoting KDMD. } 
\label{fig:20d_cyd_pod_70}
\end{figure}
\begin{figure}
\centering
\includegraphics[width=\textwidth]{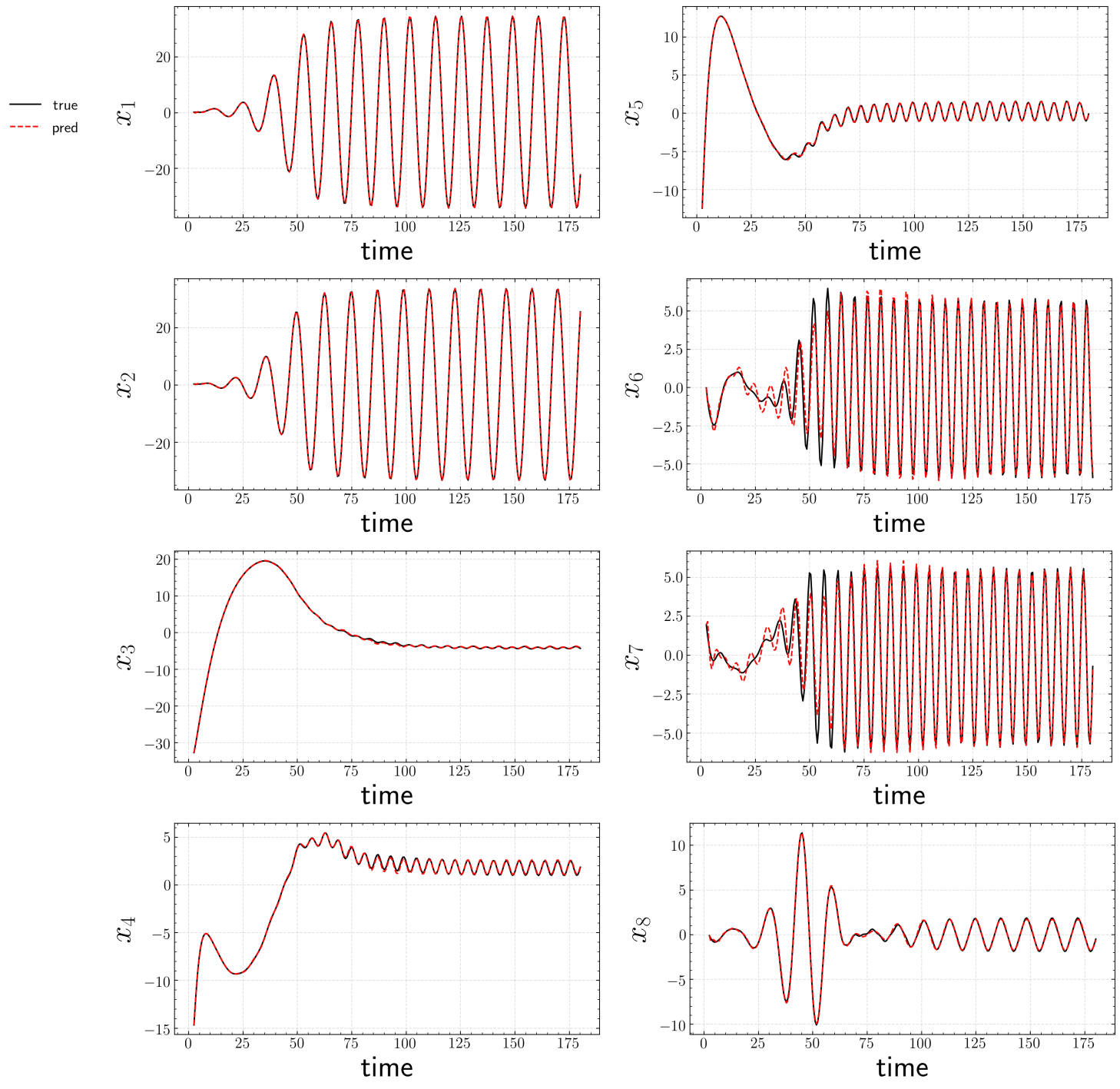}
\caption{A posteriori prediction of testing trajectory for $\Rey = 100$ in terms of top 8 POD coefficients with sparsity-promoting KDMD.} 
\label{fig:20d_cyd_pod_100}
\end{figure}
\begin{figure}
\centering
\includegraphics[width=\textwidth]{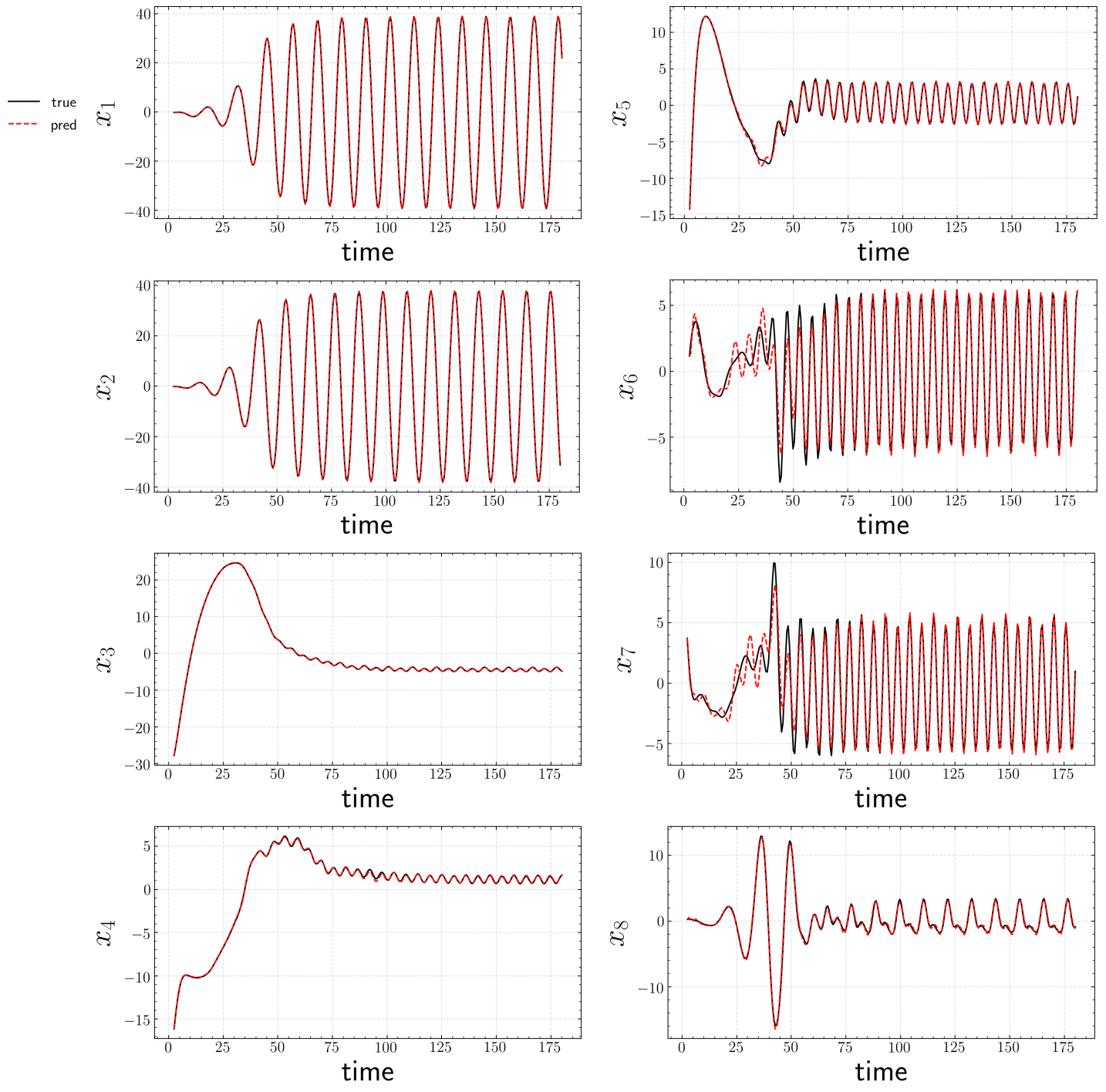}
\caption{A posteriori prediction of testing trajectory for $\Rey = 130$ in terms of top 8 POD coefficients with sparsity-promoting KDMD. } 
\label{fig:20d_cyd_pod_130}
\end{figure}

\begin{figure}
\centering
\includegraphics[width=1\textwidth]{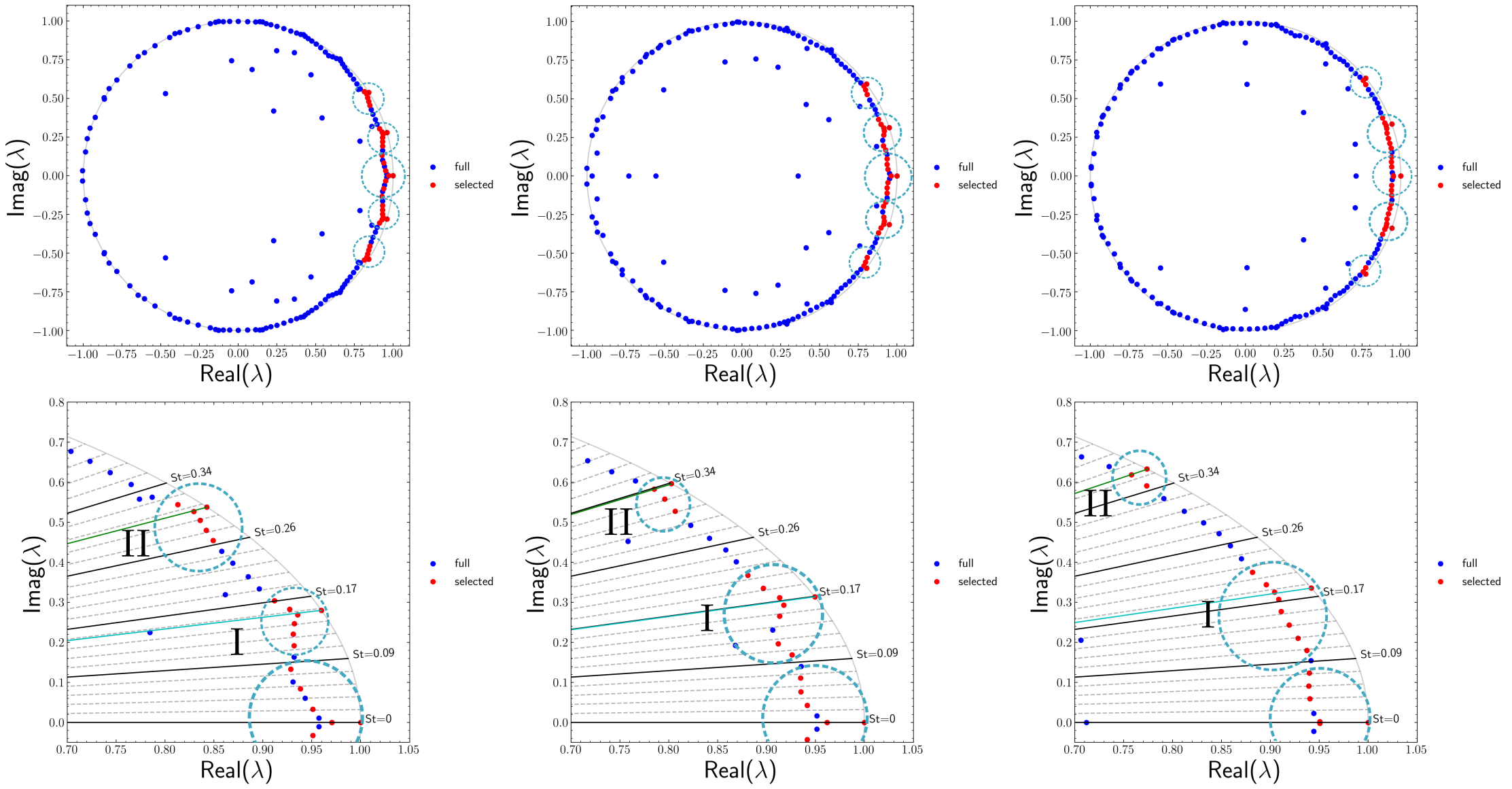}
\caption{Discrete-time eigenvalue distribution of full KDMD and spKDMD. Left: $\Rey=70$. Center: $\Rey =100$. Right: $\Rey=130$. Blue dot: full KDMD eigenvalues. Red dot: spKDMD eigenvalues. Bottom row: zoomed. $\textrm{I}$ and $\textrm{II}$ correspond to two types of eigenvalue clusters of distinct frequencies, with each of them enclosed by cyan dashed circles. Green/cyan solid line correspond to $St_D$/$St_{L}$.} 
\label{fig:20d_cyd_eigenvals}
\end{figure}

The predictions of the top 8 POD coefficients (denoted as $x_1$ to $x_8$) on testing data are displayed in \cref{fig:20d_cyd_pod_70,fig:20d_cyd_pod_100,fig:20d_cyd_pod_130}. The results match very well with ground truth for all three cases. 
\Cref{fig:20d_cyd_eigenvals} shows that there appear to be five clusters of selected eigenvalues while most of the modes are removed by the proposed algorithm. \textcolor{cyan}{Similar observations were also made in \citep{bagheri2013koopman} when DMD is applied to the \emph{full} transient dynamics.} This pattern consisting of a stable eigenvalue \emph{on} the unit circle surrounded by several decaying eigenvalues is observed for all clusters. The stable eigenvalue contributes to limit-cycle behavior, while the decaying eigenvalues account for the transient phenomenon. Due to symmetry, only eigenvalues in the first quadrant are shown in the bottom row of \cref{fig:20d_cyd_eigenvals}. It is observed that the frequency associated with the type-$\textrm{II}$ cluster is approximately twice that of type-$\textrm{I}$. \textcolor{cyan}{This is in good agreement with previous analytical results from the weakly non-linear theory~\citep{bagheri2013koopman}.} The frequency $f$ is normalized  as $\textrm{St} = fD/U_{\infty}$, where St is the Strouhal number. 

\begin{figure}
    \centering
    \includegraphics[width=0.4\textwidth]{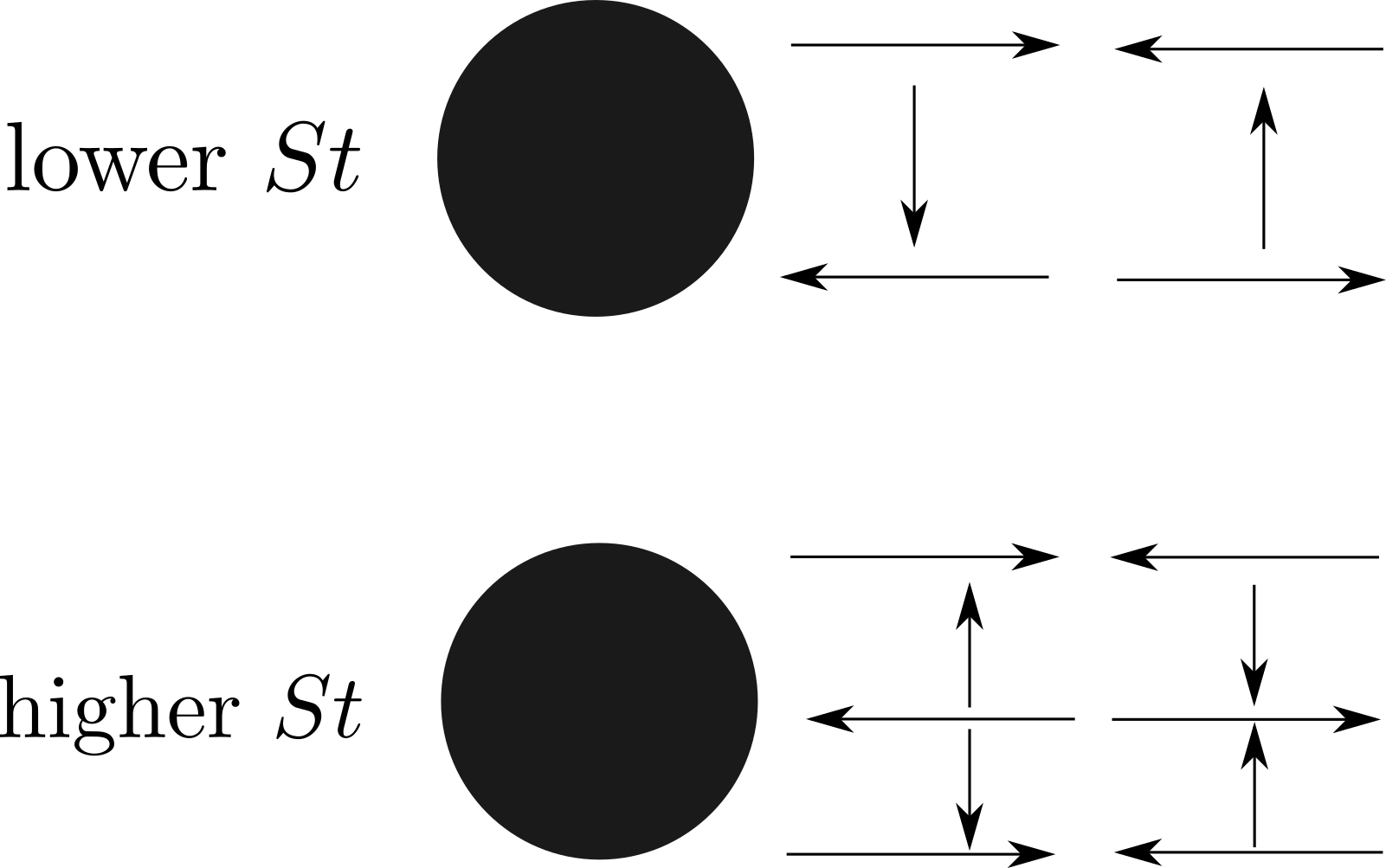}
    \caption{Illustration of the structure of velocity field for the lower (top) and  higher frequency (bottom) Koopman modes. The arrow roughly indicates the velocity direction. }
    \label{fig:reviewer_1_png}
\end{figure}

\begin{figure}
\centering
\includegraphics[width=1\textwidth]{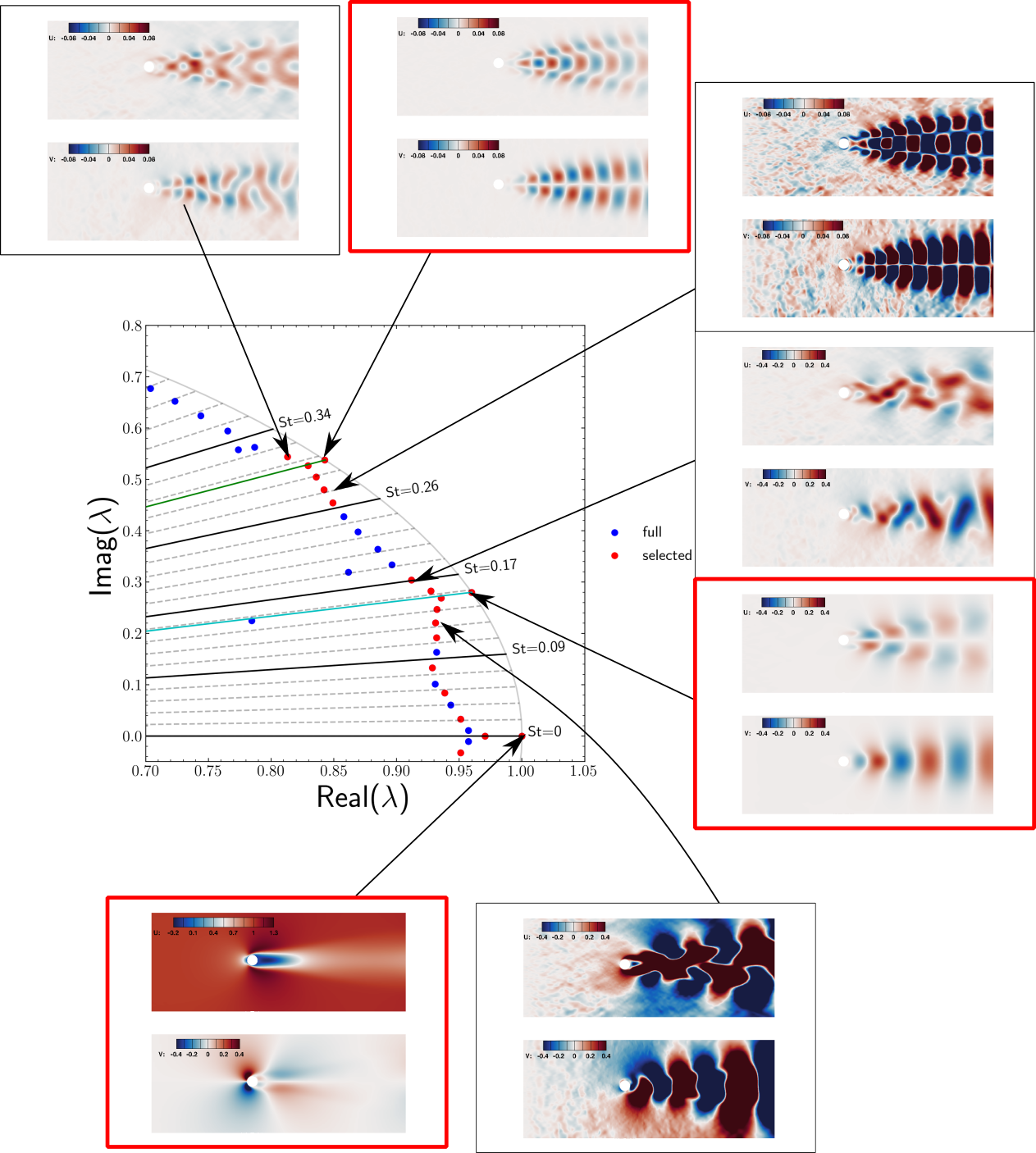}
\caption{Contours of Koopman modes of $\Rey=70$ cylinder wake flow at $t=0$. Red squares indicate stable modes.} 
\label{fig:20d_cyd_re70_modes}
\end{figure}

Recall that in the laminar parallel shedding region ($47 < \Rey < 180$), \textcolor{cyan}{the characteristic Strouhal number $St$ scales with $-1/\sqrt{\Rey}$}~\citep{fey1998new}. Therefore, it is expected that $St$ of both types tend toward higher frequency as $Re$ increases from 70 to 130. Further, it is interesting to note that the corresponding Strouhal numbers for lift and drag  when the system is on the limit-cycle - $St_L$ and $St_D$\footnote{we observe that each of lift and drag coefficients exhibits only one frequency at limit-cycle regime for the range of $\Rey$ studied in this work.} - coincide with the stable frequency of type-$\textrm{I}$ and $\textrm{II}$ respectively as indicated in \cref{fig:20d_cyd_eigenvals}. \textcolor{brown}{This is due to the anti-symmetrical/symmetrical structure of the velocity field of type-$\textrm{I}$/$\textrm{II}$ Koopman mode respectively as can be inferred from  \cref{fig:20d_cyd_re70_modes}.  A schematic is is also shown in \cref{fig:reviewer_1_png}.}

\textcolor{brown}{ The higher frequency mode is symmetrical (along the freestream direction) in $U$ and anti-symmetrical in $V$. As a consequence, this only contributes to the oscillation of drag. The lower frequency mode is anti-symmetrical in $U$ and symmetrical in $V$, and only contributes to the oscillation of lift.
}  Thus,  the fluctuation in the lift mostly results from the stable mode in type-$\textrm{I}$, while that for drag results from the stable mode in type-$\textrm{II}$ with twice the frequency.

\begin{figure}
\centering
\includegraphics[width=1\textwidth]{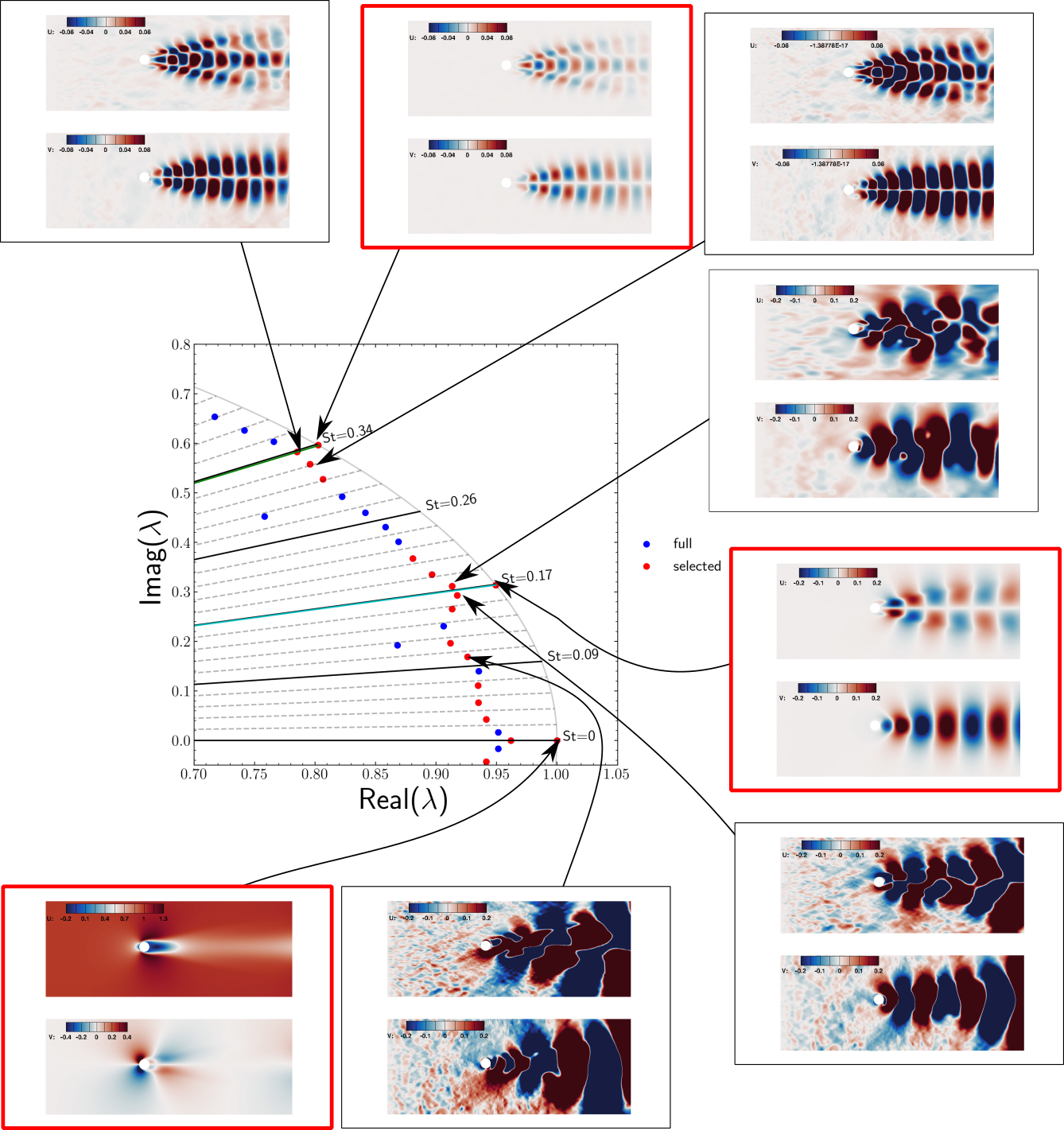}
\caption{Contours of Koopman modes of $\Rey=100$ cylinder wake flow at $t=0$. Red squares indicate stable modes.} 
\label{fig:20d_cyd_re100_modes}
\end{figure}

\begin{figure}
\centering
\includegraphics[width=1\textwidth]{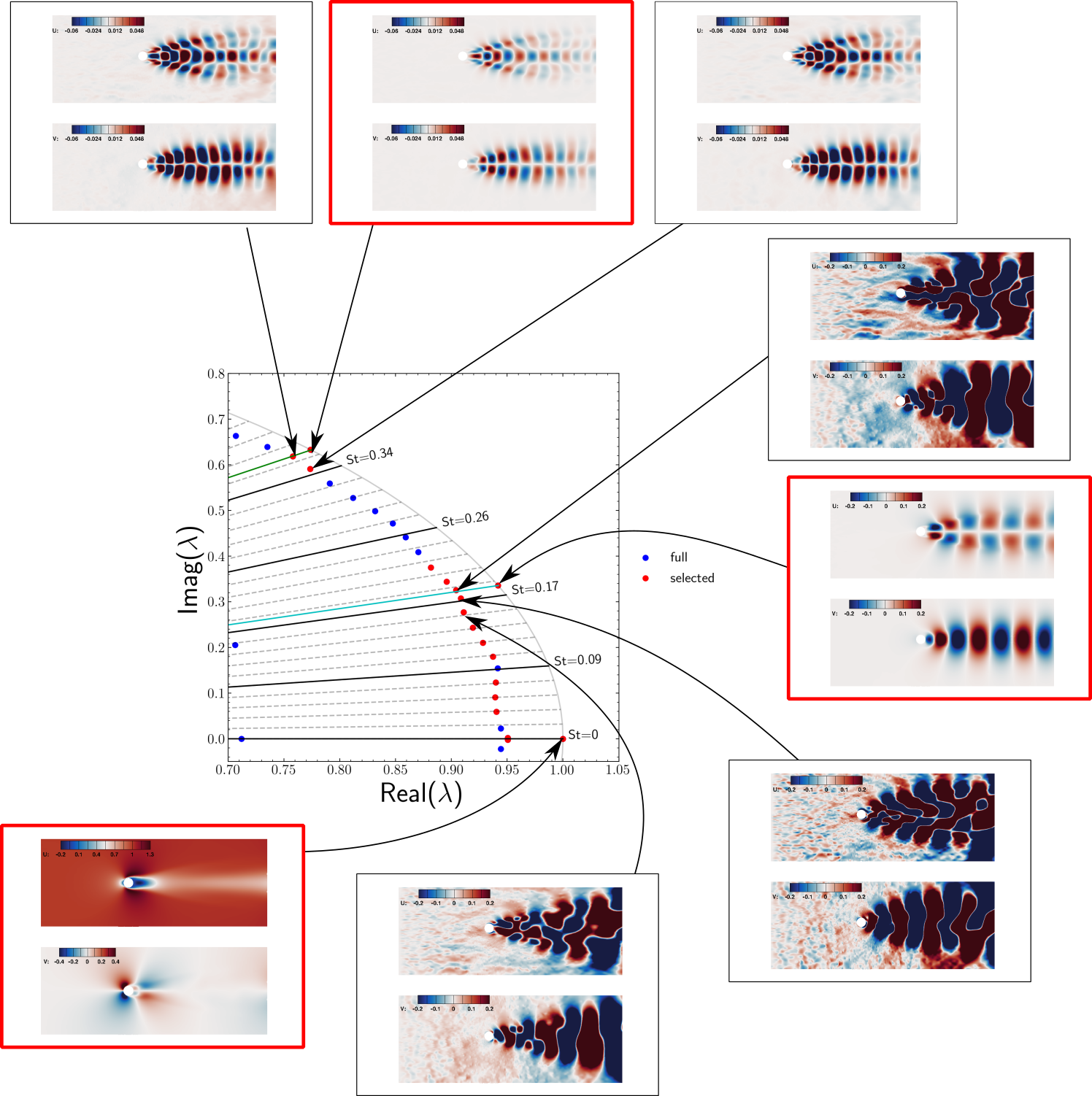}
\caption{Contours of Koopman modes  of $\Rey=130$ cylinder wake flow at $t=0$. Red squares indicate stable modes.} 
\label{fig:20d_cyd_re130_modes}
\end{figure}

\begin{figure}
\centering
\includegraphics[width=1\textwidth]{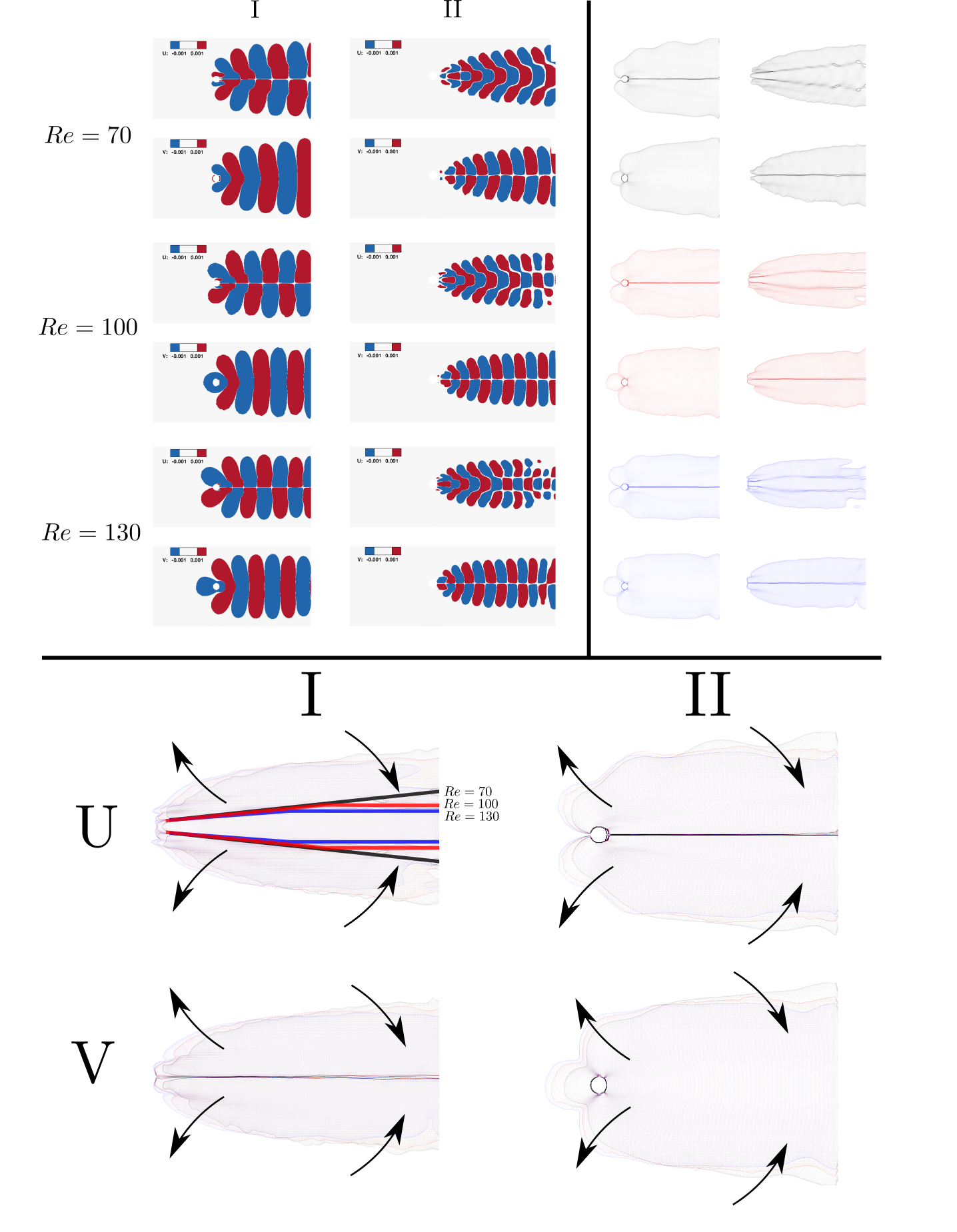}
\caption{Top left: contribution of stable Koopman modes corresponding to type-I and type-II  cluster for $\Rey = 70,100,130$ at $t=0$ visualized with threshold 0.001. Top right: time-averaged iso-contour of top left plot. Bottom: tendency of ``envelope" of type-I and II modes as $\Rey$ increases. Separation lines in $U$ component of type-I are drawn for $\Rey = 70$ (black), $\Rey = 100$ (red) and $\Rey = 130$ (blue).} 
\label{fig:20d_cyd_phase_plot}
\end{figure}

Finally, several representative Koopman modes from spKDMD for three $\Rey$ are shown in \textcolor{pink}{\cref{fig:20d_cyd_re70_modes,fig:20d_cyd_re100_modes,fig:20d_cyd_re130_modes}}. For a better comparison of mode shapes, contributions from the stable modes of type-I and II with a threshold of 0.001 at $t=0$ is displayed in the top left of \cref{fig:20d_cyd_phase_plot}. To remove the effect of time, the ``envelope" of the mode shape, i.e., time average of the iso-contours is shown in the top right of \cref{fig:20d_cyd_phase_plot}. From these results, we observe the following interesting phenomena:
\begin{itemize}
    \item the minimal dimension of the Koopman-invariant subspace that approximately captures the limit cycle attractor for all three $\Rey$ that fall into laminar vortex shedding regime \citep{white2006viscous} is five, \textcolor{cyan}{which is consistent with previous multi-scale expansion analysis near the critical $\Rey$~\citep{bagheri2013koopman}}.
    \item the lobes of stable Koopman modes in the type-I cluster show an approximately 50\% larger width than those in a type-II cluster. 
    \item similarity\textcolor{brown}{/colinearity} among Koopman modes within each cluster is observed. \textcolor{cyan}{Such a finding is previously reported in the theoretical analysis by \cite{bagheri2013koopman}. A similarity in spatial structure exists among the Koopman modes belonging to the same family, even though the structures are clearly  phase lagged.}
    \item as $\Rey$ increases from 70 to 130, mode shapes  flatten downstream while expand upstream. 
    \item at $\Rey=70$, the shear layer in the stable Koopman modes continues to grow within the computational domain. However, at $\Rey = 100, 130$, the shear layer stops growing after a certain distance that is negatively correlated with $\Rey$. 
\end{itemize}

\subsubsection{Net contribution of clustered Koopman modes}

\textcolor{brown}{By examining \cref{fig:20d_cyd_re70_modes,fig:20d_cyd_re100_modes,fig:20d_cyd_re130_modes}, we observe that the colinearity of the spatial structures among each cluster can cause cancellations. A famous example of such non-oscillatory cancellation is the ``shift mode'' defined by~\cite{noack2003hierarchy}, in conjunction with  two oscillating modes. As the ``shift mode'' decays, the stationary component of the flow transitions from the unstable equilibrium to the time-averaged mean. The existence of such non-oscillatory decaying Koopman modes is also confirmed by weakly non-linear analysis~\citep{bagheri2013koopman}. }
\textcolor{brown}{Interestingly, our algorithm is able to identify not only the non-oscillatory cancellation (from the ``shift mode'') but also oscillatory cancellations from two clusters with distinct frequencies. Such cancellations elegantly explain why no unstable Koopman eigenvalues appear in this flow given the co-existence of an attracting limit cycle and an unstable equilibrium. These modes could be understood as ``oscillatory shift modes'', as a generalization of the model proposed by~\cite{noack2003hierarchy}.}

\begin{figure}
\centering
\includegraphics[width=1\textwidth]{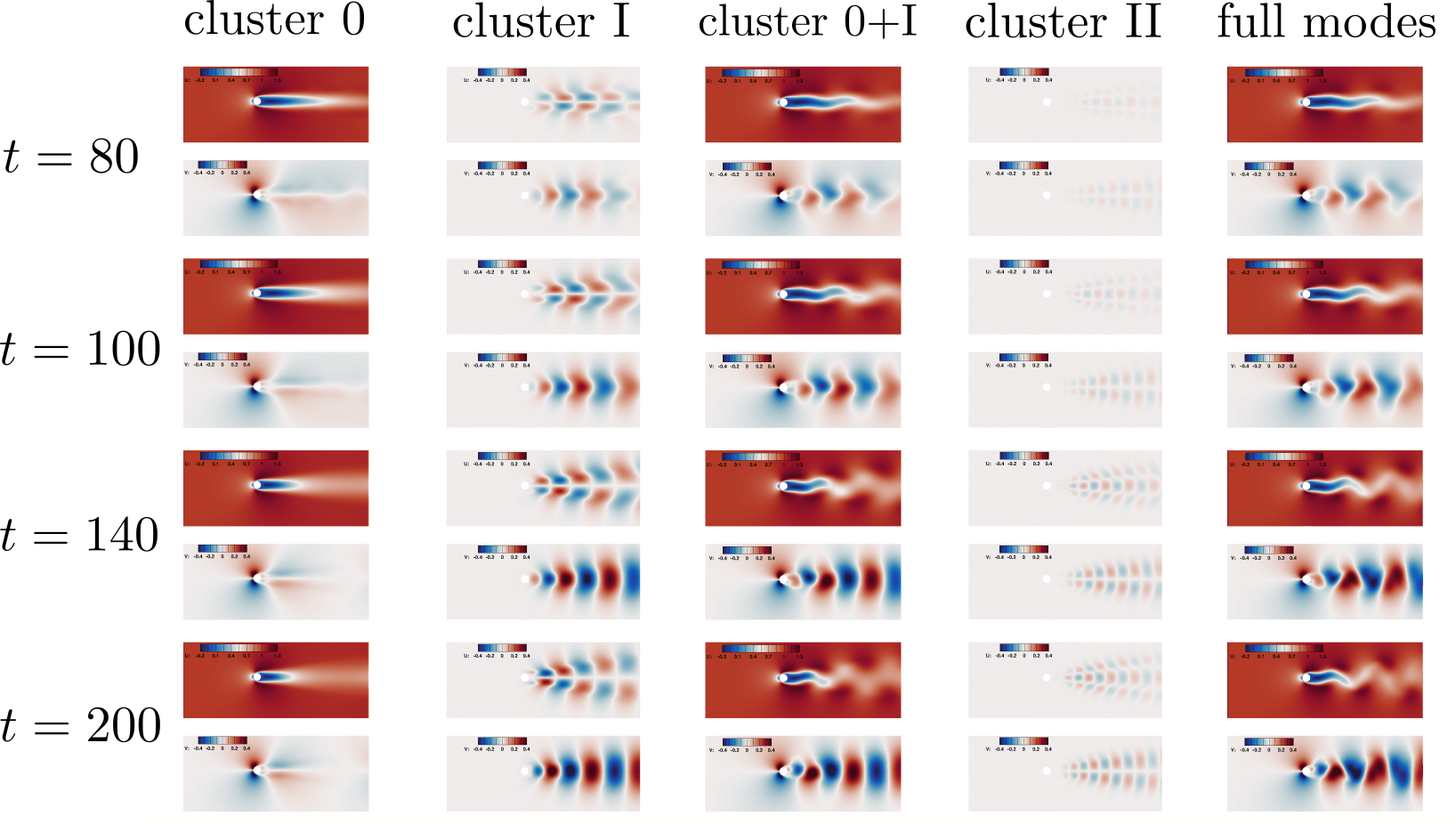}
\caption{Contribution of Koopman modes  at cluster level in the transient regime \textcolor{green}{of $Re=70$ case}. ``cluster 0" denotes  the cluster near the real axis in \cref{fig:20d_cyd_eigenvals}. ``cluster $I$" /``cluster $II$" takes the effect of mirror cluster in fourth quadrant into account. ``full modes" denotes the aggregated contribution of Koopman modes.} 
\label{fig:20d_cyd_cluster_modes}
\end{figure}

\textcolor{brown}{Since modes within each cluster appear to be colinear to each other, it is intriguing to investigate the net contribution from each cluster.}  For the above $\Rey = 70$ case, effects from different clusters for different time in the transient regime are shown in \cref{fig:20d_cyd_cluster_modes}. There are several interesting observations throughout the transient period from $t=80$ to $t=200$:
\begin{itemize}
    \item The net contribution from ``cluster 0" does not exhibit strong oscillations. For the contribution from ``cluster 0 ", the $U$ component shows a decrease in the length of the reverse flow region behind the cylinder with an increase in the low speed wake layer thickness while the $V$ component remains unchanged. \textcolor{cyan}{This is similar to the effect of ``shift mode'' which also characterizes the decay of recirculation behind the cylinder.}
    \item Initially at $t=80$, the net contribution from ``cluster $I$" is rather weak \textcolor{cyan}{primarily due to the ``cancellation'' from the lagged phase. Further,} the development of vortex shedding downstream from the top and bottom surfaces of the cylinder is nearly parallel. This corresponds to the initial wiggling of the low speed wake flow downstream. As time increases, the pattern of corresponding vortices develops away from the center line.
    \item Although the net contribution from ``cluster 0+$I$" captures most of the flow features from ``full modes" throughout the transient regime, with increasing time, the net contribution from ``cluster $II$" becomes more important and contributes to the strength of vortex shedding downstream. 
\end{itemize}

\subsubsection{Comparison with DMD and sparsity-promoting DMD at $\Rey=70$}
\label{sec:compare}

To confirm the advantage of sparsity-promoting KDMD over DMD~\citep{schmid2010dynamic} and spDMD~\citep{jovanovic2014sparsity}, we compare the following three models on the unsteady cylinder wake flow at $\Rey=70$:
\begin{enumerate}
\item sparsity-promoting KDMD on top 20 POD coefficients with 34 modes selected,
\item DMD on top 20 POD coefficients,
\item spDMD~\footnote{\textcolor{cyan}{We used the original Matlab code from http://people.ece.umn.edu/users/mihailo/}} on top 200 POD coefficients with  $\alpha$ chosen carefully such that only 34 modes are selected.
\end{enumerate}

Note that DMD with top 200 POD coefficients, i.e., $r=200$ in SVD-DMD~\citep{schmid2010dynamic}, contains 10 times stable/decaying harmonics as DMD on the top 20 modes. Hence, it is not surprising to expect that the corresponding prediction of the evolution of the top 20 POD coefficients to be very good (not shown for clarity). However, to make a fair comparison against sparsity-promoting KDMD, we consider spDMD~\citep{jovanovic2014sparsity} on the top 200 POD coefficients with 34 modes selected\footnote{Although there are 200 POD coefficients used for spDMD and 20 for KDMD, it is not an  unfair comparison given that the same number of spatial modes are selected. Further,  these are reduced order models of the same the full order dynamical system.}. 

\begin{figure}
\centering
\includegraphics[width=1\textwidth]{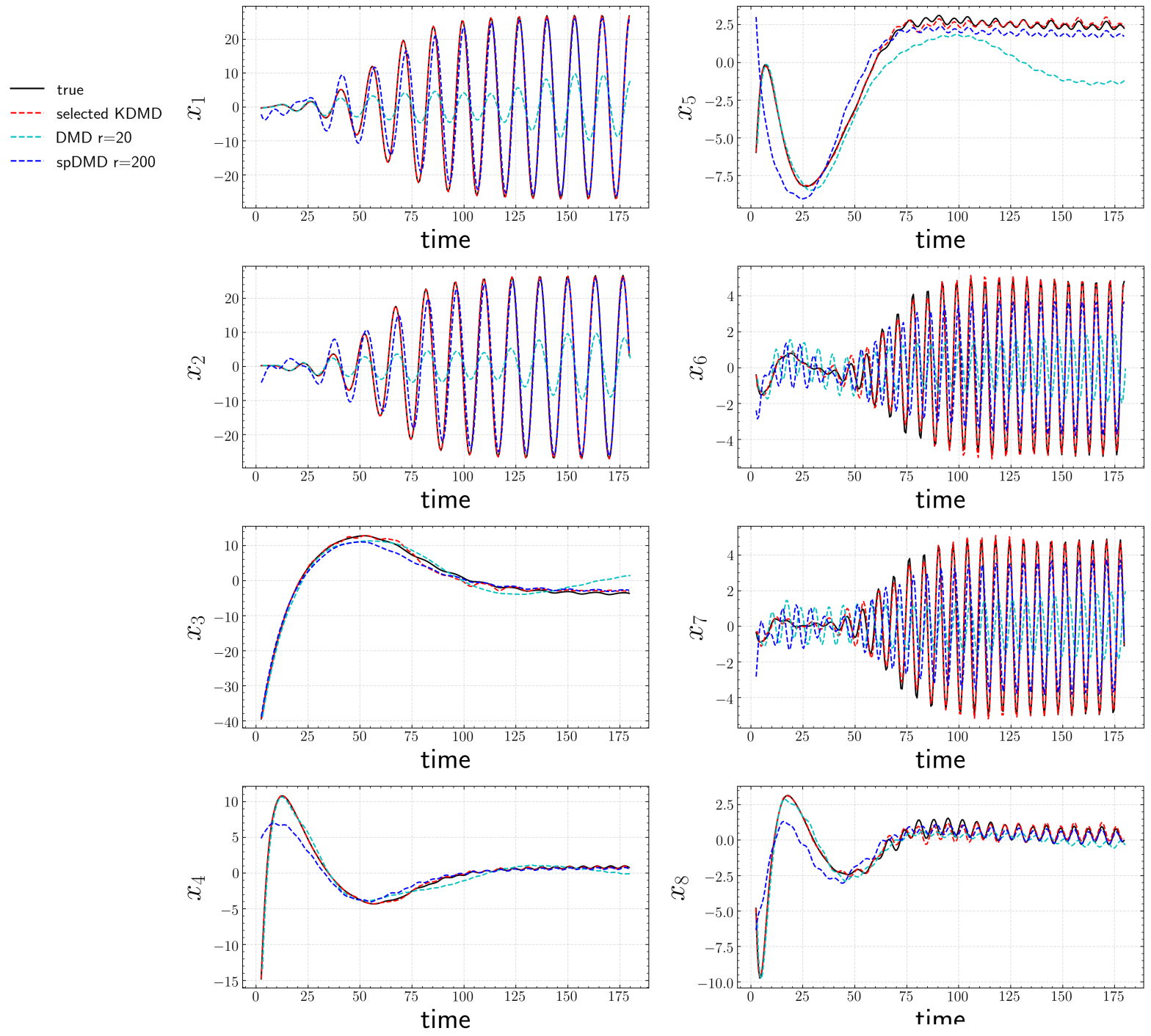}
\caption{Comparison of a posteriori prediction on the top 8 POD coefficients of the testing trajectory between sparsity-promoting KDMD, DMD~\citep{schmid2010dynamic} and spDMD~\citep{jovanovic2014sparsity} for the 2D cylinder flow at $\Rey=70$. $x_i$ denotes $i$-th POD coefficient.} 
\label{fig:20d_cyd_compare_spdmd}
\end{figure}

As shown in \cref{fig:20d_cyd_compare_spdmd}, given the same number of eigenmodes, sparsity-promoting KDMD performs remarkably better than spDMD, especially in the transient regime. \textcolor{cyan}{This is likely due to the inability of DMD in capturing non-linear algebraic growth in the transient regime~\citep{bagheri2013koopman}.} spDMD overestimates  the growth in $x_1$ and $x_2$ while ignoring a turnaround near the onset of transient regime in $x_5$ and $x_8$. As expected, DMD with 20 POD coefficients performs the worst especially for the modes where transient effects are dominant. Given the results in \cref{fig:20d_cyd_compare_spdmd}, among all of the top 8 POD coefficients, $x_6$ and $x_7$ appear to be most challenging to model: DMD and spDMD cannot match the limit-cycle while spKDMD performs very well. Notice that the frequency in $x_6$ and $x_7$ corresponds to $St_D$. Hence, there will be a difference in predicting the fluctuation of $C_D$ between spDMD and sparsity-promoting KDMD. 

\begin{figure}
\centering
\includegraphics[width=0.6\textwidth]{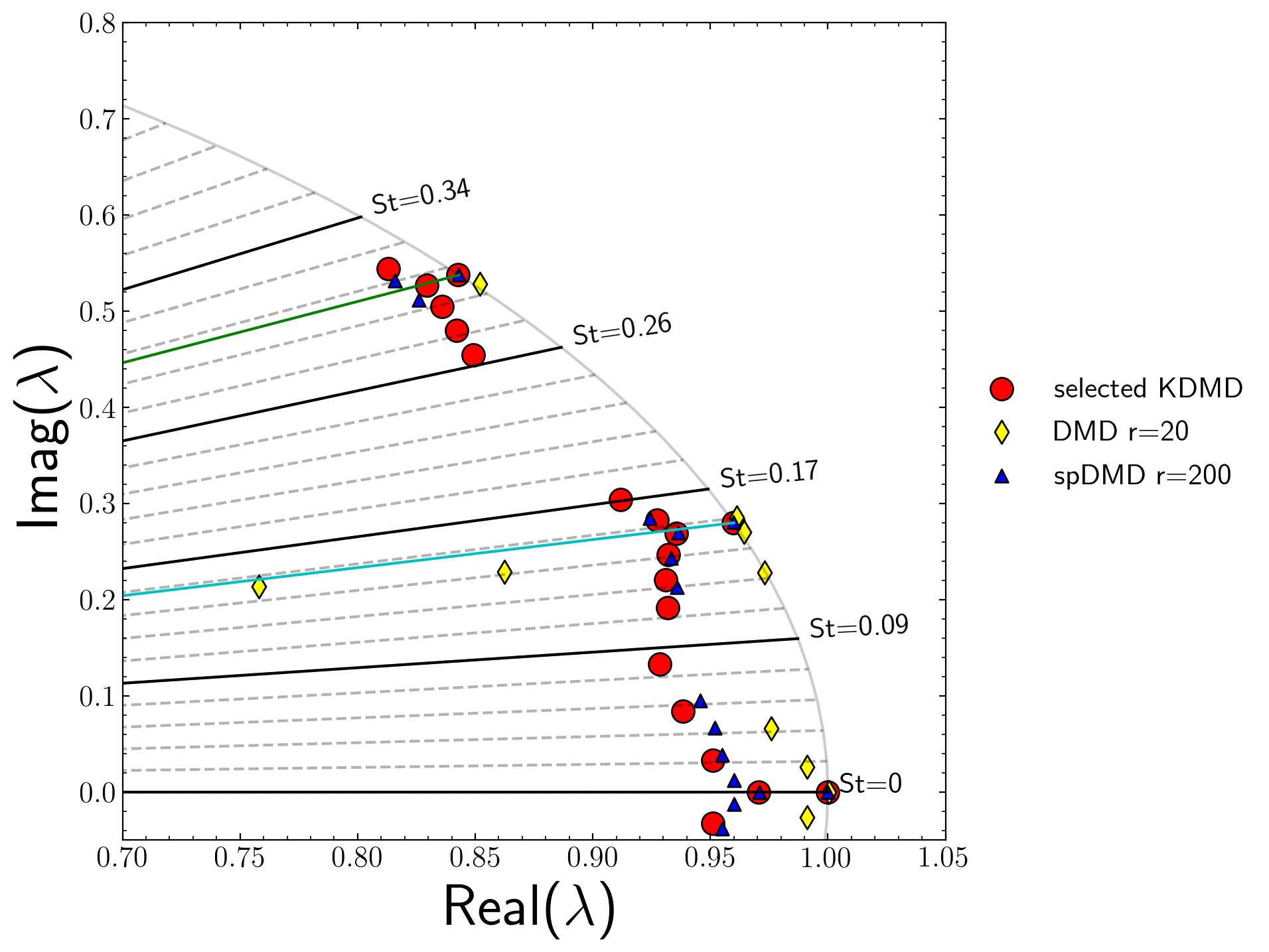}
\caption{Comparison of identified eigenvalues between spKDMD, DMD~\citep{schmid2010dynamic} and spDMD~\citep{jovanovic2014sparsity} for the 2D cylinder flow at $\Rey=70$}
\label{fig:20d_cyd_compare_spdmd_eigenvalue}
\end{figure}

Finally, comparison of the identified Koopman eigenvalues between DMD, spDMD and spKDMD is shown in \cref{fig:20d_cyd_compare_spdmd_eigenvalue}. On one hand, both spDMD and spKDMD exactly capture the stable eigenmodes that correspond to  $St_D$ and $St_L$. This is expected since DMD with 200 POD coefficients represents the dynamics very well, and deviation  from limit-cycle behavior would be penalized in spDMD. On the other hand, several erroneous stable DMD modes are obtained by DMD. This explains the deviation of a posteriori prediction  from the ground truth limit cycle in \cref{fig:20d_cyd_compare_spdmd}. For those decaying modes,  similarity/co-linearity is observed between two clusters of eigenvalue from spDMD and spKDMD. However, spKDMD contains more high frequency modes than spDMD. Finally, it is interesting to note that, although the correct stable eigenvalues are captured accurately by both spDMD and spKDMD, the former does not capture  accurate amplitudes for stable eigenvalues of type-II as seen in \cref{fig:20d_cyd_compare_spdmd}. 

As a side note, when the temporal mean was used instead of maximal error in the definition of $Q$ in \cref{eq:linear_evolve_2}, spKDMD with the above setting was found to not find a stable eigenvalue corresponding to  $ST_D$. 

\subsection{Three-dimensional transient turbulent ship airwake}

\begin{figure}
\centering
\includegraphics[width=\textwidth]{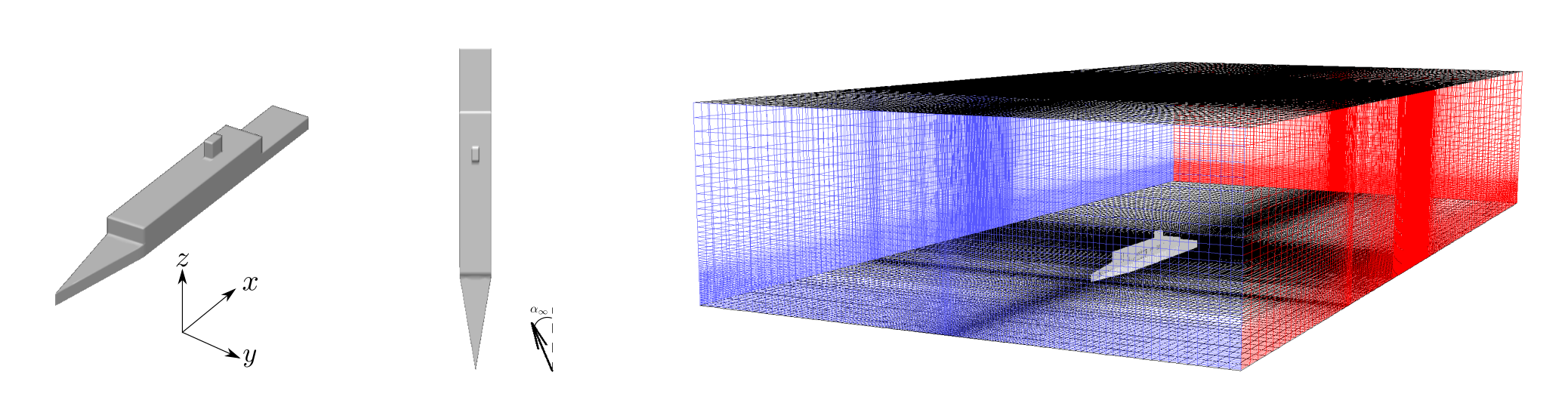}
\caption{Left: geometry of the ship (SFS2). Right: generated computational mesh.}
\label{fig:ship_geom}
\end{figure}

Understanding unsteady ship airwake flows is critical to design safe shipboard operations, such as takeoff and landing of fixed or rotary wing aircraft  \citep{forrest2010investigation}, especially when the wind direction becomes stochastic. Here we obtain snapshots from an unsteady Reynolds Averaged Navier--Stokes (URANS) simulation of a ship airwake using FLUENT \citep{ansys} with the shear layer corrected $k$-$\omega$ two-equation turbulence model. Unsteadiness arises from both bluff-body separation, and an abrupt change in wind direction. We consider a conceptual ship geometry called simple frigate shape 2 (SFS2). For the details of the geometry, readers are referred to \citep{yuan2018combined}. Configuration of the simulation setup is shown in \cref{fig:ship_geom} where $\alpha_{\infty}$ denotes the angle of side wind. 

To prepare a proper initial condition, a URANS simulation for $U_{\infty}=15 \mathrm{m/s}$ with $\alpha_{\infty}=0^{\circ}$, i.e., no side wind, is conducted to reach a physical initial condition. Following this, the last snapshot is used as the initial condition for a new run with an impulsive change in the wind direction from $\alpha_{\infty} = 0^{\circ}$ to $\alpha_{\infty}= \alpha_0=5^{\circ}$. 
The boundary conditions for outlet/input is pressure outlet/velocity inlet while top and bottom are set as symmetry for simplicity. No-slip conditions are used at surface of the ship. Further details on the simulation setup are provided in \citep{sharma2019simulation}.

The sampling time interval is $\Delta t = 0.1\mathrm{s}$ with 500 consecutive samples of the three velocity components. This corresponds to several flow through times over the ship length. The domain of interest is a cartesian region of mesh size $24\times 40\times 176$ starting on the rear landing deck. For dimension reduction, the trajectory of the top 40 POD coefficients (temporal mean subtracted) are collected, yielding $> 99\%$ kinetic energy preservation.
Discrete-time KDMD \textcolor{blue}{with an isotropic Gaussian kernel} is employed to perform nonlinear Koopman analysis where $\sigma = 200$, $r=135$ is chosen. 
Details of hyperparameter selection  are provided in \cref{apdx:hyper_ship}.

\subsubsection{Results of discrete-time KDMD with mode selection}

\begin{figure}
\centering
\includegraphics[width=\textwidth]{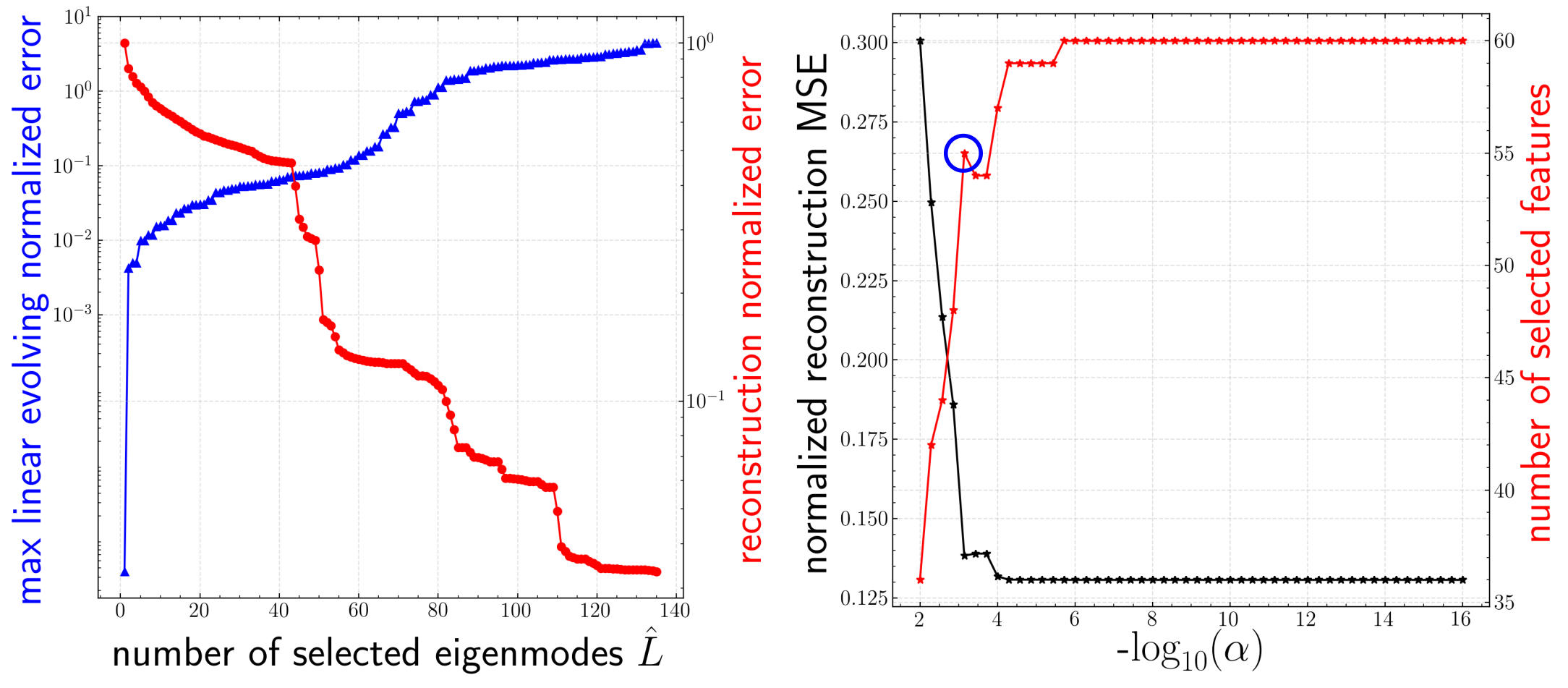}
\caption{Left: Trend of linearly evolving error $Q$ and reconstruction error $R$ from discrete-time KDMD for the ship airwake. Right: Trend of linearly evolving error $Q$ and reconstruction error $R$ from discrete-time KDMD.}
\label{fig:ship_qr_mse_5}
\end{figure}

First, the error analysis of eigenfunctions is shown in \cref{fig:ship_qr_mse_5}, where we choose $\hat{L} \approx 60$ for good reconstruction. 
However, the level of deviation from linear evolution is around 10\%. This error will be reflected later as deviation in a posteriori prediction on the testing trajectory\footnote{This implies difficulties in finding an accurate yet informative Koopman operator with isotropic Gaussian kernels. However, choosing an optimal kernel type is beyond the scope of this work. }.

\begin{figure}
\centering
\includegraphics[width=\textwidth]{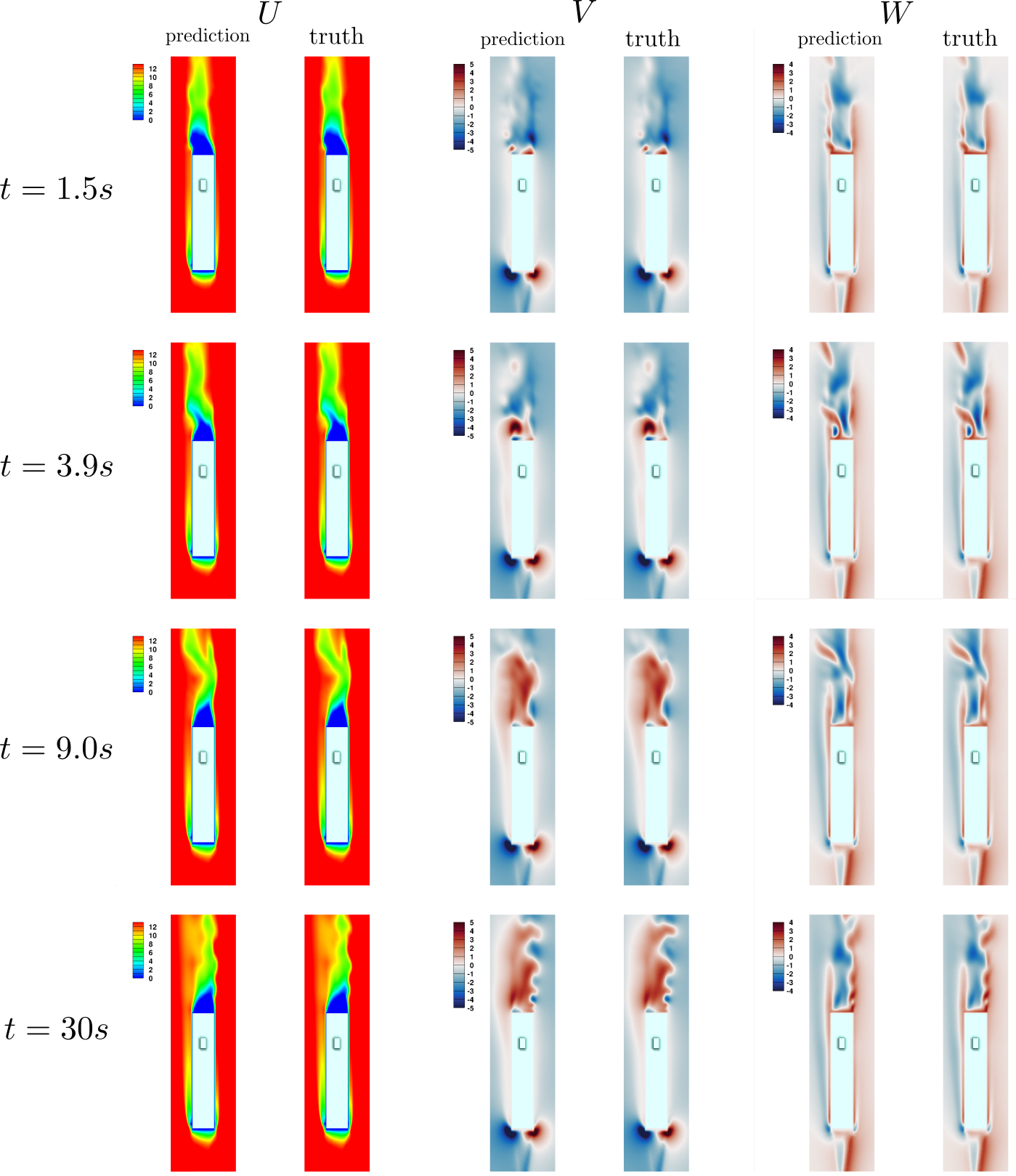}
\caption{Contour of velocity components near the ship on $z-$plane slice at $t=1.5s$, 3.9s, 9.0s, 30s. For each subfigure, left: prediction from KDMD; right: ground truth.}
\label{fig:ship_5_pred_vs_true}
\end{figure}

Second, the result of mode selection is summarized in \cref{tab:modes_selection_ship}. Note that nearly 2/3 modes are removed. Furthermore, model performance in terms of a posteriori prediction on the testing trajectory is evaluated. Comparison between KDMD and the ground truth on contours of velocity components on a special $z-$plane (1.2 meters above the landing deck) is shown in \cref{fig:ship_5_pred_vs_true}. Effects of an impulse change in wind direction in the following are observed from $t=1.5s$ to $t=30s$ and well-captured by spKDMD:
\begin{itemize}
    \item growth of a large shear layer in $U$ on the rear (left) side of the superstructure on the ship
    \item a strong side wind sweep in $V$ above the landing deck propagating  downstream
    \item development of vortex on the upwind (right) downstream side of the ship
\end{itemize}

\begin{table}
 \begin{center}
\def~{\hphantom{0}}
 \begin{tabular}{lc}
  $\alpha_{\textrm{select}}$ & $7.19\times 10^{-4}$  \\
  number of selected modes & 55   \\
  number of total modes & 167 \\
  normalized reconstruction error $R$ & 0.133  \\ 
 \end{tabular}
\caption{Summary of mode selection for discrete-time KDMD on ship airwake.}
\label{tab:modes_selection_ship}
 \end{center}
\end{table}

\begin{figure}
\centering
\includegraphics[width=1\textwidth]{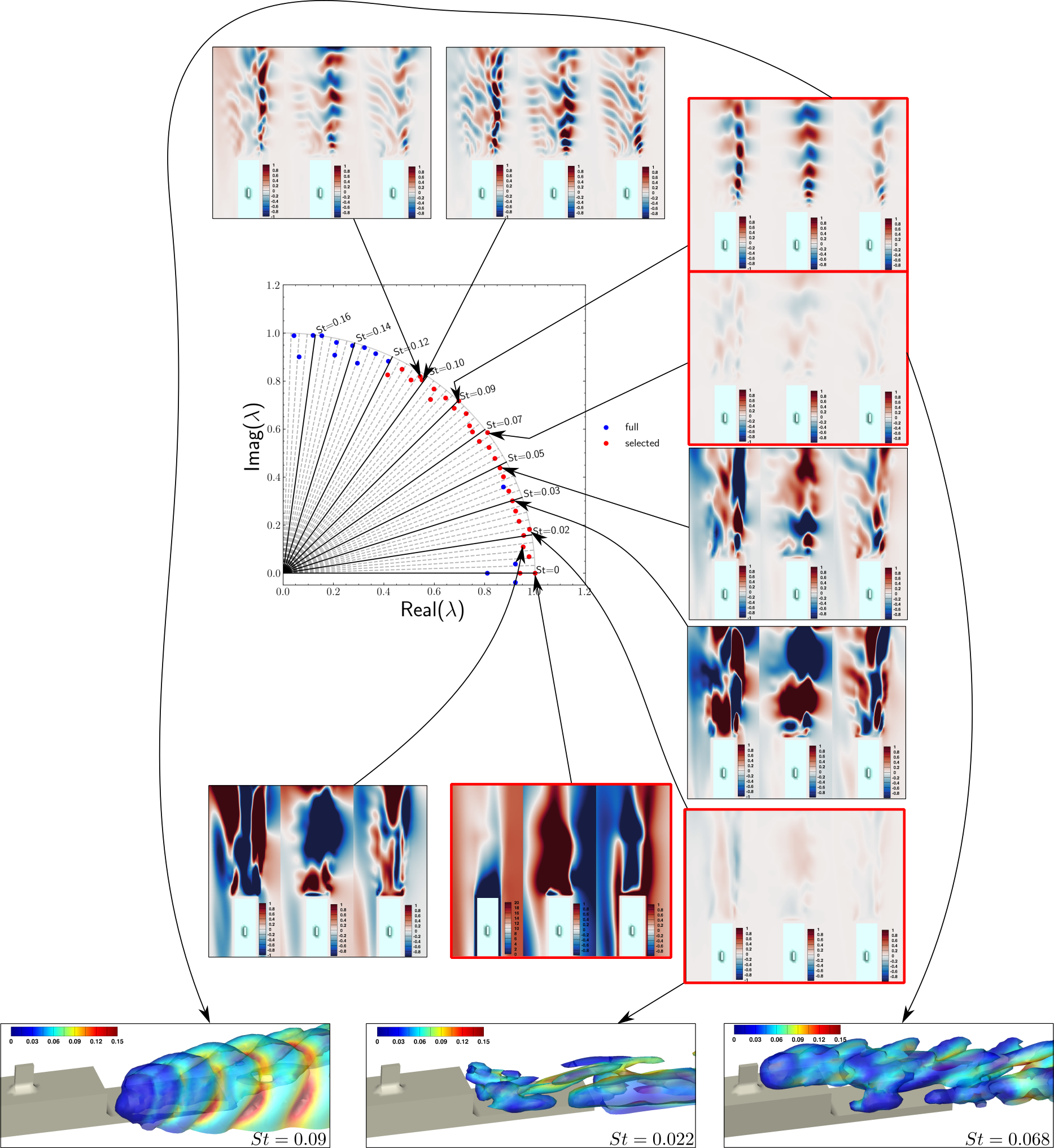}
\caption{Contours of Koopman modes of ship airwake on the $z$-plane at $t=0$. For each subfigure, left: $U$, middle: $V$, right: $W$. Red squares indicate stable modes. Bottom: iso-contour of vorticity colored by velocity magnitude \textcolor{green}{zoomed near the landing deck}. } 
\label{fig:40d_ship_alpha_5_modes}
\end{figure}

Further, three velocity components of the Koopman mode decomposition on the previously mentioned $z-$plane is shown in \cref{fig:40d_ship_alpha_5_modes} together with the isocontour of vorticity colored by the velocity magnitude for the two stable harmonics.
Note that frequency is normalized using $U_{\infty}$ as the reference velocity and funnel width of the ship $L=3.048 \mathrm{m}$ as the characteristic length scale \citep{forrest2010investigation}. 
As expected, the spKDMD  yields a large number of decaying modes with only three nontrivial stable harmonics, since significant transient effects are present in the data. From the Koopman mode decomposition in \cref{fig:40d_ship_alpha_5_modes}, we observe the following:
\begin{itemize}
    \item  modes with eigenvalues close to each other exhibit similar spatial structure,
    \item modes associated with higher frequency are dominated by smaller scales,
    \item the stable harmonic mode near $St = 0.09$ associated with a strong cone-shape vortex originating from the upwind (right) rear edge of the superstructure on the ship,
    \item the stable harmonic mode near $St = 0.068$ corresponds to vortex shedding induced by the funnel,
    \item \textcolor{green}{the slowly-decaying mode near $St=0.022$ shows unsteady circulation directly behind the superstructure,} 
    \item the steady  mode ($St=0$) is consistent with the large circulation behind the superstructure, the roll-up wind on the side of landing deck, and vertical suction towards the floor on the landing deck.
\end{itemize}

\subsubsection{Comparison with sparsity-promoting DMD}

\begin{figure}
\centering
\includegraphics[width=1\textwidth]{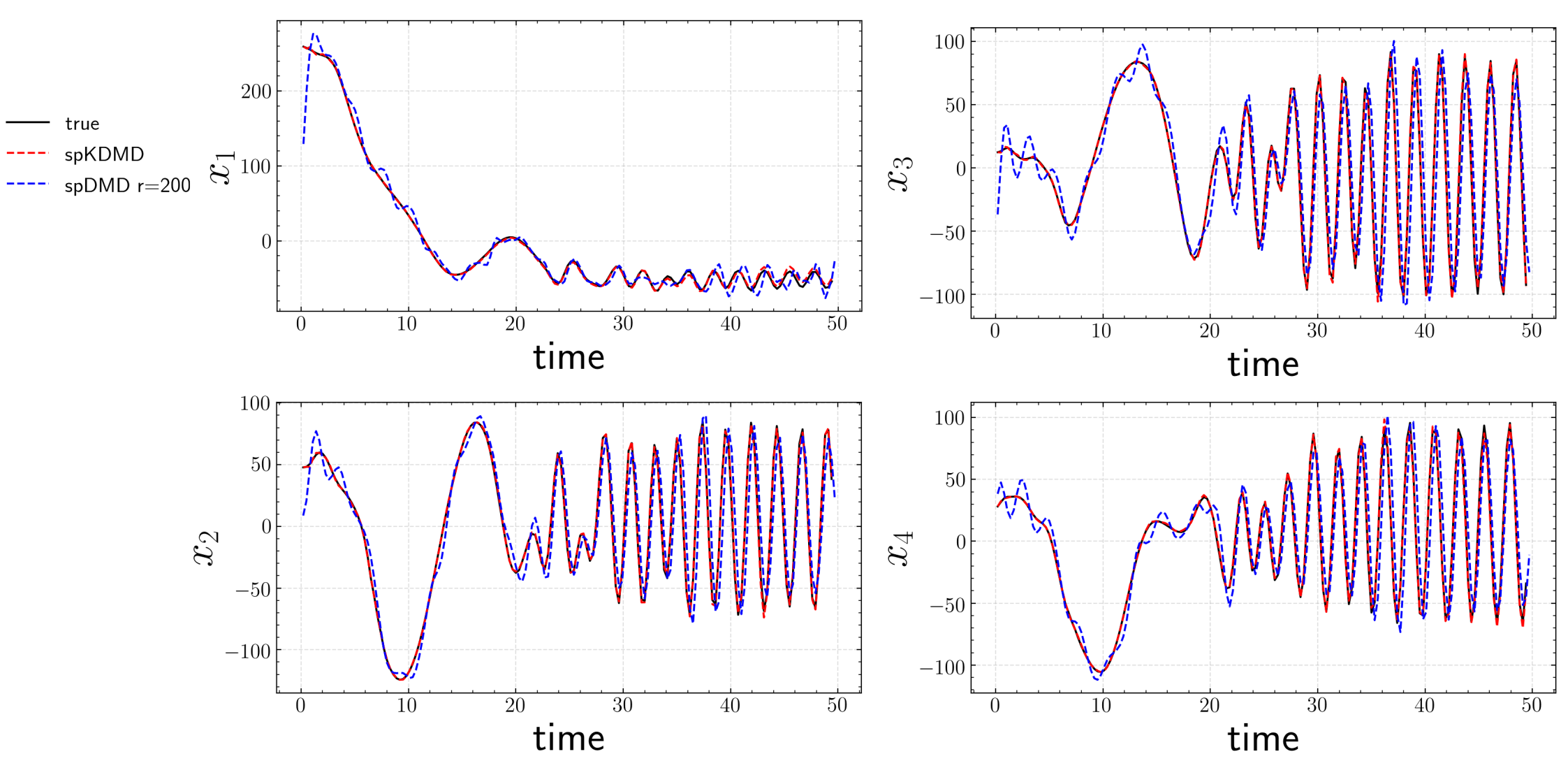}
\caption{Comparison of a posteriori prediction of the  4 most significant POD coefficients of the testing trajectory between sparsity-promoting KDMD and spDMD~\citep{jovanovic2014sparsity} for the 3D ship-airwake flow. $x_i$ denotes $i$-th POD coefficient.} 
\label{fig:40d_ship_compare_spdmd}
\end{figure}

\textcolor{cyan}{We again repeat the comparison between our spKDMD and spDMD~\citep{jovanovic2014sparsity}. Note that DMD on the first 40 POD modes performs poorly similar to \cref{sec:compare} and therefore not shown here. To make a fair comparison against the spKDMD from previous subsection, however, we collect the first 200 POD coefficients for spDMD to ensure that a sufficient number of modes are used to fit the trajectory well. We then carefully choose the penalty coefficient in spDMD to ensure that the same number of modes are retained  as in spKDMD. As shown in \cref{fig:40d_ship_compare_spdmd}, within the time horizon $t<50$, a posteriori evaluation shows that spKDMD offers much improved predictions compared  to spDMD~\citep{jovanovic2014sparsity} on the testing trajectory.}

\begin{figure}
\centering
\includegraphics[width=\textwidth]{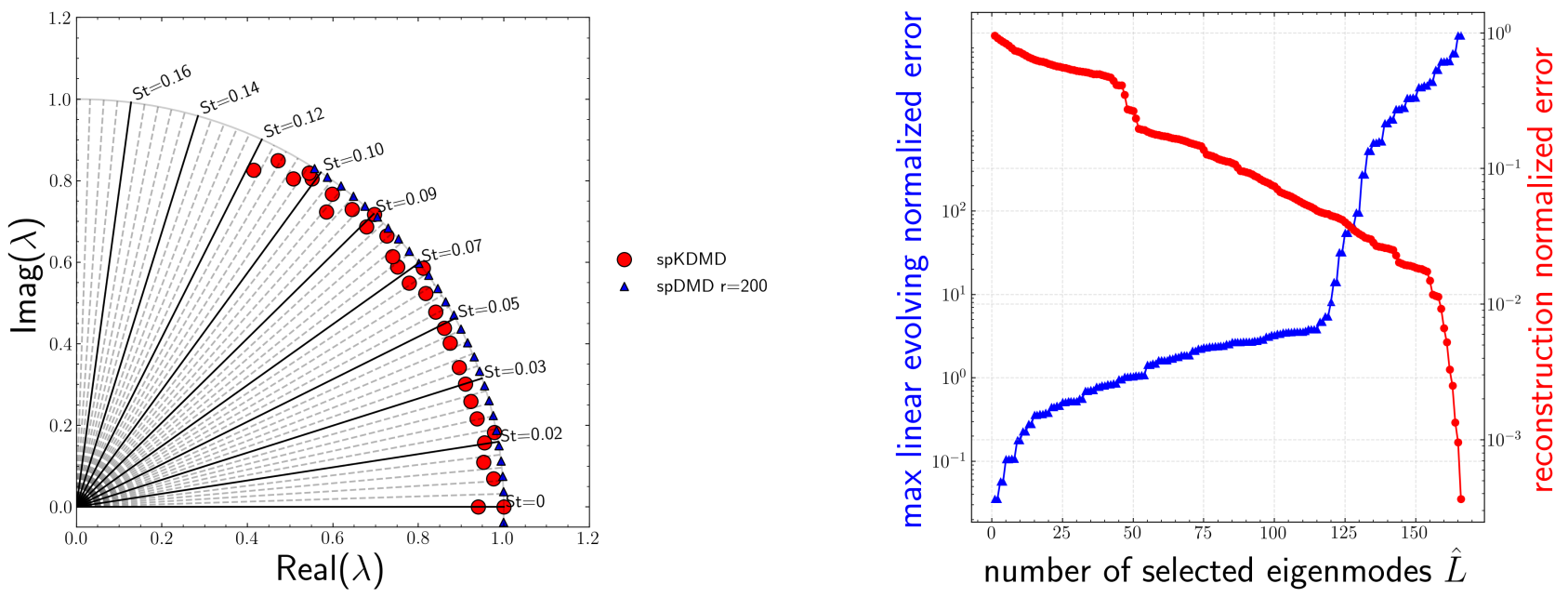}
\caption{Left: Comparison of identified eigenvalues between spKDMD and spDMD~\citep{jovanovic2014sparsity} for the 3D ship-airwake flow. Right: Trend of linear evolving error $Q$ and reconstruction error $R$ from DMD for the 3D ship-airwake flow. }
\label{fig:40d_ship_compare_spdmd_eigenvalue}
\end{figure}

\textcolor{cyan}{
Moreover, as further illustrated in the left subfigure of \cref{fig:40d_ship_compare_spdmd_eigenvalue}, eigenvalues identified from spKDMD only contain two stable modes while nearly all eigenvalues from spDMD are located near the unit circle, among which there are 30 out of 56 slightly unstable modes. These unstable modes inevitably lead to \emph{the identified system being numerically unstable} when predicting beyond the current training time horizon, whereas spKDMD predicts a ``physically consistent'' limit cycle behavior. As indicated in the right subfigure of \cref{fig:40d_ship_compare_spdmd_eigenvalue}, such instability is related to the inability of the (linear) features to approximate the Koopman-invariant subspace, where only 8 modes are within 10\% of maximal deviation from linear evolution, compared to 60 modes in KDMD. We note that similar observations of the drastically different eigenvalue distribution were reported in the original KDMD paper~\citep{williams2014kernel}.}


\section{Conclusions}
\label{sec:conclusion}

Two classical nonlinear approaches for the approximation of Koopman operator: EDMD and KDMD are revisited. From an algorithmic perspective, the main contributions of this work are a) Sparsity promoting techniques based on a posteriori error analysis, and b) multi-task learning techniques for mode selection as an extension of spDMD into the nonlinear variants. Further, analytical relationships between spDMD, Kou's criterion and the proposed method are derived from the viewpoint of optimization. The algorithm~\footnote{code available at: https://github.com/pswpswpsw/SKDMD} is first evaluated in detail on a simple two state dynamical system, for which the Koopman decomposition is known analytically.

If one is only interested in the post-transient dynamics of the system   \emph{on} an attractor, linear observables with time delays are sufficient to extract an informative Koopman-invariant subspace. Thus, the present techniques are evaluated on two fluid flows which involve strong transients: 2D flow over a cylinder at different Reynolds numbers, and a 3D ship airwake. We demonstrate the effectiveness of discovering \emph{accurate} and \emph{informative} Koopman-invariant subspaces from data and constructing  accurate reduced order models from the viewpoint of Koopman theory. Furthermore, with the proposed algorithm, the parametric dependency of Koopman mode shapes on the Reynolds number is investigated for the cylinder flows. In this case,  as $\Rey$ increases from 70 to 130, the shape of stable modes become flatter downstream and larger upstream. Moreover, the similarity of mode shapes between Koopman modes with similar eigenvalues is observed in both fluid flows. Specifically, five clusters of eigenvalues  are observed in the case of 2D cylinder wake flow \textcolor{cyan}{which is confirmed with weakly nonlinear theoretical analysis from \cite{bagheri2013koopman}.}  Type-I, II clusters are found to correspond to  fluctuations in lift and drag respectively. {\color{brown}We identify non-oscillatory as well as  oscillatory cancellations from the above two clusters with distinct frequencies. These modes could be understood as ``oscillatory shift modes'', as a generalization of the model proposed by~\cite{noack2003hierarchy}.}
For the 3D ship airwake case, two stable modes, \textcolor{green}{and one slowly-decaying}  mode with distinct frequencies and mode shapes resulting from vortex shedding are extracted, and  accurate predictive performance is observed in the transient regime. 

    
    
    

\section*{Acknowledgements}
\textcolor{brown}{We thank the three anonymous reviewers whose comments and suggestions helped improve and clarify this manuscript.} This research was supported by the DARPA Physics of AI Program under
the grant ``Physics Inspired Learning and Learning the Order and Structure of Physics,'' and by the NSF CMMI program under the grant    ``A Diagnostic Modeling Methodology for Dual Retrospective Cost Adaptive Control of Complex Systems''. Computing resources were provided by the NSF via the grant ``MRI: Acquisition of ConFlux, A Novel
Platform for Data-Driven Computational Physics.''


\section*{Declaration of interests}
The authors report no conflict of interest.

\appendix

\section{Proof for the identity for \cref{apdx:equation}}\label{apdx:idtt}

Consider the economy SVD of $\mathbf{A}^{'} = \mathbf{U} \mathbf{S} \mathbf{V}^{\textrm{H}}$. 
Since the column space of $\mathbf{A}^{'}$ spans the column space of $\mathbf{Z}_r$, there exists an invertible matrix $\mathbf{P}$ such that $\mathbf{Z}_r \mathbf{P} = \mathbf{U}$. Hence, $$\mathbf{A}^{'} \mathbf{Z}_r = \mathbf{V} \mathbf{S} \mathbf{U}^{\textrm{H}} \mathbf{Z}_r = \mathbf{V} \mathbf{S} \mathbf{P}^{\textrm{H}} \mathbf{Z}_r^{\textrm{H}} \mathbf{Z}_r = \mathbf{V}\mathbf{S} \mathbf{P}^{\textrm{H}},$$ and $((\mathbf{A}^{'} \mathbf{Z}_r)^{\textrm{H}} \mathbf{A}^{'} \mathbf{Z}_r)^{+} = (\mathbf{P}^{-1})^{\textrm{H}}  \mathbf{S}^{-2} \mathbf{P}^{-1}$.
Thus $$
 \mathbf{A}^{'} \mathbf{A}^{'\textrm{H}} \mathbf{Z}_r (\mathbf{Z}^{\textrm{H}}_r  \mathbf{A}^{'} \mathbf{A}^{'\textrm{H}} \mathbf{Z}_r)^{+} = (\mathbf{U} \mathbf{S} \mathbf{V}^{\textrm{H}})  (\mathbf{V} \mathbf{S} \mathbf{P}^{\textrm{H}} ) (\mathbf{P}^{-1})^{\textrm{H}} \mathbf{S}^{-2} \mathbf{P}^{-1} 
  = \mathbf{U} \mathbf{P}^{-1} = \mathbf{Z}_r.$$

\section{Computational complexity}
\label{apdx:cc}

\textcolor{red}{
Computational cost of the proposed framework can be divided in three parts:
\begin{enumerate}
    \item Distributed SVD
    \item  KDMD/EDMD algorithms
    \item Multi-task feature learning  with the parameter sweep (solving \cref{eq:multi_task_enet} with different values of $\alpha$)
\end{enumerate}
\quad
The computational complexity for each step in the algorithm is summarized in \cref{tab:1}, where $n$ is the dimension of the system state, $M$ is the number of snapshots in the training~\footnote{for conciseness, we assume the number of training snapshots equals the number of validation snapshots.}, $r$ is the rank of the reduced system after SVD, $\hat{L}$ is the user-defined cut off for ``accurate'' features, and $N_{iter}$ is the maximal number of iterations user-defined to achieve a residual threshold, e.g. $10^{-12}$.}

\textcolor{red}{
As shown in the Table.~\ref{tab:1}, the theoretical computational complexity for multi-task ElasticNet with an $\alpha$ sweep is $O(N_{\alpha}N_{iter}\hat{L}^2r)$. Note that this is a worst case simply because - except for the first $\alpha$ - we \emph{reuse the previous optimal solution} as the initial condition for the new objective function to start the iterative optimization process. Also, thanks to SVD-based dimension reduction, the cost  scales linearly with the reduced system rank $r$. Moreover, the user-defined linearly evolving error truncation $\hat{L}$ helps reduce that complexity as well instead of scaling with the number of snapshots $M$. Lastly, there is a cubic theoretical complexity associated with the number of snapshots when applying KDMD. The number of snapshots in a typical high fidelity simulation is $O(10^3)$. That is to say,  $r<10^3$ and $\hat{L} < 10^3$. We note that computational efficiency can be further improved,  but this will be left for future work.
\begin{table}
 \begin{center}
\def~{\hphantom{0}}
 \begin{tabular}{lc}
 Step & Computational complexity \\[3pt]
 SVD (QR iteration with vectors)& $O(Mn^2)$~\citep{dongarra2018singular}  \\ 
 KDMD  & $O(M^3)$   \\
  error evaluation \& pruning & $O(M^2n)$ \\
 multi-task ElasticNet with a fixed penalty coef. $\alpha$ & $O(N_{iter}\hat{L}^2r)$  \\  
 \textrm{multi-task ElasticNet with $N_{\alpha}$ coefs.} (worst case) & $O(N_{\alpha}N_{iter}\hat{L}^2r)$  \\
 \end{tabular}
 \caption{Computational complexity of each step in the proposed sparsity-promoting framework.}
 \label{tab:1}
 \end{center}
\end{table}
}

\section{Hyperparameter selection}
\subsection{2D fixed point attractor}
\label{apdx:hyper_2d}

\begin{figure}
\centering
\includegraphics[width=0.6\textwidth]{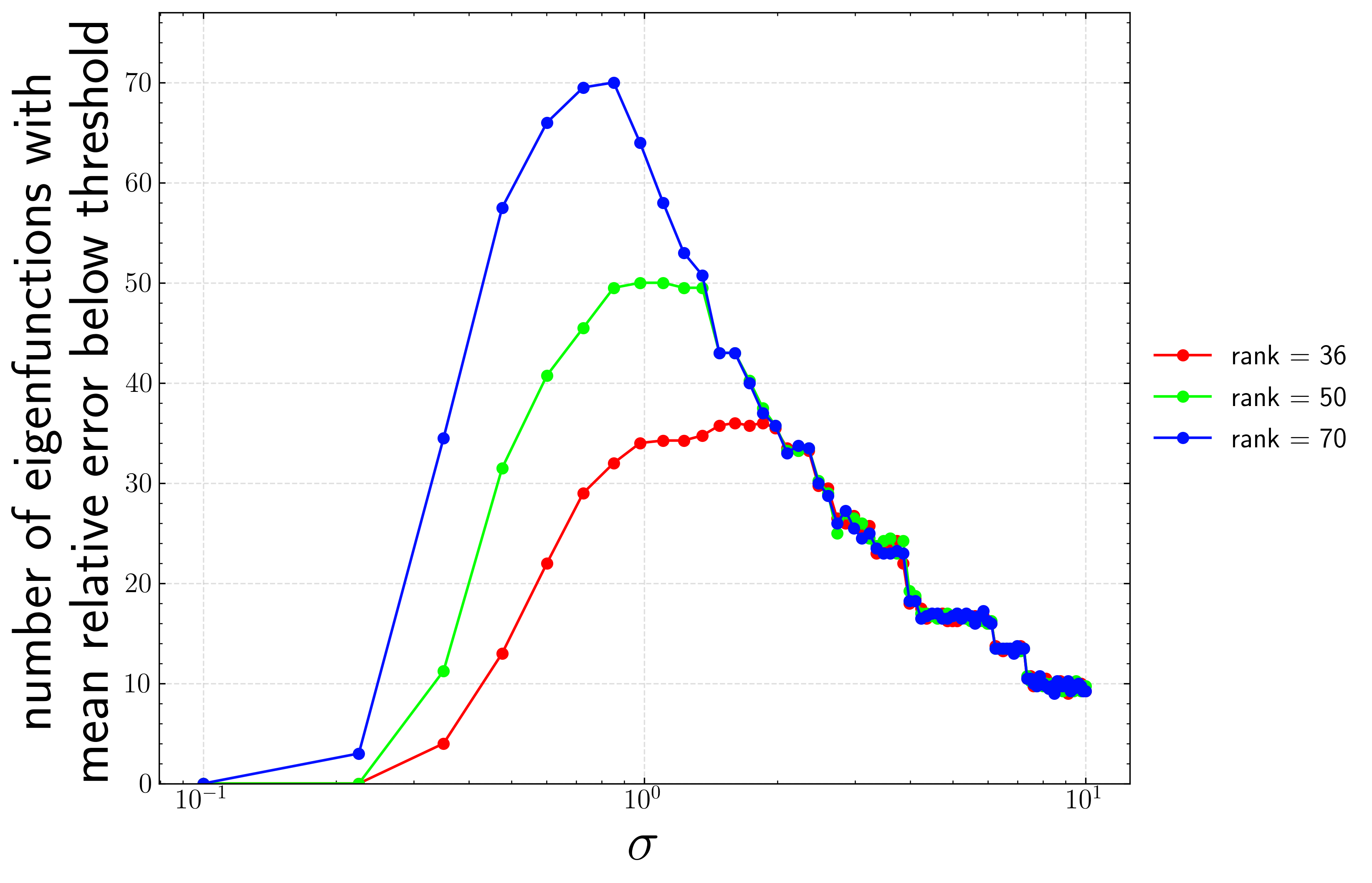}
\caption{Hyperparameter search for isotropic Gaussian KDMD on the 2D fixed point attractor.}
\label{fig:hyp_2d_fp_data}
\end{figure}

We perform grid search in parallel for the selection of $\sigma$ and the truncated rank $r$ over the range: $\sigma \in [10^{-1}, 10]$ with 80 points uniformly distributed in the log sense and $r=36,50,70$ to find a proper combination of $r$ and $\sigma$. As shown in \cref{fig:hyp_2d_fp_data}, the higher rank $r$ leads to larger number of linearly evolving eigenfunctions. Thus, it is more crucial to choose a proper scale $\sigma$ than $r$ from \cref{fig:hyp_2d_fp_data}. However, considering the simplicity of this problem, $\sigma=2$ and $r=36$ would suffice. 

\subsection{Cylinder flow case}
\label{apdx:hyper_cyd}

\begin{figure}
\centering
\includegraphics[width=0.6\textwidth]{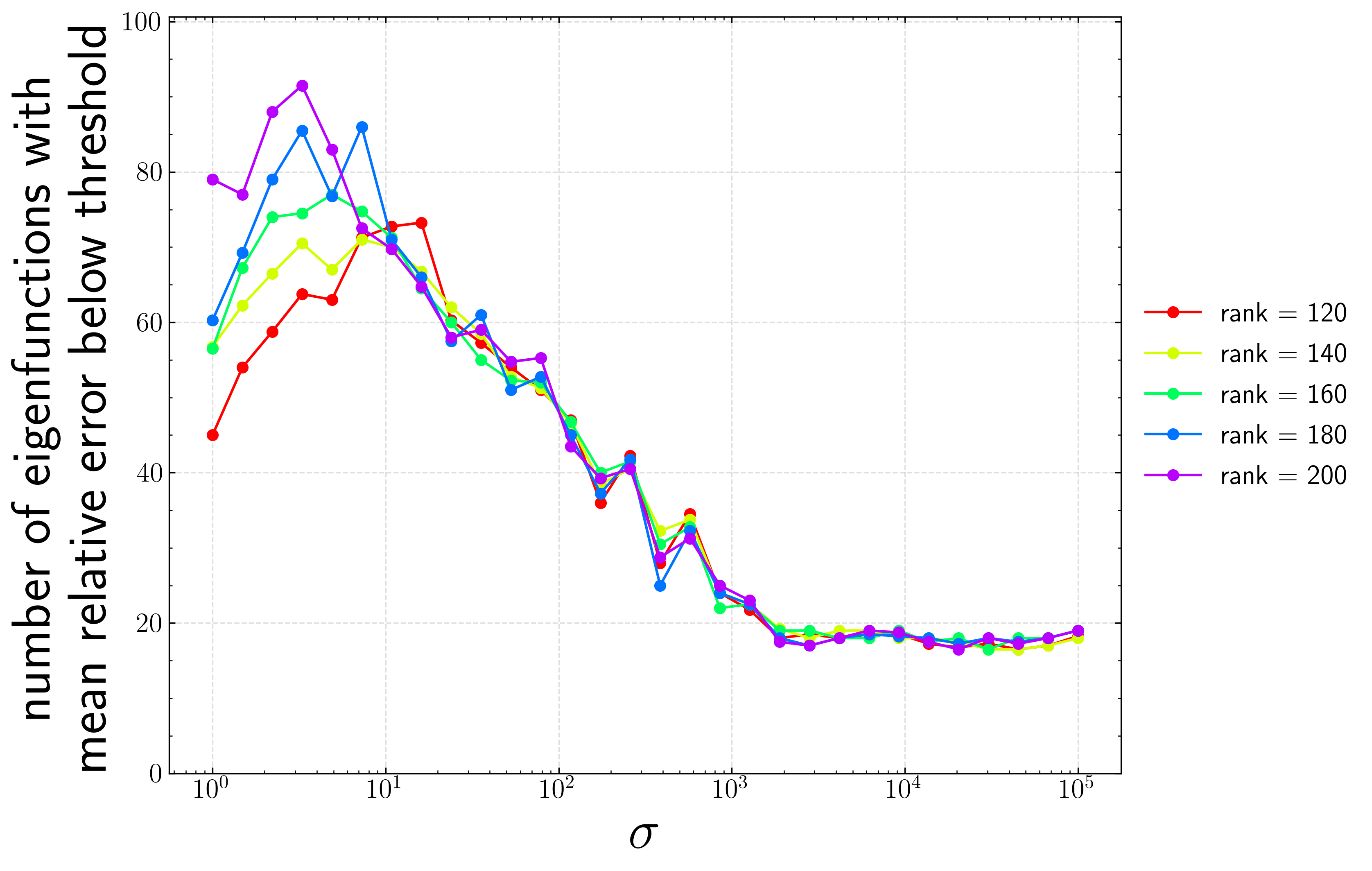}
\caption{Hyperparameter search for isotropic Gaussian KDMD on transient cylinder wake flows.}
\label{fig:hyp_20d_cyd_data}
\end{figure}

As the flow is not turbulent, we choose  hyperparameters for the $Re=70$ case and fix them for all the other cases. We sweep $\sigma$ from $[1,10^5]$ with 30 points uniformly distributed in the log sense and $r = 120,140,160,180, 200$ as shown in figure~\ref{fig:hyp_20d_cyd_data}. From the plot, we choose $r=180$ and $\sigma=3$ for an approximate  choice of the hyperparameter.
Again, we observe that the number of accurate eigenfunctions first increases then decreases with $\sigma$ increasing and the saturation of rank truncation at higher $\sigma$ which is related to the variation in the characteristic scale of the features with respect to $\sigma$.

\subsection{Turbulent Ship-airwake case}
\label{apdx:hyper_ship}

\begin{figure}
\centering
\includegraphics[width=0.6\textwidth]{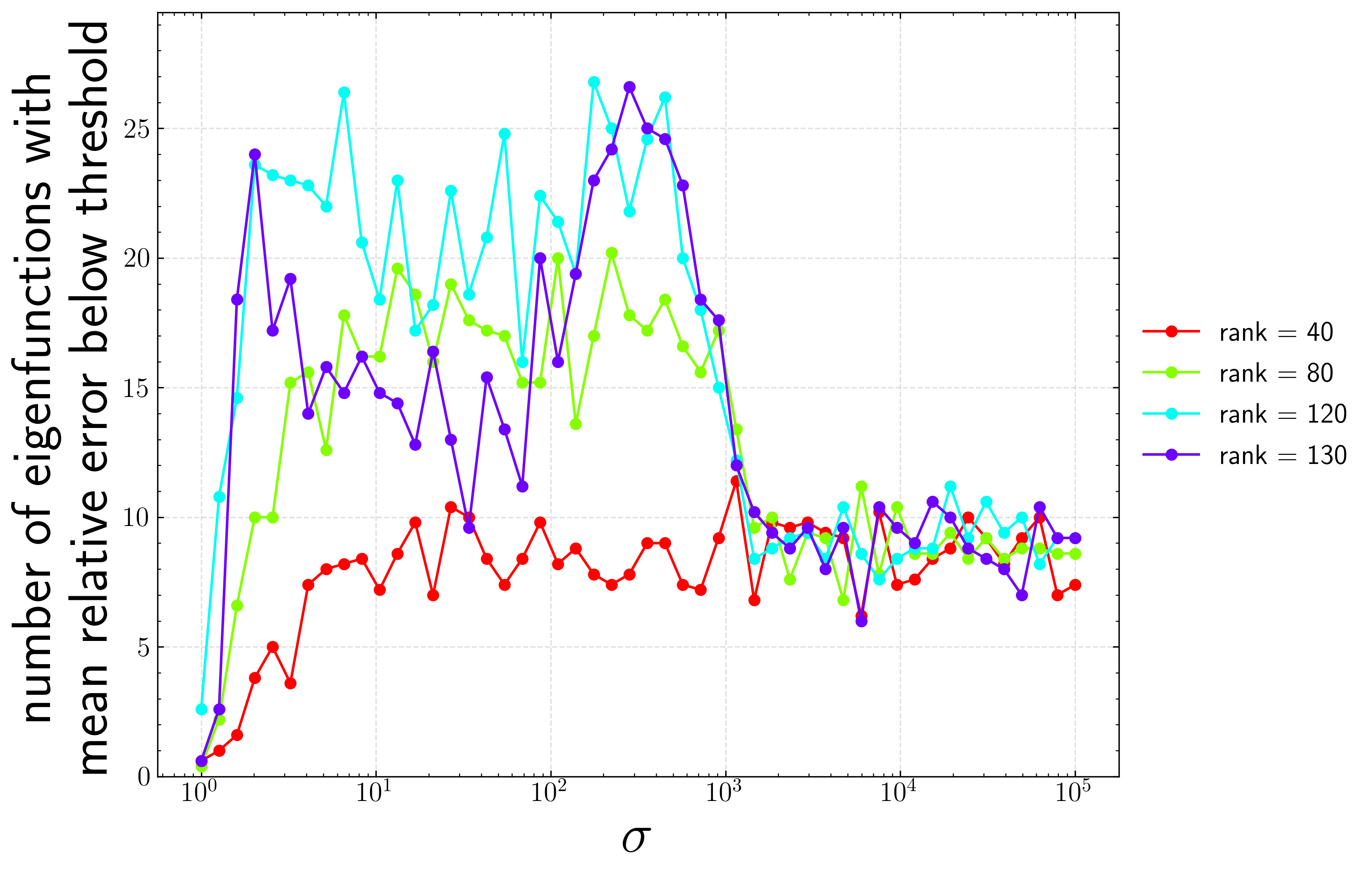}
\caption{Hyperparameter search for isotropic Gaussian KDMD on transient ship airwake.}
\label{fig:hyp_30d_ship_data}
\end{figure}

Grid search in parallel for the selection of $\sigma$ and $r$ is performed over the range $\sigma \in [1,10^5]$ with 50 points uniformly distributed in the log sense, $r=40,80,120, 130$. As shown in \cref{fig:hyp_30d_ship_data}, a good choice of $\sigma$ can be 200 for the case of $\alpha_{\infty}=5^{\circ}$. Note that since the hyperparameter selection is performed with a 5-fold cross validation on the training data, we only have \textcolor{brown}{up to} $166*0.8 \approx 132$ data points, i.e., maximal possible rank is 132. While in the actual training, we have maximal rank up to 166. Note that as long as the system is  well-conditioned, the higher the rank, the richer the subspace. Here we take $\sigma = 200$ and $r=135$.

\section{Parallel Data Processing (PDP) Tool}
\label{apdx:pdp}
The application of data-driven methodologies to increasingly large data sets has exposed the need for scalable tools.  Method development is easily achieved in scripting environments, but for application to large multidimensional problems, memory and data I/O becomes a limiting factor. 
To illustrate this, consider a state-of-the art 3D computation. Given  the fact that there are multiple variables associated with each grid cell, and that reasonable meshes can be expected to contain millions to hundreds of millions of cells, work-space memory requirements of Terabytes may be required.
Utilizing available distributed tools and I/O strategies a generalized toolset, Parallel Data Processing(PDP), has been developed for application of complex methods to arbitrarily large data sets. Originally developed for use in reduced-order modeling, the generality of the toolset allows for a range of complex methodologies to applied. This is was leveraged for use in this work for the application of methods, as the problem size is intractable on desktop machines. The toolset leverages commonly available packages in the form of ScaLAPACK, PBLAS, and PBLACS~\cite{slug} routines. ScaLAPACK is widely available on most computing resources where the original simulations are run.

\bibliographystyle{jfm}
\bibliography{ref}

\begin{thebibliography}{87}
\expandafter\ifx\csname natexlab\endcsname\relax\def\natexlab#1{#1}\fi
\def\au#1{#1} \def\ed#1{#1} \def\yr#1{#1}\def\at#1{#1}\def\jt#1{\textit{#1}}
  \def\bt#1{#1}\def\bvol#1{\textbf{#1}} \def\vol#1{#1} \def\pg#1{#1}
  \def\publ#1{#1}\def\arxiv#1{#1}\def\org#1{#1}\def\st#1{\textit{#1}}

\bibitem[Ansys(2016)]{ansys}
{\sc \au{Ansys}} \yr{2016}  \at{Ansys fluent user’s guide, release 17.0}.
  \jt{Ansys} .

\bibitem[Arbabi {\em et~al.\/}(2018)Arbabi, Korda \& Mezic]{arbabi1804data}
{\sc \au{Arbabi, H}, \au{Korda, M} \& \au{Mezic, I}} \yr{2018}  \at{A
  data-driven koopman model predictive control framework for nonlinear flows
  (2018)}.  \jt{arXiv preprint arXiv:1804.05291} .

\bibitem[Arbabi \& Mezi{\'c}(2017{\natexlab{{\em a\/}}})]{arbabi2017ergodic}
{\sc \au{Arbabi, Hassan} \& \au{Mezi{\'c}, Igor}} \yr{2017{\natexlab{{\em
  a\/}}}}  \at{Ergodic theory, dynamic mode decomposition, and computation of
  spectral properties of the koopman operator}.  \jt{SIAM Journal on Applied
  Dynamical Systems}  \bvol{16}~(4),  \pg{2096--2126}.

\bibitem[Arbabi \& Mezi{\'c}(2017{\natexlab{{\em b\/}}})]{arbabi2017study}
{\sc \au{Arbabi, Hassan} \& \au{Mezi{\'c}, Igor}} \yr{2017{\natexlab{{\em
  b\/}}}}  \at{Study of dynamics in post-transient flows using koopman mode
  decomposition}.  \jt{Physical Review Fluids}  \bvol{2}~(12),  \pg{124402}.

\bibitem[Argyriou {\em et~al.\/}(2008{\natexlab{{\em a\/}}})Argyriou, Evgeniou
  \& Pontil]{argyriou2008convex}
{\sc \au{Argyriou, Andreas}, \au{Evgeniou, Theodoros} \& \au{Pontil,
  Massimiliano}} \yr{2008{\natexlab{{\em a\/}}}}  \at{Convex multi-task feature
  learning}.  \jt{Machine Learning}  \bvol{73}~(3),  \pg{243--272}.

\bibitem[Argyriou {\em et~al.\/}(2008{\natexlab{{\em b\/}}})Argyriou, Pontil,
  Ying \& Micchelli]{argyriou2008spectral}
{\sc \au{Argyriou, Andreas}, \au{Pontil, Massimiliano}, \au{Ying, Yiming} \&
  \au{Micchelli, Charles~A}} \yr{2008{\natexlab{{\em b\/}}}} A spectral
  regularization framework for multi-task structure learning.  \bt{In {\em
  Advances in neural information processing systems\/}},  \pg{pp. 25--32}.

\bibitem[Aronszajn(1950)]{aronszajn1950theory}
{\sc \au{Aronszajn, Nachman}} \yr{1950}  \at{Theory of reproducing kernels}.
  \jt{Transactions of the American mathematical society}  \bvol{68}~(3),
  \pg{337--404}.

\bibitem[Askham \& Kutz(2018)]{askham2018variable}
{\sc \au{Askham, Travis} \& \au{Kutz, J~Nathan}} \yr{2018}  \at{Variable
  projection methods for an optimized dynamic mode decomposition}.  \jt{SIAM
  Journal on Applied Dynamical Systems}  \bvol{17}~(1),  \pg{380--416}.

\bibitem[Azencot {\em et~al.\/}(2019)Azencot, Yin \&
  Bertozzi]{azencot2019consistent}
{\sc \au{Azencot, Omri}, \au{Yin, Wotao} \& \au{Bertozzi, Andrea}} \yr{2019}
  \at{Consistent dynamic mode decomposition}.  \jt{SIAM Journal on Applied
  Dynamical Systems}  \bvol{18}~(3),  \pg{1565--1585}.

\bibitem[Bach {\em et~al.\/}(2012)Bach, Jenatton, Mairal, Obozinski {\em
  et~al.\/}]{bach2012optimization}
{\sc \au{Bach, Francis}, \au{Jenatton, Rodolphe}, \au{Mairal, Julien},
  \au{Obozinski, Guillaume} \& \au{others}} \yr{2012}  \at{Optimization with
  sparsity-inducing penalties}.  \jt{Foundations and Trends{\textregistered} in
  Machine Learning}  \bvol{4}~(1),  \pg{1--106}.

\bibitem[Bagheri(2013)]{bagheri2013koopman}
{\sc \au{Bagheri, Shervin}} \yr{2013}  \at{Koopman-mode decomposition of the
  cylinder wake}.  \jt{Journal of Fluid Mechanics}  \bvol{726},  \pg{596--623}.

\bibitem[Baudin {\em et~al.\/}(2015)]{baudin2015pydoe}
{\sc \au{Baudin, M} \& \au{others}} \yr{2015} pydoe.

\bibitem[Blackford {\em et~al.\/}(1997)Blackford, Choi, Cleary, D'Azevedo,
  Demmel, Dhillon, Dongarra, Hammarling, Henry, Petitet, Stanley, Walker \&
  Whaley]{slug}
{\sc \au{Blackford, L.~S.}, \au{Choi, J.}, \au{Cleary, A.}, \au{D'Azevedo, E.},
  \au{Demmel, J.}, \au{Dhillon, I.}, \au{Dongarra, J.}, \au{Hammarling, S.},
  \au{Henry, G.}, \au{Petitet, A.}, \au{Stanley, K.}, \au{Walker, D.} \&
  \au{Whaley, R.~C.}} \yr{1997} {\em {ScaLAPACK} Users' Guide\/}.
  \publ{Philadelphia, PA: Society for Industrial and Applied Mathematics}.

\bibitem[Bruce {\em et~al.\/}(2019)Bruce, Zeidan \&
  Bernstein]{bruce2019koopman}
{\sc \au{Bruce, Adam~L}, \au{Zeidan, Vera~M} \& \au{Bernstein, Dennis~S}}
  \yr{2019} What is the koopman operator? a simplified treatment for
  discrete-time systems.  \bt{In {\em 2019 American Control Conference
  (ACC)\/}},  \pg{pp. 1912--1917}. IEEE.

\bibitem[Brunton {\em et~al.\/}(2017)Brunton, Brunton, Proctor, Kaiser \&
  Kutz]{brunton2017chaos}
{\sc \au{Brunton, Steven~L}, \au{Brunton, Bingni~W}, \au{Proctor, Joshua~L},
  \au{Kaiser, Eurika} \& \au{Kutz, J~Nathan}} \yr{2017}  \at{Chaos as an
  intermittently forced linear system}.  \jt{Nature communications}
  \bvol{8}~(1),  \pg{19}.

\bibitem[Brunton {\em et~al.\/}(2016{\natexlab{{\em a\/}}})Brunton, Brunton,
  Proctor \& Kutz]{brunton2016koopman}
{\sc \au{Brunton, Steven~L}, \au{Brunton, Bingni~W}, \au{Proctor, Joshua~L} \&
  \au{Kutz, J~Nathan}} \yr{2016{\natexlab{{\em a\/}}}}  \at{Koopman invariant
  subspaces and finite linear representations of nonlinear dynamical systems
  for control}.  \jt{PloS one}  \bvol{11}~(2),  \pg{e0150171}.

\bibitem[Brunton {\em et~al.\/}(2016{\natexlab{{\em b\/}}})Brunton, Proctor \&
  Kutz]{brunton2016discovering}
{\sc \au{Brunton, Steven~L}, \au{Proctor, Joshua~L} \& \au{Kutz, J~Nathan}}
  \yr{2016{\natexlab{{\em b\/}}}}  \at{Discovering governing equations from
  data by sparse identification of nonlinear dynamical systems}.
  \jt{Proceedings of the National Academy of Sciences}  \bvol{113}~(15),
  \pg{3932--3937}.

\bibitem[Budi{\v{s}}i{\'c} {\em et~al.\/}(2012)Budi{\v{s}}i{\'c}, Mohr \&
  Mezi{\'c}]{budivsic2012applied}
{\sc \au{Budi{\v{s}}i{\'c}, Marko}, \au{Mohr, Ryan} \& \au{Mezi{\'c}, Igor}}
  \yr{2012}  \at{Applied koopmanism}.  \jt{Chaos: An Interdisciplinary Journal
  of Nonlinear Science}  \bvol{22}~(4),  \pg{047510}.

\bibitem[Carlberg {\em et~al.\/}(2013)Carlberg, Farhat, Cortial \&
  Amsallem]{carlberg2013gnat}
{\sc \au{Carlberg, Kevin}, \au{Farhat, Charbel}, \au{Cortial, Julien} \&
  \au{Amsallem, David}} \yr{2013}  \at{The gnat method for nonlinear model
  reduction: effective implementation and application to computational fluid
  dynamics and turbulent flows}.  \jt{Journal of Computational Physics}
  \bvol{242},  \pg{623--647}.

\bibitem[Chen {\em et~al.\/}(2012)Chen, Tu \& Rowley]{chen2012variants}
{\sc \au{Chen, Kevin~K}, \au{Tu, Jonathan~H} \& \au{Rowley, Clarence~W}}
  \yr{2012}  \at{Variants of dynamic mode decomposition: boundary condition,
  koopman, and fourier analyses}.  \jt{Journal of nonlinear science}
  \bvol{22}~(6),  \pg{887--915}.

\bibitem[Cho \& Saul(2009)]{cho2009kernel}
{\sc \au{Cho, Youngmin} \& \au{Saul, Lawrence~K}} \yr{2009} Kernel methods for
  deep learning.  \bt{In {\em Advances in neural information processing
  systems\/}},  \pg{pp. 342--350}.

\bibitem[Dalcin {\em et~al.\/}(2011)Dalcin, Paz, Kler \& Cosimo]{mpi4py}
{\sc \au{Dalcin, Lisandro~D}, \au{Paz, Rodrigo~R}, \au{Kler, Pablo~A} \&
  \au{Cosimo, Alejandro}} \yr{2011}  \at{Parallel distributed computing using
  python}.  \jt{Advances in Water Resources}  \bvol{34}~(9),  \pg{1124--1139}.

\bibitem[Dawson {\em et~al.\/}(2016)Dawson, Hemati, Williams \&
  Rowley]{dawson2016characterizing}
{\sc \au{Dawson, Scott~TM}, \au{Hemati, Maziar~S}, \au{Williams, Matthew~O} \&
  \au{Rowley, Clarence~W}} \yr{2016}  \at{Characterizing and correcting for the
  effect of sensor noise in the dynamic mode decomposition}.  \jt{Experiments
  in Fluids}  \bvol{57}~(3),  \pg{42}.

\bibitem[DeGennaro \& Urban(2019)]{degennaro2019scalable}
{\sc \au{DeGennaro, Anthony~M} \& \au{Urban, Nathan~M}} \yr{2019}  \at{Scalable
  extended dynamic mode decomposition using random kernel approximation}.
  \jt{SIAM Journal on Scientific Computing}  \bvol{41}~(3),  \pg{A1482--A1499}.

\bibitem[Dongarra {\em et~al.\/}(2018)Dongarra, Gates, Haidar, Kurzak,
  Luszczek, Tomov \& Yamazaki]{dongarra2018singular}
{\sc \au{Dongarra, Jack}, \au{Gates, Mark}, \au{Haidar, Azzam}, \au{Kurzak,
  Jakub}, \au{Luszczek, Piotr}, \au{Tomov, Stanimire} \& \au{Yamazaki,
  Ichitaro}} \yr{2018}  \at{The singular value decomposition: Anatomy of
  optimizing an algorithm for extreme scale}.  \jt{SIAM review}  \bvol{60}~(4),
   \pg{808--865}.

\bibitem[Dowell \& Hall(2001)]{dowell2001modeling}
{\sc \au{Dowell, Earl~H} \& \au{Hall, Kenneth~C}} \yr{2001}  \at{Modeling of
  fluid-structure interaction}.  \jt{Annual review of fluid mechanics}
  \bvol{33}~(1),  \pg{445--490}.

\bibitem[Drazin \& Reid(2004)]{drazin2004hydrodynamic}
{\sc \au{Drazin, Philip~G} \& \au{Reid, William~Hill}} \yr{2004} {\em
  Hydrodynamic stability\/}.  \publ{Cambridge university press}.

\bibitem[Fey {\em et~al.\/}(1998)Fey, K{\"o}nig \& Eckelmann]{fey1998new}
{\sc \au{Fey, Uwe}, \au{K{\"o}nig, Michael} \& \au{Eckelmann, Helmut}}
  \yr{1998}  \at{A new strouhal--reynolds-number relationship for the circular
  cylinder in the range 47< re< 2$\times$ 10 5}.  \jt{Physics of Fluids}
  \bvol{10}~(7),  \pg{1547--1549}.

\bibitem[Forrest \& Owen(2010)]{forrest2010investigation}
{\sc \au{Forrest, James~S} \& \au{Owen, Ieuan}} \yr{2010}  \at{An investigation
  of ship airwakes using detached-eddy simulation}.  \jt{Computers \& Fluids}
  \bvol{39}~(4),  \pg{656--673}.

\bibitem[Haseli \& Cort{\'e}s(2019)]{haseli2019efficient}
{\sc \au{Haseli, Masih} \& \au{Cort{\'e}s, Jorge}} \yr{2019}  \at{Efficient
  identification of linear evolutions in nonlinear vector fields: Koopman
  invariant subspaces}.  \jt{arXiv preprint arXiv:1909.01419} .

\bibitem[Holmes {\em et~al.\/}(2012)Holmes, Lumley, Berkooz \&
  Rowley]{holmes2012turbulence}
{\sc \au{Holmes, Philip}, \au{Lumley, John~L}, \au{Berkooz, Gahl} \&
  \au{Rowley, Clarence~W}} \yr{2012} {\em Turbulence, coherent structures,
  dynamical systems and symmetry\/}.  \publ{Cambridge university press}.

\bibitem[Hua {\em et~al.\/}(2017)Hua, Noorian, Moss, Leong \&
  Gunaratne]{hua2017high}
{\sc \au{Hua, Jia-Chen}, \au{Noorian, Farzad}, \au{Moss, Duncan}, \au{Leong,
  Philip~HW} \& \au{Gunaratne, Gemunu~H}} \yr{2017}  \at{High-dimensional time
  series prediction using kernel-based koopman mode regression}.  \jt{Nonlinear
  Dynamics}  \bvol{90}~(3),  \pg{1785--1806}.

\bibitem[Huang {\em et~al.\/}(2018)Huang, Duraisamy \&
  Merkle]{huang2018challenges}
{\sc \au{Huang, Cheng}, \au{Duraisamy, Karthik} \& \au{Merkle, Charles}}
  \yr{2018} Challenges in reduced order modeling of reacting flows.  \bt{In
  {\em 2018 Joint Propulsion Conference\/}},  \pg{p. 4675}.

\bibitem[Jasak {\em et~al.\/}(2007)Jasak, Jemcov, Tukovic {\em
  et~al.\/}]{jasak2007openfoam}
{\sc \au{Jasak, Hrvoje}, \au{Jemcov, Aleksandar}, \au{Tukovic, Zeljko} \&
  \au{others}} \yr{2007} Openfoam: A c++ library for complex physics
  simulations.  \bt{In {\em International workshop on coupled methods in
  numerical dynamics\/}}, ,  \vol{vol. 1000},  \pg{pp. 1--20}. IUC Dubrovnik
  Croatia.

\bibitem[Jones {\em et~al.\/}(2001--)Jones, Oliphant, Peterson {\em
  et~al.\/}]{scipy}
{\sc \au{Jones, Eric}, \au{Oliphant, Travis}, \au{Peterson, Pearu} \&
  \au{others}} \yr{2001--} {SciPy}: Open source scientific tools for {Python}.

\bibitem[Jovanovi{\'c} {\em et~al.\/}(2014)Jovanovi{\'c}, Schmid \&
  Nichols]{jovanovic2014sparsity}
{\sc \au{Jovanovi{\'c}, Mihailo~R}, \au{Schmid, Peter~J} \& \au{Nichols,
  Joseph~W}} \yr{2014}  \at{Sparsity-promoting dynamic mode decomposition}.
  \jt{Physics of Fluids}  \bvol{26}~(2),  \pg{024103}.

\bibitem[Jungers \& Tabuada(2019)]{jungers2019non}
{\sc \au{Jungers, Rapha{\"e}l~M} \& \au{Tabuada, Paulo}} \yr{2019} Non-local
  linearization of nonlinear differential equations via polyflows.  \bt{In {\em
  2019 American Control Conference (ACC)\/}},  \pg{pp. 1--6}. IEEE.

\bibitem[Kaiser {\em et~al.\/}(2017)Kaiser, Kutz \& Brunton]{kaiser2017data}
{\sc \au{Kaiser, Eurika}, \au{Kutz, J~Nathan} \& \au{Brunton, Steven~L}}
  \yr{2017}  \at{Data-driven discovery of koopman eigenfunctions for control}.
  \jt{arXiv preprint arXiv:1707.01146} .

\bibitem[Kamb {\em et~al.\/}(2018)Kamb, Kaiser, Brunton \& Kutz]{kamb2018time}
{\sc \au{Kamb, Mason}, \au{Kaiser, Eurika}, \au{Brunton, Steven~L} \& \au{Kutz,
  J~Nathan}} \yr{2018}  \at{Time-delay observables for koopman: Theory and
  applications}.  \jt{arXiv preprint arXiv:1810.01479} .

\bibitem[Koopman(1931)]{koopman1931hamiltonian}
{\sc \au{Koopman, Bernard~O}} \yr{1931}  \at{Hamiltonian systems and
  transformation in hilbert space}.  \jt{Proceedings of the National Academy of
  Sciences of the United States of America}  \bvol{17}~(5),  \pg{315}.

\bibitem[Korda \& Mezi{\'c}(2018{\natexlab{{\em a\/}}})]{korda2018linear}
{\sc \au{Korda, Milan} \& \au{Mezi{\'c}, Igor}} \yr{2018{\natexlab{{\em a\/}}}}
   \at{Linear predictors for nonlinear dynamical systems: Koopman operator
  meets model predictive control}.  \jt{Automatica}  \bvol{93},  \pg{149--160}.

\bibitem[Korda \& Mezi{\'c}(2018{\natexlab{{\em b\/}}})]{korda2018convergence}
{\sc \au{Korda, Milan} \& \au{Mezi{\'c}, Igor}} \yr{2018{\natexlab{{\em b\/}}}}
   \at{On convergence of extended dynamic mode decomposition to the koopman
  operator}.  \jt{Journal of Nonlinear Science}  \bvol{28}~(2),  \pg{687--710}.

\bibitem[Kou \& Zhang(2017)]{kou2017improved}
{\sc \au{Kou, Jiaqing} \& \au{Zhang, Weiwei}} \yr{2017}  \at{An improved
  criterion to select dominant modes from dynamic mode decomposition}.
  \jt{European Journal of Mechanics-B/Fluids}  \bvol{62},  \pg{109--129}.

\bibitem[Kutz {\em et~al.\/}(2016)Kutz, Brunton, Brunton \&
  Proctor]{kutz2016dynamic}
{\sc \au{Kutz, J~Nathan}, \au{Brunton, Steven~L}, \au{Brunton, Bingni~W} \&
  \au{Proctor, Joshua~L}} \yr{2016} {\em Dynamic mode decomposition:
  data-driven modeling of complex systems\/}.  \publ{SIAM}.

\bibitem[Le~Clainche \& Vega(2017)]{le2017higher}
{\sc \au{Le~Clainche, Soledad} \& \au{Vega, Jos{\'e}~M}} \yr{2017}  \at{Higher
  order dynamic mode decomposition}.  \jt{SIAM Journal on Applied Dynamical
  Systems}  \bvol{16}~(2),  \pg{882--925}.

\bibitem[Li {\em et~al.\/}(2017)Li, Dietrich, Bollt \&
  Kevrekidis]{li2017extended}
{\sc \au{Li, Qianxiao}, \au{Dietrich, Felix}, \au{Bollt, Erik~M} \&
  \au{Kevrekidis, Ioannis~G}} \yr{2017}  \at{Extended dynamic mode
  decomposition with dictionary learning: A data-driven adaptive spectral
  decomposition of the koopman operator}.  \jt{Chaos: An Interdisciplinary
  Journal of Nonlinear Science}  \bvol{27}~(10),  \pg{103111}.

\bibitem[Lietz {\em et~al.\/}(2017)Lietz, Johnsen \& Kushner]{lietz2017plasma}
{\sc \au{Lietz, Amanda~M}, \au{Johnsen, Eric} \& \au{Kushner, Mark~J}}
  \yr{2017}  \at{Plasma-induced flow instabilities in atmospheric pressure
  plasma jets}.  \jt{Applied Physics Letters}  \bvol{111}~(11),  \pg{114101}.

\bibitem[Lusch {\em et~al.\/}(2018)Lusch, Kutz \& Brunton]{lusch2018deep}
{\sc \au{Lusch, Bethany}, \au{Kutz, J~Nathan} \& \au{Brunton, Steven~L}}
  \yr{2018}  \at{Deep learning for universal linear embeddings of nonlinear
  dynamics}.  \jt{Nature communications}  \bvol{9}~(1),  \pg{4950}.

\bibitem[Mamakoukas {\em et~al.\/}(2019)Mamakoukas, Castano, Tan \&
  Murphey]{mamakoukas2019local}
{\sc \au{Mamakoukas, Giorgos}, \au{Castano, Maria~L}, \au{Tan, Xiaobo} \&
  \au{Murphey, Todd}} \yr{2019} Local koopman operators for data-driven control
  of robotic systems.  \bt{In {\em Robotics: science and systems\/}}.

\bibitem[Mercer(1909)]{mercer1909functions}
{\sc \au{Mercer, J}} \yr{1909} Functions of positive and negative type and
  their connection with the theory of integral equations, philosophical
  transsaction of the royal society of london, ser.

\bibitem[Mezi{\'c}(2005)]{mezic2005spectral}
{\sc \au{Mezi{\'c}, Igor}} \yr{2005}  \at{Spectral properties of dynamical
  systems, model reduction and decompositions}.  \jt{Nonlinear Dynamics}
  \bvol{41}~(1-3),  \pg{309--325}.

\bibitem[Mezi{\'c}(2013)]{mezic2013analysis}
{\sc \au{Mezi{\'c}, Igor}} \yr{2013}  \at{Analysis of fluid flows via spectral
  properties of the koopman operator}.  \jt{Annual Review of Fluid Mechanics}
  \bvol{45},  \pg{357--378}.

\bibitem[Minh(2010)]{minh2010some}
{\sc \au{Minh, Ha~Quang}} \yr{2010}  \at{Some properties of gaussian
  reproducing kernel hilbert spaces and their implications for function
  approximation and learning theory}.  \jt{Constructive Approximation}
  \bvol{32}~(2),  \pg{307--338}.

\bibitem[Noack {\em et~al.\/}(2003)Noack, Afanasiev, MORZY{\'N}SKI, Tadmor \&
  Thiele]{noack2003hierarchy}
{\sc \au{Noack, Bernd~R}, \au{Afanasiev, Konstantin}, \au{MORZY{\'N}SKI,
  MAREK}, \au{Tadmor, Gilead} \& \au{Thiele, Frank}} \yr{2003}  \at{A hierarchy
  of low-dimensional models for the transient and post-transient cylinder
  wake}.  \jt{Journal of Fluid Mechanics}  \bvol{497},  \pg{335--363}.

\bibitem[Otto \& Rowley(2019)]{otto2019linearly}
{\sc \au{Otto, Samuel~E} \& \au{Rowley, Clarence~W}} \yr{2019}  \at{Linearly
  recurrent autoencoder networks for learning dynamics}.  \jt{SIAM Journal on
  Applied Dynamical Systems}  \bvol{18}~(1),  \pg{558--593}.

\bibitem[Pan \& Duraisamy(2018)]{pan2018long}
{\sc \au{Pan, Shaowu} \& \au{Duraisamy, Karthik}} \yr{2018}  \at{Long-time
  predictive modeling of nonlinear dynamical systems using neural networks}.
  \jt{Complexity}  \bvol{2018}.

\bibitem[Pan \& Duraisamy(2019)]{pan2019structure}
{\sc \au{Pan, Shaowu} \& \au{Duraisamy, Karthik}} \yr{2019}  \at{On the
  structure of time-delay embedding in linear models of non-linear dynamical
  systems}.  \jt{arXiv preprint arXiv:1902.05198} .

\bibitem[Pan \& Duraisamy(2020)]{pan2019physics}
{\sc \au{Pan, Shaowu} \& \au{Duraisamy, Karthik}} \yr{2020}
  \at{Physics-informed probabilistic learning of linear embeddings of nonlinear
  dynamics with guaranteed stability}.  \jt{SIAM Journal on Applied Dynamical
  Systems}  \bvol{19}~(1),  \pg{480--509}.

\bibitem[Parikh {\em et~al.\/}(2014)Parikh, Boyd {\em
  et~al.\/}]{parikh2014proximal}
{\sc \au{Parikh, Neal}, \au{Boyd, Stephen} \& \au{others}} \yr{2014}
  \at{Proximal algorithms}.  \jt{Foundations and Trends{\textregistered} in
  Optimization}  \bvol{1}~(3),  \pg{127--239}.

\bibitem[Parish {\em et~al.\/}(2020)Parish, Wentland \&
  Duraisamy]{parish2020adjoint}
{\sc \au{Parish, Eric~J}, \au{Wentland, Christopher~R} \& \au{Duraisamy,
  Karthik}} \yr{2020}  \at{The adjoint petrov--galerkin method for non-linear
  model reduction}.  \jt{Computer Methods in Applied Mechanics and Engineering}
   \bvol{365},  \pg{112991}.

\bibitem[Parry(2004)]{parry2004topics}
{\sc \au{Parry, William}} \yr{2004} {\em Topics in ergodic theory\/}, ,
  \vol{vol.~75}.  \publ{Cambridge university press}.

\bibitem[Pedregosa {\em et~al.\/}(2011)Pedregosa, Varoquaux, Gramfort, Michel,
  Thirion, Grisel, Blondel, Prettenhofer, Weiss, Dubourg, Vanderplas, Passos,
  Cournapeau, Brucher, Perrot \& Duchesnay]{scikit-learn}
{\sc \au{Pedregosa, F.}, \au{Varoquaux, G.}, \au{Gramfort, A.}, \au{Michel,
  V.}, \au{Thirion, B.}, \au{Grisel, O.}, \au{Blondel, M.}, \au{Prettenhofer,
  P.}, \au{Weiss, R.}, \au{Dubourg, V.}, \au{Vanderplas, J.}, \au{Passos, A.},
  \au{Cournapeau, D.}, \au{Brucher, M.}, \au{Perrot, M.} \& \au{Duchesnay, E.}}
  \yr{2011}  \at{Scikit-learn: Machine learning in {P}ython}.  \jt{Journal of
  Machine Learning Research}  \bvol{12},  \pg{2825--2830}.

\bibitem[Pope(2001)]{pope2001turbulent}
{\sc \au{Pope, Stephen~B}} \yr{2001} Turbulent flows.

\bibitem[Rasmussen(2003)]{rasmussen2003gaussian}
{\sc \au{Rasmussen, Carl~Edward}} \yr{2003} Gaussian processes in machine
  learning.  \bt{In {\em Summer School on Machine Learning\/}},  \pg{pp.
  63--71}. Springer.

\bibitem[R{\"o}jsel(2017)]{rojsel2017koopman}
{\sc \au{R{\"o}jsel, Jimmy}} \yr{2017}  \bt{Koopman mode analysis of the
  side-by-side cylinder wake}. {\em Tech. Rep.\/}.  \org{Royal Institute of
  Technology, Department of Mechanics}.

\bibitem[Rowley \& Dawson(2017)]{rowley2017model}
{\sc \au{Rowley, Clarence~W} \& \au{Dawson, Scott~TM}} \yr{2017}  \at{Model
  reduction for flow analysis and control}.  \jt{Annual Review of Fluid
  Mechanics}  \bvol{49},  \pg{387--417}.

\bibitem[Rowley {\em et~al.\/}(2009)Rowley, Mezi{\'c}, Bagheri, Schlatter \&
  Henningson]{rowley2009spectral}
{\sc \au{Rowley, Clarence~W}, \au{Mezi{\'c}, Igor}, \au{Bagheri, Shervin},
  \au{Schlatter, Philipp} \& \au{Henningson, Dan~S}} \yr{2009}  \at{Spectral
  analysis of nonlinear flows}.  \jt{Journal of fluid mechanics}  \bvol{641},
  \pg{115--127}.

\bibitem[Schmid(2007)]{schmid2007nonmodal}
{\sc \au{Schmid, Peter~J}} \yr{2007}  \at{Nonmodal stability theory}.
  \jt{Annu. Rev. Fluid Mech.}  \bvol{39},  \pg{129--162}.

\bibitem[Schmid(2010)]{schmid2010dynamic}
{\sc \au{Schmid, Peter~J}} \yr{2010}  \at{Dynamic mode decomposition of
  numerical and experimental data}.  \jt{Journal of fluid mechanics}
  \bvol{656},  \pg{5--28}.

\bibitem[Schmid {\em et~al.\/}(2012)Schmid, Violato \&
  Scarano]{schmid2012decomposition}
{\sc \au{Schmid, Peter~J}, \au{Violato, Daniele} \& \au{Scarano, Fulvio}}
  \yr{2012}  \at{Decomposition of time-resolved tomographic piv}.
  \jt{Experiments in Fluids}  \bvol{52}~(6),  \pg{1567--1579}.

\bibitem[Sharma {\em et~al.\/}(2019)Sharma, Xu, Padthe, Friedmann \&
  Duraisamy]{sharma2019simulation}
{\sc \au{Sharma, Abhinav}, \au{Xu, Jiayang}, \au{Padthe, Ashwani~K},
  \au{Friedmann, Peretz~P} \& \au{Duraisamy, Karthik}} \yr{2019} Simulation of
  maritime helicopter dynamics during approach to landing with time-accurate
  wind-over-deck.  \bt{In {\em AIAA Scitech 2019 Forum\/}},  \pg{p. 0861}.

\bibitem[Taira {\em et~al.\/}(2019)Taira, Hemati, Brunton, Sun, Duraisamy,
  Bagheri, Dawson \& Yeh]{taira2019modal}
{\sc \au{Taira, Kunihiko}, \au{Hemati, Maziar~S}, \au{Brunton, Steven~L},
  \au{Sun, Yiyang}, \au{Duraisamy, Karthik}, \au{Bagheri, Shervin}, \au{Dawson,
  Scott~TM} \& \au{Yeh, Chi-An}} \yr{2019}  \at{Modal analysis of fluid flows:
  Applications and outlook}.  \jt{AIAA Journal}  \pg{pp. 1--25}.

\bibitem[Takeishi {\em et~al.\/}(2017)Takeishi, Kawahara \&
  Yairi]{takeishi2017learning}
{\sc \au{Takeishi, Naoya}, \au{Kawahara, Yoshinobu} \& \au{Yairi, Takehisa}}
  \yr{2017} Learning koopman invariant subspaces for dynamic mode
  decomposition.  \bt{In {\em Advances in Neural Information Processing
  Systems\/}},  \pg{pp. 1130--1140}.

\bibitem[Tissot {\em et~al.\/}(2014)Tissot, Cordier, Benard \&
  Noack]{tissot2014model}
{\sc \au{Tissot, Gilles}, \au{Cordier, Laurent}, \au{Benard, Nicolas} \&
  \au{Noack, Bernd~R}} \yr{2014}  \at{Model reduction using dynamic mode
  decomposition}.  \jt{Comptes Rendus M{\'e}canique}  \bvol{342}~(6-7),
  \pg{410--416}.

\bibitem[Tu {\em et~al.\/}(2013)Tu, Rowley, Luchtenburg, Brunton \&
  Kutz]{tu2013dynamic}
{\sc \au{Tu, Jonathan~H}, \au{Rowley, Clarence~W}, \au{Luchtenburg, Dirk~M},
  \au{Brunton, Steven~L} \& \au{Kutz, J~Nathan}} \yr{2013}  \at{On dynamic mode
  decomposition: theory and applications}.  \jt{arXiv preprint arXiv:1312.0041}
  .

\bibitem[Wang {\em et~al.\/}(2019)Wang, Wang, Xiang \& hui Yue]{wang2019re}
{\sc \au{Wang, Linyu}, \au{Wang, Junyan}, \au{Xiang, Jianhong} \& \au{hui Yue,
  Hui}} \yr{2019}  \at{A re-weighted smoothed l 0-norm regularized sparse
  reconstructed algorithm for linear inverse problems}.  \jt{Journal of Physics
  Communications} .

\bibitem[White \& Corfield(2006)]{white2006viscous}
{\sc \au{White, Frank~M} \& \au{Corfield, Isla}} \yr{2006} {\em Viscous fluid
  flow\/}, ,  \vol{vol.~3}.  \publ{McGraw-Hill New York}.

\bibitem[Williams {\em et~al.\/}(2015)Williams, Kevrekidis \&
  Rowley]{williams2015data}
{\sc \au{Williams, Matthew~O}, \au{Kevrekidis, Ioannis~G} \& \au{Rowley,
  Clarence~W}} \yr{2015}  \at{A data--driven approximation of the koopman
  operator: Extending dynamic mode decomposition}.  \jt{Journal of Nonlinear
  Science}  \bvol{25}~(6),  \pg{1307--1346}.

\bibitem[Williams {\em et~al.\/}(2014)Williams, Rowley \&
  Kevrekidis]{williams2014kernel}
{\sc \au{Williams, Matthew~O}, \au{Rowley, Clarence~W} \& \au{Kevrekidis,
  Ioannis~G}} \yr{2014}  \at{A kernel-based approach to data-driven koopman
  spectral analysis}.  \jt{arXiv:1411.2260} .

\bibitem[Xu \& Duraisamy(2019)]{xu2019multi}
{\sc \au{Xu, Jiayang} \& \au{Duraisamy, Karthik}} \yr{2019}  \at{Multi-level
  convolutional autoencoder networks for parametric prediction of
  spatio-temporal dynamics}.  \jt{arXiv preprint arXiv:1912.11114} .

\bibitem[Xu {\em et~al.\/}(2020)Xu, Huang \& Duraisamy]{xu2020reduced}
{\sc \au{Xu, Jiayang}, \au{Huang, Cheng} \& \au{Duraisamy, Karthik}} \yr{2020}
  \at{Reduced-order modeling framework for combustor instabilities using
  truncated domain training}.  \jt{AIAA Journal}  \bvol{58}~(2),
  \pg{618--632}.

\bibitem[Yeung {\em et~al.\/}(2019)Yeung, Kundu \& Hodas]{yeung2019learning}
{\sc \au{Yeung, Enoch}, \au{Kundu, Soumya} \& \au{Hodas, Nathan}} \yr{2019}
  Learning deep neural network representations for koopman operators of
  nonlinear dynamical systems.  \bt{In {\em 2019 American Control Conference
  (ACC)\/}},  \pg{pp. 4832--4839}. IEEE.

\bibitem[Yuan {\em et~al.\/}(2018)Yuan, Wall \& Lee]{yuan2018combined}
{\sc \au{Yuan, Weixing}, \au{Wall, Alanna} \& \au{Lee, Richard}} \yr{2018}
  \at{Combined numerical and experimental simulations of unsteady ship
  airwakes}.  \jt{Computers \& Fluids}  \bvol{172},  \pg{29--53}.

\bibitem[Zhang {\em et~al.\/}(2017)Zhang, Dawson, Rowley, Deem \&
  Cattafesta]{zhang2017evaluating}
{\sc \au{Zhang, Hao}, \au{Dawson, Scott}, \au{Rowley, Clarence~W}, \au{Deem,
  Eric~A} \& \au{Cattafesta, Louis~N}} \yr{2017}  \at{Evaluating the accuracy
  of the dynamic mode decomposition}.  \jt{arXiv preprint arXiv:1710.00745} .

\bibitem[Zhang \& Schaeffer(2019)]{zhang2019convergence}
{\sc \au{Zhang, Linan} \& \au{Schaeffer, Hayden}} \yr{2019}  \at{On the
  convergence of the sindy algorithm}.  \jt{Multiscale Modeling \& Simulation}
  \bvol{17}~(3),  \pg{948--972}.

\bibitem[Zhao {\em et~al.\/}(2015)Zhao, Sun, Ye, Chen, Lu \&
  Ramakrishnan]{zhao2015multi}
{\sc \au{Zhao, Liang}, \au{Sun, Qian}, \au{Ye, Jieping}, \au{Chen, Feng},
  \au{Lu, Chang-Tien} \& \au{Ramakrishnan, Naren}} \yr{2015} Multi-task
  learning for spatio-temporal event forecasting.  \bt{In {\em Proceedings of
  the 21th ACM SIGKDD International Conference on Knowledge Discovery and Data
  Mining\/}},  \pg{pp. 1503--1512}. ACM.

\bibitem[Zou \& Hastie(2005)]{zou2005regularization}
{\sc \au{Zou, Hui} \& \au{Hastie, Trevor}} \yr{2005}  \at{Regularization and
  variable selection via the elastic net}.  \jt{Journal of the royal
  statistical society: series B (statistical methodology)}  \bvol{67}~(2),
  \pg{301--320}.

\end{thebibliography}

\end{document}